\documentclass[10pt]{article}
\usepackage{graphicx}
\usepackage{amsmath,amsthm,amssymb}
\newtheorem{remark}{Remark}
\usepackage{booktabs}
\usepackage{hyperref}

\usepackage{color}

\usepackage{geometry}
\usepackage{fancyhdr}
\usepackage[ruled,linesnumbered]{algorithm2e}

\begin{document}
\title{Model order reduction for parameterized electromagnetic 
       problems using matrix decomposition and deep neural networks}
       
\date{}

\author{Xiao-Feng~He$^{a}$, \mbox{  } Liang~Li$^{a}$,\mbox{  } St\'ephane~Lanteri$^{b}$\mbox{  } Kun~Li$^{c}$\\
{\small{\it a. School of Mathematical Sciences,}}\\
{\small{\it University of Electronic Science and Technology of China, Chengdu, Sichuan, P. R. China}}\\
{\small{\it b. Universit\'e C\^ote d'Azur, Inria, CNRS, LJAD, Inria, }}\\
{\small{\it 2004 Route des Lucioles, BP 93, 06902, Sophia Antipolis Cedex, France}} \\ 
{\small{\it c. School of Economic Mathematics, }}\\
{\small{\it Southwestern University of Finance and Economics, Chengdu, P.R. China}}}

\maketitle

\begin{abstract}
A  non-intrusive  model  order  reduction  (MOR)  method  for  solving
parameterized electromagnetic scattering problems  is proposed in this
paper.  A database collecting  snapshots of high-fidelity solutions is
built by  solving the parameterized time-domain  Maxwell equations for
some values of  the material parameters using a  fullwave solver based
on a high order discontinuous  Galerkin time-domain (DGTD) method.  To
perform a prior dimensionality reduction,  a set of reduced basis (RB)
functions  are  extracted from  the  database  via a  two-step  proper
orthogonal decomposition (POD) method.  Projection coefficients of the
reduced basis functions are further compressed through a convolutional
autoencoder (CAE) network. Singular  value decomposition (SVD) is then
used to extract the principal components of the reduced-order matrices
generated  by  CAE,  and  a  cubic  spline  interpolation-based  (CSI)
approach  is  employed  for  approximating the  dominating  time-  and
parameter-modes of  the reduced-order matrices. The  generation of the
reduced basis  and the training of  the CAE and CSI are accomplished in
the  offline stage,  thus  the RB solution  for given  time/parameter
values can be quickly recovered via outputs of the interpolation model
and decoder network.  In particular,  the offline and online stages of
the proposed  RB method  are completely  decoupled, which  ensures the
validity of the  method.  The performance of the  proposed CAE-CSI ROM
is illustrated  with numerical experiments  for scattering of  a plane
wave by a 2-D dielectric disk and a multi-layer heterogeneous medium.

\end{abstract}

\section{Introduction}
\label{sec:intro}

Electromagnetic wave interaction with heterogeneous media is generally
modeled by  the system of parameterized  time-domain Maxwell equations
\cite{Jin2011}.   The  Finite  Difference  Time-Domain  (FDTD)  method
\cite{Yee1966,  TafloveHagness2005, Gedney2011}  is  the most  popular
method   for  simulating   these   time-domain  electromagnetic   wave
propagation problems.   The discontinuous Galerkin  time-domain (DGTD)
method is an alternative approach that  has emerged during the last 20
years  \cite{Hesthaven2007,  Fezoui2005}.   Compared with  the  FDTD
method,  the DGTD  method has  several attractive  advantages such  as
local approximation  strategy, easy  adaption to complex  geometry and
material composition \cite{Hesthaven2007}.  With  the DGTD method, the
number of  degree of  freedom (DOFs) is  determined by  the underlying
mesh  and the  polynomial order  used for  approximating the  EM field
components. This number  is usually high to guarantee  the accuracy of
the solution  of Maxwell's  equations.  Therefore,  it makes  sense to
study the  model order  reduction method  to reduce  the computational
burden of  CPU time and  memory, when repeatedly  simulating Maxwell's
equations with different parameters.
Model  order  reduction  (MOR)  techniques  \cite{Benner2017}  aim  to
replace  a full-order  model  (FOM) by  a  reduced-order model  (ROM),
featuring a much lower  dimension, thereby reducing computational cost
under acceptable  accuracy conditions. A reduced  basis (RB) method
based  on a  offline-online procedure  is  a widely  used model  order
reduction technique  \cite{Hesthaven2016, Guo2019,  Li2021, Wang2019}.
During the  offline stage, a reduced  space spanned by a  set of time-
and  parameter-independent  RB  functions   is  constructed  from  the
full-order solutions (snapshots), which  are for instance generated by
the  DGTD  method at  different  parameter  values. The  reduced-order
solution for  a new  time and  parameter value can  be expressed  as a
linear  combination   of  the  RB  functions   where  the  combination
coefficients, also  called the intrinsic coordinates,  need to be
calculated during the online stage.   The high-fidelity solver is only
used  to   generate  the  snapshots,  thus   guaranteeing  a  complete
decoupling   between   the   online   evaluation   and   the   offline
training \cite{Xiao2015, Li2021}.

More recently,  non-intrusive ROM combining RB  techniques and machine
learning have been developed  for solving large-scale complex physical
problems. The advantage of these  methods is that they are data-driven
and do not  require access to the governing equations  of the original
FOM. For  instance, a non-intrusive  RB method combining POD  and feed
forward  neural networks  (FNNs) has  been proposed  for parameterized
unsteady flows \cite{Hesthaven2018, Wang2019},  where a neural network
is  trained  to approximate  the  mapping  between the  time/parameter
values   and  the   projection   coefficients.   Moreover,   different
regression  and   interpolation  methods  such  as   Gaussian  process
regression (GPR)  \cite{Guo2018, Guo2019, Zhao2021,  Kast2020}, radial
basis  functions (RBF)  \cite{Xiao2017,  Li2019RBF}  and cubic  spline
interpolation  (CSI) \cite{Li2021,  Bui2003, Oulghelou2021}  have also
been considered in  place of FNNs to approximate the  mapping from the
time/parameter to the reduced coefficients.  By combining an ANN-based
regression  model  with  a   physics-informed  neural  network  (PINN)
\cite{Raissi2019,   Lu2021},  a   hybrid   strategy   is  devised   in
\cite{Chen2021} to train a network  by minimizing the residual loss of
the reduced-order equation at a set of points in the parameter space.

However,  linear reduction  methods  such as  POD  hardly capture  the
complex  dynamics of  highly  nonlinear  systems, therefore  nonlinear
manifold  learning  methods such  as  Kernel  PCA \cite{Zhou2020}  and
Hessian feature  maps \cite{Ye2015}  have attracted much  attention in
the  past several  years.   Assuming that  data points  lie  in a  low
dimensional  manifold  embedded  in an  higher  dimensional  Euclidean
space,  manifold  learning  \cite{Melas2020}   aims  to  identify  the
intrinsic  dimensionality,  equal to  the  number  of parameters  that
describe the  system, and thus obtain  low dimensional representations
of the data points.  Although kernel  PCA and Hessian feature maps are
effective  in  providing  low  dimensional  representations  for  high
dimensional  data points,  their main  drawback  is that  they do  not
provide an  analytical formula to  decode the compressed data  back to
their   high  dimensional   representation  in   the  original   space
\cite{Nikolopoulos2022}.    An   autoencoder   (AE)   overcomes   this
disadvantage by learning  how to compress (encode)  a high dimensional
data to a low dimensional code  and then reconstruct (decode) the code
to  a  representation as  close  to  the  original input  as  possible
\cite{Nikolopoulos2022}.  The   encoder  and   decoder  parts   of  an
autoencoder  are trained  simultaneously but  can be  used separately,
which provides  an opportunity  to build a  mapping between  the input
time/parameter and encoded representation.

Directly  applying  large-scale  simulation   data  (snapshots)  to  a
fully-connected autoencoder  is not only  computationally prohibitive,
but also ignores the opportunity  to exploit the structure of features
in  high  dimensional  data \cite{Gonzalez2018,  Fresca2022}.   As  an
extension of ordinary  autoencoders, convolutional autoencoders (CAEs)
\cite{Lecun1998, Guo2016}  are characterized by shared  parameters and
local  connectivity  which  help  to  reduce the  memory  as  well  as
computational costs.   For example, projection-based ROMs  proposed in
\cite{Lee2020}  project dynamic  systems onto  nonlinear manifolds  by
means of autoencoders, which still require the assembling and solution
of a ROM as in traditional POD-Galerkin ROMs.  In \cite{Gonzalez2018},
a  reduced  trial  manifold  is generated  via  a  deep  convolutional
recurrent AE,  which is then  used to  train a long  short-term memory
(LSTM) network that models the reduced dynamics. A POD-DL-ROM approach
proposed in  \cite{Fresca2022} combines and improves  the previous two
methods as the nonlinear trial manifold is learnt by using the decoder
function  of a  CAE  while the  dynamics on  the  reduced manifold  is
modeled  through  a DFNN  and  the  encoder  function  of a  CAE.   In
\cite{Nikolopoulos2022}, the authors utilize a CAE in conjunction with
a FNN  to establish a mapping  from the problem's parametric  space to
its  solution space.   The  differences  between \cite{Fresca2022}  and
\cite{Nikolopoulos2022} are  that the  former uses  the POD  method for
pre-dimensionality reduction  to reduce  training time and  trains CAE
and FNN at the same time.

In  this paper,  a  non-intrusive reduced-order  modeling strategy  is
proposed  to solve  electromagnetic  scattering  problems governed  by
parameterized  time-domain Maxwell  equations.  Firstly,  the approach
performs a  prior dimensionality  reduction by  a two-step  POD method
\cite{Li2021}.  Then, relying on its powerful dimensionality reduction
properties,  we  use  a  CAE   to  encode  and  decode  the  intrinsic coordinates.  Furthermore, a CSI-based model is devised to approximate
a mapping  from the problem's  parameter space to its  low dimensional
encoded vector space.  Using this approach, the encoded representation
of the solution at a new  time/parameter value is recovered by the CSI
model, while the RB solution in the original high dimensional space is
obtained by the decoder and  projection matrix.  The resulting CAE-CSI
ROM provides a very fast and  accurate evaluation of the entire system
response.

The paper is organized as follows.  In section \ref{sec:maxw_dgtd}, we
briefly  introduce  the time-domain  Maxwell  equations  and the  DGTD
scheme to generate the snapshots. Section \ref{sec:2step_pod} provides
a  description  for  the  two-step  POD, from  which  a  RB  basis  is
constructed.   In  section \ref{sec:cae},  we  present  how a  CAE  is
employed to further reduce the dimension of intrinsic coordinates. In section \ref{sec:csi}, we
show  how to  use  CSI to  build an  approximate  mapping between  the
problem’s  parameters and  the low  dimensional encoded  vectors.  The
overall  CAE-CSI  ROM  for  parameterized  electromagnetic  scattering
problems   is   also   presented   in  this   section.    In   Section
\ref{sec:results}, two  numerical experiments  for testing  the method
are   provided.   And   conclusion  remarks   are  drawn   in  section
\ref{sec:concl}.
\section{Full-order model and DGTD method}
\label{sec:maxw_dgtd}

Unsteady  electromagnetic  scattering  problems are  governed  by  the
following normalized form of the time-domain Maxwell equations
\begin{equation}
\begin{cases}
\nu_{r} \dfrac{\partial \mathbf{H}(\mathbf{x}, t)}{\partial t} + 
\operatorname{curl}(\mathbf{E}(\mathbf{x}, t)) = 0, & 
\forall(\mathbf{x}, t) \in \Omega \times \mathcal{T}, \\
\varepsilon_{r} \dfrac{\partial \mathbf{E}(\mathbf{x}, t)}{\partial t} -
\operatorname{curl}(\mathbf{H}(\mathbf{x}, t)) = 0, &
\forall(\mathbf{x}, t) \in \Omega \times \mathcal{T}, \\
\mathcal{L}(\mathbf{E}(\mathbf{x}, t), \mathbf{H}(\mathbf{x}, t)) = 
\mathcal{L}\left(\mathbf{E}^{\text {inc }}(\mathbf{x}, t), 
\mathbf{H}^{\text {inc }}(\mathbf{x}, t)\right), & 
\forall(\mathbf{x}, t) \in \partial \Omega \times \mathcal{T}, \\
\mathbf{E}(\mathbf{x}, 0)=\mathbf{E}_{0}(\mathbf{x}), 
\mathbf{H}(\mathbf{x}, 0)=\mathbf{H}_{0}(\mathbf{x}), &
\forall \mathbf{x} \in \Omega,
\end{cases}
\label{eq:maxw}
\end{equation}
where $\Omega$ is the spatial domain and 
$\mathcal{T} = \left[0, T_{f}\right)$ 
is the time interval,
$\mathbf{E}=\left(E_{x}, E_{y}, E_{z}\right)^{T}$ and 
$\mathbf{H}=\left(H_{x}, H_{y}, H_{z}\right)^{T}$
are the electric field and the magnetic field respectively, 
$\varepsilon_{r}$ and $\nu_{r}$
denote the  relative electric  permittivity and  magnetic permeability
parameters  respectively.  In  this work  we consider  the first-order
Silver-Müller    absorbing    boundary   condition    $(\mathrm{ABC})$
\cite{Monk2003}
\begin{equation}
\mathcal{L}(\mathbf{E}(\mathbf{x}, t), \mathbf{H}(\mathbf{x}, t)) = 
\mathbf{n} \times \mathbf{E}(\mathbf{x}, t) + 
Z\mathbf{n} \times(\mathbf{n} \times \mathbf{H}(\mathbf{x}, t)), 
\text { on } \partial \Omega,
\end{equation}
where  $\partial \Omega$  is  the boundary  of $\Omega$,  $\mathbf{n}$
denotes  the  outer  unit   normal  vector  along  
$\partial \Omega$, $\mathbf{E}^{\rm{inc}}$ and  $\mathbf{H}^{\rm{inc}}$ 
are  the incident
fields,  and  $Z = \sqrt{\nu_{r}/{\varepsilon_{r}}}$.   Our goal  is  to
solve equation (\ref{eq:maxw}) with varying parameter $\varepsilon_{r}
\in \mathcal{P}$, where $\mathcal{P}$ denotes the parameter domain.

System   (\ref{eq:maxw})  is   discretized  in   space  by   means  of
discontinuous Galerkin method, and in  time by means of a second-order
leapfrog   scheme.  The   time   interval   $\mathcal{T}  =   \left[0,
  T_{f}\right)$ is  divided into $m$ equal  subintervals as $0=t_{0}<$
  $t_{1}<\cdots<t_{m}=T_{f}$   with   $t_{n}=n  \Delta   t\left(n=0,1,
  \cdots, m\right)$, and  $\Delta t$ denotes the time  step size.  The
  fully  discrete  scheme  of  DGTD   based  on  centered  fluxes  and
  second-order Leap-Frog time stepping \cite{Fezoui2005} is given by
\begin{equation}
\left\{
\begin{array}{l}
\mathbb{M}^{\varepsilon_{r}} 
\dfrac{\underline{\mathbf{E}}_{h}\left(t_{n+1}\right) - 
       \underline{\mathbf{E}}_{h}\left(t_{n}\right)}{\Delta t} = 
\left(\mathbb{K}-\mathbb{S}^{i}\right) 
\underline{\mathbf{H}}_{h}(t_{n+\frac{1}{2}})-\mathbb{S}^{h} 
\underline{\mathbf{H}}_{h}(t_{n+\frac{1}{2}})-\mathbf{B}^{h}(n \Delta t), 
\\ [0.25cm]
\mathbb{M}^{\nu_{r}} 
\dfrac{\underline{\mathbf{H}}_{h}(t_{n+\frac{3}{2}}) - 
       \underline{\mathbf{H}}_{h}(t_{n+\frac{1}{2}})}{\Delta t} = 
\left(-\mathbb{K}+\mathbb{S}^{i}\right) 
\underline{\mathbf{E}}_{h}(t_{n+1}) + 
\mathbb{S}^{e} \widehat{\mathbf{E}}_{h}(t_{n+1}) + 
\mathbf{B}^{e}((n+\frac{1}{2}) \Delta t).
\end{array}
\right.
\end{equation}
Here $\mathbb{M}^{\varepsilon_{r}}$ and $\mathbb{M}^{\nu_{r}}$ are the
mass matrices, $\mathbb{K}$ is  the stiffness matrix, $\mathbb{S}^{i}$
is the surface matrix for the interior faces, and $\mathbb{S}^{h}$ and
$\mathbb{S}^{e}$ are  the boundary  face matrices.  For  more detailed
definition of these matrices we refer to \cite{Viquerat2015, Li2019}.
\section{Two-step POD method}
\label{sec:2step_pod}

\subsection{Snapshot-based data collection}

For the parameterized time-dependent problem, 
let $\mathcal{P}_h^{tr} = \{\boldsymbol{\mu}_{1}, \boldsymbol{\mu}_{2}, 
\cdots, \boldsymbol{\mu}_{\mathcal{N}_p}\}$ be a parameter sampling over
the parameter domain $\mathcal{P}$. 
A collection of high-fidelity solutions can be obtained
by solving system (\ref{eq:maxw}) with the DGTD solver for different 
parameter values $\boldsymbol{\mu} \in \mathcal{P}_h^{tr}$
in a time sampling $\mathcal{T}_h = \{t_{1},t_{2},\cdots,t_{m}\}$. 
Then we uniformly select
$\mathcal{N}_t$ transient solutions in 
$\mathcal{T}_h^{tr} = \{t_{n_1},t_{n_2},\cdots,t_{n_{\mathcal{N}_t}}\} 
\subseteq \mathcal{T}_{h}$ to form the time trajectory matrix 
$\mathbf{S}_{\mathbf{u}}^{j} \in \mathbb{R}^{\mathcal{N}_{h}\times \mathcal{N}_t}$, 
which takes the form
\begin{equation}
\mathbf{S}_{\mathbf{u}}^{j} =
\left(
\begin{array}{cccc}
\mathbf{u}_{h, 1}\left(t_{n_1}, \boldsymbol{\mu}_{j}\right) & 
\mathbf{u}_{h, 1}\left(t_{n_2}, \boldsymbol{\mu}_{j}\right) & \cdots & 
\mathbf{u}_{h, 1}(t_{n_{\mathcal{N}_t}}, \boldsymbol{\mu}_{j})  \\
\mathbf{u}_{h, 2}\left(t_{n_1}, \boldsymbol{\mu}_{j}\right) & 
\mathbf{u}_{h, 2}\left(t_{n_2}, \boldsymbol{\mu}_{j}\right) & \cdots & 
\mathbf{u}_{h, 2}(t_{n_{\mathcal{N}_t}}, \boldsymbol{\mu}_{j})  \\
\vdots & \vdots & \ddots & \vdots \\
\mathbf{u}_{h, \mathcal{N}_{h}}\left(t_{n_1}, \boldsymbol{\mu}_{j}\right) & 
\mathbf{u}_{h, \mathcal{N}_{h}}\left(t_{n_2}, \boldsymbol{\mu}_{j}\right) & 
\cdots & 
\mathbf{u}_{h, \mathcal{N}_{h}} (t_{n_{\mathcal{N}_t}}, \boldsymbol{\mu}_j)
\end{array}\right) 
\text{with} \ j=1,2, \cdots, \mathcal{N}_{p}, 
\mathbf{u} \in \{\mathbf{E}, \mathbf{H}\},
\end{equation}
and the snapshot matrix for all parameters in $\mathcal{P}_{h}^{tr}$ is
\begin{equation}
\mathbf{S}_{\mathbf{u}} = 
\left[\mathbf{S}_{\mathbf{u}}^{1} \ |\ \mathbf{S}_{\mathbf{u}}^{2} 
     \ |\ \cdots \ | \ \mathbf{S}_{\mathbf{u}}^{\mathcal{N}_{p}}\right] 
\in \mathbb{R}^{\mathcal{N}_{h} \times \mathcal{N}_{s}}, 
    \mathbf{u} \in \{\mathbf{E}, \mathbf{H}\},
\end{equation}
where $\mathcal{N}_s = \mathcal{N}_t \cdot \mathcal{N}_p$, 
$\mathbf{u}_h(t,\boldsymbol{\mu}) \in \mathbb{R}^{\mathcal{N}_h}$ 
is the  high-fidelity solution, and  $\mathcal{N}_h$ is the  number of
DOFs, which is determined by  the underlying mesh and polynomial order
of the discretization scheme.

Directly reducing the dimensionality of the snapshots through a neural
network will result  in an overwhelming computational  burden when the
FOM dimension  $\mathcal{N}_h$ becomes  moderately large.  The CAE-CSI
technique  proposed here  can  be considered  as  a non-intrusive  ROM
technique   in  which   a  three-step   dimensionality  reduction   is
performed. First, a  two-step POD strategy is applied on  a set of FOM
snapshots. Then a convolutional autoencoder  is utilized to reduce the
dimensionality of  the projection coefficients (also  called intrinsic
coordinates) generated by  POD. Lastly, a CSI-based model  is built to
approximate     the      mapping     between      input     parameters
$(t,\boldsymbol{\mu})$ and the reduced-order matrices.

\subsection{Two-step POD for dimensionality reduction}

In  this  subsection,  we  perform a  low-rank  approximation  to  the
snapshot   matrix  $\mathbf{S}_{\mathbf{u}}$   and  construct   a  low
dimensional   space   $\mathcal{V}_{\mathbf{u},rb}$   with   dimension
$\mathcal{N}  \ll  \operatorname{min}\{\mathcal{N}_h,\mathcal{N}_s\}$.
Spanned by  a group of  time- and parameter-independent  RB functions,
the reduced space is expressed as
\begin{equation}
\mathcal{V}_{\mathbf{u}, rb} = 
\operatorname{span}\left\{\mathbf{v}_{\mathbf{u}, 1}, 
                          \mathbf{v}_{\mathbf{u}, 2}, \cdots, 
                          \mathbf{v}_{\mathbf{u}, \mathcal{N}}\right\}, 
\mathbf{u} \in \{\mathbf{E}, \mathbf{H}\}.
\end{equation}
The POD method is a popular  technique to compress data and extract an
optimal set of RB functions in the least-squares sense.  The POD basis
of size $\mathcal{N}$ is the solution to the minimization problem
\begin{equation}
\begin{aligned}
 & \min _{\mathbf{V_u} \in \mathbb{R}^{\mathcal{N}_h \times \mathcal{N} }}
   \left\|\mathbf{S_u}-\mathbf{V_u}\mathbf{V}_{\mathbf{u}}^{T} 
          \mathbf{S_u}\right\|_{F}, \\
 & \text { s.t. } \quad \mathbf{V}_{\mathbf{u}}^{T} \mathbf{V_u}=\mathbf{I},
\end{aligned}
\label{eq:pod_ls}
\end{equation}
where $\|\cdot\|_{F}$ is the Frobenius norm and 
$\mathbf{I} \in \mathbb{R}^{\mathcal{N} \times \mathcal{N}}$ is an 
identity matrix. 

According to the  Eckart-Young theorem \cite{Schmidt1907, Eckart1936},
the   solution  of   $(\ref{eq:pod_ls})$   is  given   by  the   first
$\mathcal{N}$  left singular  vectors  of  the matrix  $\mathbf{S_u}$,
which can be obtained through the singular value decomposition (SVD)
\begin{equation}
\mathbf{W}_{\mathbf{u}}^{T} \mathbf{S}_{\mathbf{u}} \mathbf{Z}_{\mathbf{u}} = 
\left(
\begin{array}{cc}
\Sigma_{r_{\mathbf{u}} \times r_{\mathbf{u}}}^{\mathbf{u}} & 
\mathbf{O}_{r_{\mathbf{u}} \times\left(\mathcal{N}_{s}-r_{\mathbf{u}}\right)} \\
\mathbf{O}_{\left(\mathcal{N}_{h}-r_{\mathbf{u}}\right) \times r_{\mathbf{u}}} & 
\mathbf{O}_{\left(\mathcal{N}_{h}-r_{\mathbf{u}}\right) \times
\left(\mathcal{N}_{s}-r_{\mathbf{u}}\right)}
\end{array}
\right) = 
\mathbf{D}_{\mathbf{u}}, \mathbf{u} \in \{\mathbf{E}, \mathbf{H}\},
\end{equation}
where $\mathbf{W}_{\mathbf{u}}=[\mathbf{w}_{\mathbf{u}, 1}, 
\mathbf{w}_{\mathbf{u}, 2}, \cdots, \mathbf{w}_{\mathbf{u}, \mathcal{N}_{h}}]$
and $\mathbf{Z}_{\mathbf{u}}=[\mathbf{z}_{\mathbf{u}, 1}, 
\mathbf{z}_{\mathbf{u}, 2}, \cdots, \mathbf{z}_{\mathbf{u}, \mathcal{N}_{s}}]$
are $\mathcal{N}_{h} \times \mathcal{N}_{h}$ and 
$\mathcal{N}_{s} \times \mathcal{N}_{s}$ orthogonal matrices respectively, 
and $\sum_{r_{\mathbf{u}} \times r_{\mathbf{u}}}^{\mathbf{u}} = 
     \operatorname{diag}\left(\sigma_{\mathbf{u}, 1}, \sigma_{\mathbf{u}, 2}, 
     \cdots, \sigma_{\mathbf{u}, r_{\mathbf{u}}}\right)$.
The singular values $\sigma_{\mathbf{u}, 1} \geq \sigma_{\mathbf{u}, 2} 
\geq \cdots \geq \sigma_{\mathbf{u}, r_{\mathbf{u}}} \geq 0$ of 
$\mathbf{S}_{\mathbf{u}}$
are sorted in descending order, and $r_{\mathbf{u}}$ 
is the rank of $\mathbf{S}_{\mathbf{u}}$.
The POD basis with $\mathcal{N}\left(\mathcal{N}\ll r_{\mathbf{u}}\right)$ 
vectors is the set $\left\{\mathbf{v}_{\mathbf{u}, i}\right\}_{i=1}^{\mathcal{N}}$ 
with $\mathbf{v}_{\mathbf{u}, i}=\mathbf{w}_{\mathbf{u}, i}$,
which can minimize the projection error of the snapshots among all 
$\mathcal{N}$-dimensional orthogonal basis in $\mathbb{R}^{\mathcal{N}_{h}}$.
The error bound can be evaluated using the singular values
\begin{equation}
\sum_{i=1}^{\mathcal{N}_{s}}
\left\|\mathbf{S}_{\mathbf{u}}(:, i)-\mathbf{V_{u}} 
       \mathbf{V}_{\mathbf{u}}^{T} \mathbf{S}_{\mathbf{u}}(:, i)
\right\|_{\mathbb{R}^{\mathcal{N}_{h}}}^{2} = 
\sum_{j=\mathcal{N}+1}^{r_{\mathbf{u}}} \sigma_{\mathbf{u}, j}^{2}, 
\label{eq:pod_err}
\end{equation}
where $\mathbf{V}_{\mathbf{u}}=\left[\mathbf{v}_{\mathbf{u}, 1}, \cdots, 
\mathbf{v}_{\mathbf{u}, \mathcal{N}}\right]$, $\mathbf{u} \in \{\mathbf{E}, 
\mathbf{H}\}$.
According to (\ref{eq:pod_err}), it is clear that the error in the POD
basis is  equal to the  sum of the  squares of the  neglected singular
values, i.e., the  error can be controlled by  $\mathcal{N}$.  In this
work, POD is  utilized to perform a  moderate dimensionality reduction
of the  snapshots data, thus  yielding a linear subspace  of dimension
$\mathcal{N} \ll \mathcal{N}_h$.

However, since  $\mathcal{N}_s = \mathcal{N}_t\cdot  \mathcal{N}_p$ is
large, performing SVD  on such a large-scale  snapshot matrix directly
is extremely expensive. We adopt  a two-step POD method to effectively
save the  computational cost.  The detailed process  is as  follow and
also shown in Algorithm \ref{algo:2step_pod}:

\begin{enumerate}
\item For each time trajectory matrix  $\mathbf{S}_{\mathbf{u}}^j$, 
      $j=$ $1,2, \cdots, \mathcal{N}_{p}$,
      the POD basis $\{\gamma_{\mathbf{u},i}^{j}\}_{i=1}^{k}$ are obtained 
      through SVD, then assemble them to a matrix 
      $\mathbf{T}_{\mathbf{u}}^{j}=[\gamma_{\mathbf{u}, 1}^{j}, 
      \gamma_{\mathbf{u}, 2}^{j}, \cdots, \gamma_{\mathbf{u},k}^{j}]$;
\item Assemble a composite matrix 
      $\mathbf{T}_{\mathbf{u}} = [\mathbf{T}_{\mathbf{u}}^{1}\ | \ 
                                \mathbf{T}_{\mathbf{u}}^{2}\ | \ \cdots \ | \ 
                                \mathbf{T}_{\mathbf{u}}^{\mathcal{N}_{p}}]$ 
      with the matrices generated in the first step,
      then perform POD on $\mathbf{T}_{\mathbf{u}}$ with truncation 
      $\mathcal{N}$. The POD basis
      $\{\mathbf{v}_{\mathbf{u}, i}\}_{i=1}^{\mathcal{N}}$ is obtained, i.e., 
      $\mathbf{V}_{\mathbf{u}}=[\mathbf{v}_{\mathbf{u}, 1}, \cdots, 
       \mathbf{v}_{\mathbf{u}, \mathcal{N}}]$.
\end{enumerate}
According to the  projection theory \cite{Li2018}, one  can obtain the
reduced-order   solution  and   the  
projection coefficients or intrinsic 
coordinates as
\begin{equation}
\left\{
\begin{array}{l}
\mathbf{u}_{h}^{r}(t, \boldsymbol{\mu}) = 
\mathbf{V}_{\mathbf{u}}
\left(\mathbf{V}_{\mathbf{u}}^{T} \mathbf{V}_{\mathbf{u}}\right)^{-1} 
\mathbf{V}_{\mathbf{u}}^{T} \mathbf{u}_{h}(t, \boldsymbol{\mu}) = 
\mathbf{V}_{\mathbf{u}} \mathbf{V}_{\mathbf{u}}^{T} 
\mathbf{u}_{h}(t, \boldsymbol{\mu}), \\ [0.25cm]
\alpha_{\mathbf{u}}(t, \boldsymbol{\mu})\equiv
\mathbf{u}_{\mathcal{N}}(t,\boldsymbol{\mu}) = 
\left(\mathbf{V}_{\mathbf{u}}^{T} \mathbf{V}_{\mathbf{u}}\right)^{-1} 
\mathbf{V}_{\mathbf{u}}^{T} \mathbf{u}_{h}(t, \boldsymbol{\mu}) = 
\mathbf{V}_{\mathbf{u}}^{T} \mathbf{u}_{h}(t, \boldsymbol{\mu}),
\end{array} 
\mathbf{u} \in\{\mathbf{E}, \mathbf{H}\}.
\right.
\end{equation}
Thus, the reduced-order solution $\mathbf{u}_{h}^{r}(t, \boldsymbol{\mu})$ 
serves as an approximation to the high-fidelity solution 
$\mathbf{u}_{h}(t, \boldsymbol{\mu})$ and can be represented as
\begin{equation}
\mathbf{u}_{h}^{r}(t, \boldsymbol{\mu}) = 
\mathbf{V_u}\alpha_{\mathbf{u}}(t,\boldsymbol{\mu}) =
\sum_{i=1}^{\mathcal{N}} 
\alpha_{\mathbf{u}, i}(t, \boldsymbol{\mu}) \mathbf{v}_{\mathbf{u}, i}, 
                        \mathbf{u}\in \{\mathbf{E}, \mathbf{H}\},
\end{equation}
where $\alpha_{\mathbf{u}}(t, \boldsymbol{\mu}) =
\left[\alpha_{\mathbf{u}, 1}(t, \boldsymbol{\mu}), 
      \alpha_{\mathbf{u}, 2}(t, \boldsymbol{\mu}), \cdots, 
      \alpha_{\mathbf{u}, \mathcal{N}}(t, \boldsymbol{\mu})\right]^{T}
\in \mathbb{R}^{\mathcal{N}}$ collects the combination coefficients.

In  our  previous  work  \cite{Zhao2021, Li2021}  we  constructed  the
mapping  between  the  time/parameter  values  and  the  projection coefficients   by  decoupling   time-  and   parameters-modes  through
SVD.  However,  these ROMs  need  to  create  too many  regression  or
interpolation models.  The approach  we proposed here  firstly reduces
the   length   of  projection   coefficients   $\alpha_{\mathbf{u}}(t,
\boldsymbol{\mu})$          or          intrinsic          coordinates
$\mathbf{u}_{\mathcal{N}}(t,\boldsymbol{\mu})$ from $\mathcal{N}$ to a
fairly small $n$ through  a convolutional autoencoder, then constructs
the mapping between the  time/parameter values and the low-dimensional
coded representation  using cubic spline interpolation,  thus reducing
the number of interpolation models as well as test time online.

\begin{algorithm}[htbp]
\label{algo:2step_pod}
\SetKwProg{Fn}{Function}{}{end}
\SetKwFunction{POD}{POD}
\SetKwInOut{Input}{Input}\SetKwInOut{Output}{Output}

\caption{Two-step POD method}\label{algo:algorithm_two-step_POD}
\Input{Time trajectory matrices 
       $\mathbf{S}_{\mathbf{u}}^{j}
        \left(j=1,2, \cdots, \mathcal{N}_{p}, 
              \mathbf{u} \in\{\mathbf{E}, \mathbf{H}\}\right)$, 
       truncation parameter $k$ and size parameter $\mathcal{N}$}
\Output{POD basis matrix 
        $\mathbf{V}_{\mathbf{u}}\ (\mathbf{u} \in\{\mathbf{E}, \mathbf{H}\})$}
\BlankLine
Compute the compressed matrices 
$\mathbf{T}_{\mathbf{u}}^{j} = $\POD$(\mathbf{S}_{\mathbf{u}}^{j}, k )$ 
for $j=1,2, \cdots, \mathcal{N}_{p}, \mathbf{u} \in\{\mathbf{E}, 
     \mathbf{H}\}$; \\
Assemble the matrix 
$\mathbf{T}_{\mathbf{u}} = [\mathbf{T}_{\mathbf{u}}^{1}\ | \ 
 \mathbf{T}_{\mathbf{u}}^{2}\ | \ \cdots \ | \  
 \mathbf{T}_{\mathbf{u}}^{\mathcal{N}_{p}}]\ 
 (\mathbf{u} \in\{\mathbf{E}, \mathbf{H}\})$ for all parameter values; \\
Calculate the POD basis matrix 
$\mathbf{V}_{\mathbf{u}} = 
$\POD$\left(\mathbf{T}_{\mathbf{u}}, \mathcal{N}\right)\ 
(\mathbf{u} \in\{\mathbf{E}, \mathbf{H}\})$
\BlankLine
\Fn{$\mathbf{V}=$ \POD$(\mathbf{A}, k)$}{
Perform $\operatorname{eig}(\mathbf{A}^{T}\mathbf{A})$ 
in MATLAB, and get the eigenvalues 
$\lambda_{1} \geq \lambda_{2} \geq \cdots \geq \lambda_{r}>0$ 
and the corresponding orthogonal eigenvectors  
$\mathbf{u}_{1}, \mathbf{u}_{2}, \cdots, \mathbf{u}_{r}$ 
with $r$ being the rank of $\mathbf{A}$; \\
Compute the POD basis functions 
$\mathbf{v}_{i}=\dfrac{1}{\sqrt{\lambda_{i}}} \mathbf{A}
 \mathbf{u}_{i}, i=1,2, \cdots, k, k \ll r$; \\
Obtain the POD basis matrix $\mathbf{V}=[\mathbf{v}_{1}, 
                             \mathbf{v}_{2}, \cdots, \mathbf{v}_{k}]$.
}
\end{algorithm}

\begin{remark}
According to  (\ref{eq:pod_err}), the  error bounds  in the  first and
second steps of the two-step POD algorithm are written as
\begin{equation}
\left\{
\begin{array}{l}
\sum \limits_{i=1}^{\mathcal{N}_{t}}
\left\|\mathbf{S}_{\mathbf{u}}^{j}(:, i) - 
       \mathbf{T}_{\mathbf{u}}^{j} (\mathbf{T}_{\mathbf{u}}^{j})^{T} 
       \mathbf{S}_{\mathbf{u}}^{j}(:, i)
\right\|_{\mathbb{R}^{\mathcal{N}_{h}}}^{2} =
\sum \limits_{i=k+1}^{r_{\mathbf{u}}^{j}}(\sigma_{\mathbf{u}, i}^{j})^{2}, 
1 \leqslant j \leqslant \mathcal{N}_{p}, \\
\sum \limits_{j=1}^{\mathcal{N}_{p}} \sum \limits_{i=1}^{k}
\left\|\mathbf{T}_{\mathbf{u}}^{j}(:, i) -
       \mathbf{V}_{\mathbf{u}} \mathbf{V}_{\mathbf{u}}^{T} 
       \mathbf{T}_{\mathbf{u}}^{j}(:, i)
\right\|_{\mathbb{R}^{\mathcal{N}_{h}}}^{2} =
\sum\limits_{i=\mathcal{N}+1}^{r_{\mathbf{u}}}\left(\sigma_{\mathbf{u}, i}\right)^{2},
\end{array}
\right.
\end{equation}
where $r_{\mathbf{u}}$ and 
$r_{\mathbf{u}}^{j}\ 
 \left(j=1,2, \cdots, \mathcal{N}_{p}, 
       \mathbf{u}\in \{\mathbf{E}, \mathbf{H}\}\right)$ 
are the rank of $\mathbf{S}_{\mathbf{u}}$ and 
$\mathbf{S}_{\mathbf{u}}^{j}$, and 
$\left\{\sigma_{\mathbf{u}, i}\right\}_{i=1}^{r_{\mathbf{u}}}$ as well as 
$\{\sigma_{\mathbf{u}, i}^{j}\}_{i=1}^{r_{\mathbf{u}}^{j}}$ 
are the corresponding singular values. 
The two-step POD projection error can be bounded as
\begin{equation*}
\begin{aligned}
\sum_{j=1}^{\mathcal{N}_{p}} \sum_{i=1}^{\mathcal{N}_{t}}
\left\|\mathbf{S}_{\mathbf{u}}^{j}(:, i) - 
       \mathbf{V}_{\mathbf{u}} \mathbf{V}_{\mathbf{u}}^{T} 
       \mathbf{S}_{\mathbf{u}}^{j}(:, i)
\right\|_{\mathbb{R}^{\mathcal{N}_{h}}} & 
\leqslant \sum_{j=1}^{\mathcal{N}_{p}} \sum_{i=1}^{\mathcal{N}_{t}}
\left\|\mathbf{S}_{\mathbf{u}}^{j}(:, i) - 
       \mathbf{T}_{\mathbf{u}}^{j} (\mathbf{T}_{\mathbf{u}}^{j})^{T} 
       \mathbf{S}_{\mathbf{u}}^{j}(:, i)\right\|_{\mathbb{R}^{\mathcal{N}_{h}}} \\
& 
+ \sum_{j=1}^{\mathcal{N}_{p}} \sum_{i=1}^{\mathcal{N}_{t}}
\left\|\mathbf{T}_{\mathbf{u}}^{j} (\mathbf{T}_{\mathbf{u}}^{j})^{T} 
       \mathbf{S}_{\mathbf{u}}^{j}(:, i) - 
       \mathbf{V}_{\mathbf{u}} \mathbf{V}_{\mathbf{u}}^{T} 
       \mathbf{T}_{\mathbf{u}}^{j} 
       (\mathbf{T}_{\mathbf{u}}^{j})^{T} \mathbf{S}_{\mathbf{u}}^{j}(:, i)
\right\|_{\mathbb{R}^{\mathcal{N}_{h}}} \\
& 
+ \sum_{j=1}^{\mathcal{N}_{p}} \sum_{i=1}^{\mathcal{N}_{t}}
\left\|\mathbf{V}_{\mathbf{u}} \mathbf{V}_{\mathbf{u}}^{T} 
       \mathbf{T}_{\mathbf{u}}^{j} (\mathbf{T}_{\mathbf{u}}^{j})^{T} 
       \mathbf{S}_{\mathbf{u}}^{j}(:, i) - 
       \mathbf{V}_{\mathbf{u}} \mathbf{V}_{\mathbf{u}}^{T} 
       \mathbf{S}_{\mathbf{u}}^{j}(:, i)
\right\|_{\mathbb{R}^{\mathcal{N}_{h}}} \\
& 
\leqslant \mathcal{L}_{1} + \mathcal{L}_{2}, 
\mathbf{u} \in \{\mathbf{E}, \mathbf{H}\},
\end{aligned}
\end{equation*}
where
\begin{equation*}
\left\{
\begin{array}{l}
\mathcal{L}_{1} = 
(1 + \left\|\mathbf{V}_{\mathbf{u}} \mathbf{V}_{\mathbf{u}}^{T}\right\|_{F}) 
     \sum\limits_{j=1}^{\mathcal{N}_{p}}(\mathcal{N}_{t} 
     \sum\limits_{i=k+1}^{r_{\mathbf{u}}^{j}}
     (\sigma_{\mathbf{u}, i}^{j})^{2})^{\frac{1}{2}}, \\
\mathcal{L}_{2} = 
\max\limits_{1 \leqslant j \leqslant \mathcal{N}_{p}} 
\sum\limits_{i=1}^{\mathcal{N}_{t}}
\|\mathbf{S}_{\mathbf{u}}^{j}(:, i)\|_{\mathbb{R}^{\mathcal{N}_{h}}} 
\max\limits_{\substack{1 \leqslant i \leqslant k \\ 
1 \leqslant j \leqslant \mathcal{N}_{p}}}
\|\mathbf{T}_{\mathbf{u}}^{j}(:, i)\|_{\mathbb{R}^{\mathcal{N}_{h}}}(k\mathcal{N}_p 
\sum \limits_{i=\mathcal{N}+1}^{r_{\mathbf{u}}}(\sigma_{\mathbf{u}, i})^{2})^{\frac{1}{2}}.
\end{array}
\right.
\end{equation*}
Therefore, the accuracy  of the two-step POD can be  controlled by the
truncation parameter $k$ and size parameter $\mathcal{N}$.
\end{remark}
\section{Convolutional autoencoders for model reduction}
\label{sec:cae}

In  this section,  we present  an  approach to  construct a  nonlinear
manifold,   which   compresses   the  dimensionality   of   projection
coefficients
$\alpha_{\mathbf{u}}(t, \boldsymbol{\mu}) \equiv  
 \mathbf{u}_{\mathcal{N}}(t,\boldsymbol{\mu}) \in \mathbb{R}^{\mathcal{N}}$ 
$(\mathbf{u} \in \{\mathbf{E}, \mathbf{H}\})$ from 
$\mathcal{N}$ to a fairly small $n$. In
data-driven sciences, the manifold hypothesis presumes that real-world
high   dimensional  data   lie   near  a   low  dimensional   manifold
$\mathcal{S}$  embedded  in  $\mathbb{R}^{m}$,   where  $m$  is  large
\cite{Bengio2013}.  As  a   result  POD  has  been   widely  used  for
dimensionality  reduction  of  physical  systems.   However  the  main
drawback is  that POD  can only construct  an optimal  linear manifold
while  data  sampled  from  complex real-world  systems  tends  to  be
strongly nonlinear.

\begin{figure}[htbp]
\centering
\includegraphics[width=0.4\textwidth]{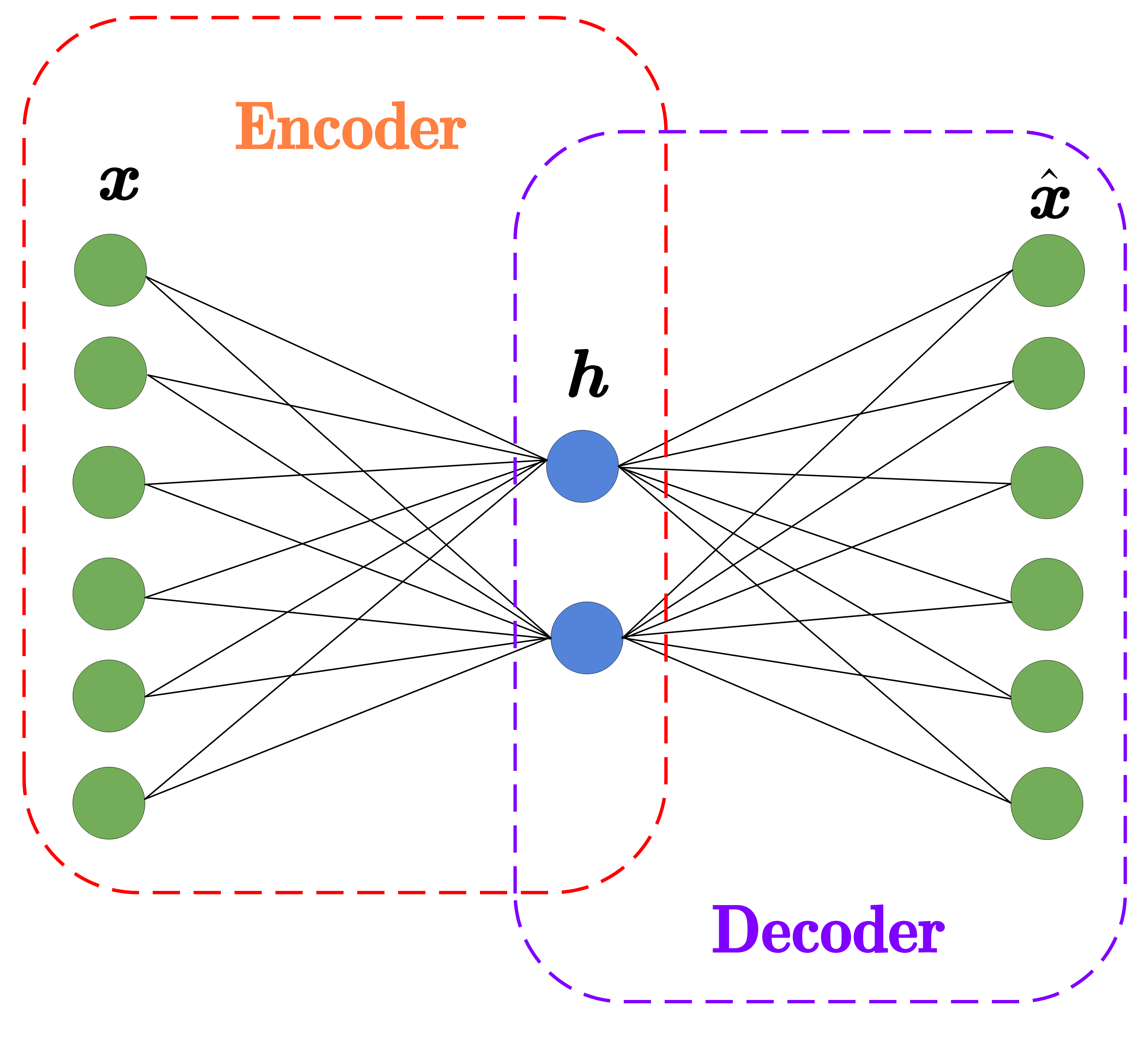}
\caption{Network architecture of a fully-connected autoencoder.}
\label{AE}
\end{figure}

A nonlinear generalization of POD  is the autoencoder, which takes the
the   form   of   a    fully-connected   neural   network   shown   in
Fig.~\ref{AE}. Specifically,  an autoencoder  includes two  parts: the
encoder  maps a  high dimensional  input to  a low-dimensional  latent
vector, known as encoding; then  the decoder learns how to reconstruct
the original input from the latent vector with a minimum error.
A basic  autoencoder consists  of a single  or multiple  layer encoder
network
\begin{equation}
\mathbf{h} = f_{E}\left(\mathbf{x} ; \boldsymbol{\theta}_{E}\right),
\end{equation}
where $\boldsymbol{\theta}_E$  denotes the  weights and the  biases of
the encoder network, $\mathbf{x}  \in \mathbb{R}^{\mathcal{N}}$ is the
input state,  $\mathbf{h} \in  \mathbb{R}^{n}$ is  the feature  or low
dimensional representation  vector with $n \ll  \mathcal{N}$. A decoder
network is then used to reconstruct $\mathbf{x}$ by
\begin{equation}
\hat{\mathbf{x}} = f_{D}\left(\mathbf{h} ; \boldsymbol{\theta}_{D}\right).
\end{equation}
The objective is how to train this autoencoder and find the parameters
that minimize the mean squared error over all training examples
\begin{equation}
\boldsymbol{\theta}_{E}^{*}, \boldsymbol{\theta}_{D}^{*} = 
\mathop{\arg\min}\limits_{\boldsymbol{\theta}_{E}, \boldsymbol{\theta}_{D}}
\|\mathbf{x} - \hat{\mathbf{x}}\|^{2}.
\end{equation}
Note that  directly applying large-scale snapshots  to fully connected
autoencoders is  not only computationally expensive,  but also ignores
opportunities to exploit feature  structures in high dimensional data.
So we  adopt an  convolutional autoencoder  which is  characterized by
shared parameters  and local connectivity,  thus reduce the  memory as
well as computational cost.

\subsection{Basics of convolutional and deconvolutional layers}

Some basic  operations in  convolutional and  transposed convolutional
layers are introduced  in this subsection.  The reader  is referred to
\cite{Dumoulin2016,    Zhang2021}   for    more    details   of    the
convolution operations.

For   simplicity,  we   consider   a  3-D   tensor  $\mathcal{X}   \in
\mathbb{R}^{c\times  h\times  w}$ as  an  input  of the  convolutional
neural network,  with $c$, $h$ and  $w$ being the number  of channels,
height  and width  respectively.   The convolution  operations aim  to
extract the most  important features from the  input $\mathcal{X}$ and
then use them to construct  a feature map $\mathcal{G}$. Specifically,
the  output   at  the   $(l-1)$th  layer  is   $\mathcal{X}^{l-1}  \in
\mathbb{R}^{c_{q}^{l-1}\times h^{l-1}\times w^{l-1}}$, and the feature
map  at  $l$th layer  is  denoted  as  a tensor  $\mathcal{G}^{l}  \in
\mathbb{R}^{c_{p}^{l}   \times  n_{1}^{l}   \times  n_{2}^{l}}$   with
$c_{p}^{l}  =c_{q}^{l-1} $  and  element  $\mathcal{G}_{i, j,  k}^{l}$
representing a unit within channel $i$ at row $j$ and column $k$.  The
filter  banks collecting  a  set  of kernels  at  $l$th  layer can  be
considered   as   a    4-dimensional   tensor   $\mathcal{W}^{l}   \in
\mathbb{R}^{c_{p}^{l}   \times  c_{q}^{l}   \times  k_{h}^{l}   \times
  k_{w}^{l}}$ with  element $\mathcal{W}_{i, j, m,  n}^{l}$ connecting
between a unit in channel $i$ of  the output and a unit in channel $j$
of the input, with  an offset of $m$ rows and  $n$ columns between the
output unit and the input unit.  $c_{p}$ denotes the number of kernels
in the  filter bank,  and the  kernel size is  denoted by  $k_{h}$ and
$k_{w}$. Convolution  between a feature map  $\mathcal{G}^{l-1}$ and a
filter bank $\mathcal{W}^{l}$ can be expressed as
\begin{equation}
\mathcal{G}_{i,j,k}^{l} = 
\sigma_{l}\left(\sum_{r, m, n} 
               \mathcal{G}_{r,(j-1) \times s+m,(k-1) \times s+n}^{l-1} 
               \mathcal{W}_{i, r, m, n}^{l} + 
               \mathcal{B}_{i, j, k}^{l}\right), 
\end{equation}
where   $\sigma_l$   is   a    nonlinear   activation   function   and
$\mathcal{B}^{l}$  is a  bias.  Here,  $s$ denotes  the stride,  which
determines the downsampling rate of each convolution. The dimension of
the next feature map can be reduced by a factor  $s$ in each direction
when $s>1$.  The  kernel size $[k_{h}, k_{w}]$, the  number of filters
$c_{p}$ and  the stride $s$  are the hyperparameters while  the filter
banks  $\mathcal{W}$  and  the   biases  $\mathcal{B}$  are  learnable
parameters. A deconvolutional layer performs the reverse operations of
convolution,  called  transposed  convolution,   and  it  is  used  to
construct   decoding   layers   \cite{Lee2020}.    In   summary,   the
architecture   of  CAE   shown  in   Fig.~\ref{fig:cae} is composed  of
convolutional,  deconvolutional and  dense layers,  which is  used for
dimensionality reduction and reconstruction.

\subsection{Data preparation}

In this subsection, we present a  procedure to generate a data set for
training.
Provided the high-fidelity snapshots and the orthogonal
basis $\mathbf{V_u} \in \mathbb{R}^{{\mathcal{N}_h}\times \mathcal{N}}$ 
      ($\mathbf{u}\in \{\mathbf{E}, \mathbf{H}\}$) generated by the 
two-step POD method,
we compute the projection coefficients or intrinsic coordinates 
$\alpha_{\mathbf{u}}(t,\boldsymbol{\mu}) = 
 \mathbf{u}_{\mathcal{N}}(t, \boldsymbol{\mu}) = 
 \mathbf{V}_{\mathbf{u}}^{T} \mathbf{u}_{h}(t, \boldsymbol{\mu})$ 
 ($\mathbf{u}\in \{\mathbf{E}, \mathbf{H}\}$) to perform a moderate 
dimensionality reduction, and the projection coefficient matrix 
$\mathbf{C_u}$ takes the form
\begin{equation}
\mathbf{C}_{\mathbf{u}} = 
\left[\mathbf{C}_{\mathbf{u}}^{1} \ | \  
      \mathbf{C}_{\mathbf{u}}^{2} \ | \ \cdots \ | \  
      \mathbf{C}_{\mathbf{u}}^{\mathcal{N}_{p}}
\right] \in \mathbb{R}^{\mathcal{N} \times \mathcal{N}_{s}}, \ 
            \mathbf{u} \in \{\mathbf{H}, \mathbf{E}\},
\end{equation}
where $\mathbf{C}_{\mathbf{u}}^{j} = 
\mathbf{V}_{\mathbf{u}}^{T} \mathbf{S}_{\mathbf{u}}^{j} \in 
\mathbb{R}^{\mathcal{N}\times \mathcal{N}_{t}}$, $j = 1, 2, \cdots \mathcal{N}_p$.
Next, we randomly shuffle $\mathbf{C_u}$ by column and split it 
into training set and validation set according to a training-validation 
splitting fraction $\lambda$, i.e.,
$\mathbf{C_{u}} = 
 [\mathbf{C}_{\mathbf{u}}^{\rm{tr}}, \mathbf{C}_{\mathbf{u}}^{\rm{val}}]$, 
 $\mathbf{u} \in \{\mathbf{H}, \mathbf{E}\}$, where 
$\mathbf{C}_{\mathbf{u}}^{\rm{tr}} \in 
 \mathbb{R}^{\mathcal{N}\times \lambda \mathcal{N}_s}$.
As the  training speed of neural  network is affected by  the range of
input/output  \cite{Ioffe2015},  feature  scaling is  required  before
feeding  the   training  data  into   the  network.  The   input  data
$\mathbf{C}_{\mathbf{u}}^{\rm{tr}}$ and
$\mathbf{C}_{\mathbf{u}}^{\rm{val}}$ are normalized to $[0,1]$ through
the following affine transformation
\begin{equation}
\mathbf{C}_{\mathbf{u}}^{\rm{tr}}(i,j) =
\frac{\mathbf{C}_{\mathbf{u}}^{\rm{tr}}(i,j) - 
      \mathop{\min}\limits_{1 \leqslant i \leqslant \mathcal{N}}^{} 
      \mathop{\min}\limits_{1 \leqslant j \leqslant \lambda\mathcal{N}_s} 
      \mathbf{C}_{\mathbf{u}}^{\rm{tr}}(i,j)}
     {\mathop{\max}\limits_{1 \leqslant i \leqslant \mathcal{N}} 
      \mathop{\max}\limits_{1 \leqslant j \leqslant \lambda\mathcal{N}_s} 
      \mathbf{C}_{\mathbf{u}}^{\rm{tr}}(i,j) - 
      \mathop{\min}\limits_{1 \leqslant i \leqslant \mathcal{N}} 
      \mathop{\min}\limits_{1 \leqslant j \leqslant \lambda\mathcal{N}_s} 
      \mathbf{C}_{\mathbf{u}}^{\rm{tr}}(i,j)}.
\end{equation}
We then reshape each column of  $\mathbf{C_u}$ into a square matrix of
dimension     $(\sqrt{\mathcal{N}}, \sqrt{\mathcal{N}})$,     where
$\mathcal{N}=m^{2}$  with   $m \in \mathbb{N}$.   Then   stack  the
components of $\mathbf{H}$ and $\mathbf{E}$  together to form a tensor
with  $d$  channels,  where  $d$   denotes  the  number  of  vectorial
components of the solution of system (\ref{eq:maxw}).  Thus, the input
of   the    autoencoder   network    is   a   tensor    of   dimension
$(\sqrt{\mathcal{N}}, \sqrt{\mathcal{N}}, d)$,  which allows to reduce
the number of  parameters, thus save the time of  training and testing
the network.

\subsection{Dimensionality reduction via convolutional autoencoders}

\begin{figure}[htbp]
\centering
\includegraphics[width=0.9\textwidth]{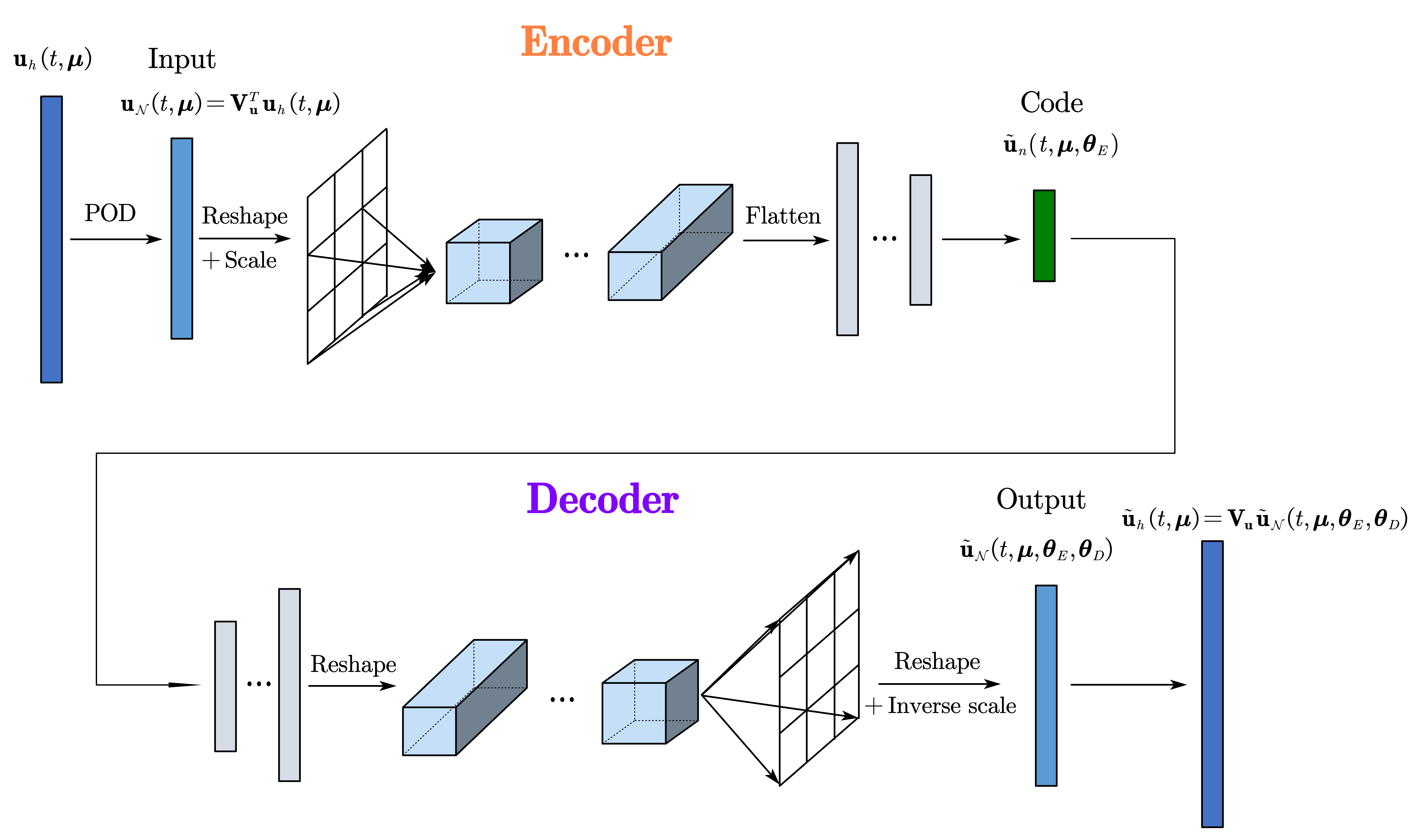}
\caption{Architecture of a deep convolutional autoencoder network used
  for parameterized electromagnetic scattering problems, which takes a
  projection  coefficient  as an  input  and  produces an  approximate
  projection  coefficient  as  an   output.  A  low  dimensional  code
  $\tilde{\mathbf{u}}_n(t, \boldsymbol{\mu}, \boldsymbol{\theta}_{E})$
  is  extracted  via  the  encoder  from  the  projection  coefficient
  $\mathbf{u}_\mathcal{N}(t, \boldsymbol{\mu})$. The  decoder
  reconstructs an approximate projection  coefficient
  $\tilde{\mathbf{u}}_{\mathcal{N}}(t, \boldsymbol{\mu},
  \boldsymbol{\theta}_{E},  \boldsymbol{\theta}_{D})$   from  the  low
  dimensional code  $\tilde{\mathbf{u}}_n(t, \boldsymbol{\mu},
  \boldsymbol{\theta}_{E})$  by performing  the inverse  operations of
  the encoder.}
\label{fig:cae}
\end{figure}

The architecture of  the CAE, employed at training stage,  is shown in
Fig.~\ref{fig:cae}. The  loss function to measure  the discrepancy 
between the input
$\mathbf{u}_{\mathcal{N}}(t, \boldsymbol{\mu}) = 
 \mathbf{V}_{\mathbf{u}}^{T}\mathbf{u}_{h}(t,\boldsymbol{\mu})$
and its reconstruction 
$\tilde{\mathbf{u}}_{\mathcal{N}}(t, \boldsymbol{\mu}, 
                                   \boldsymbol{\theta}_{E}, 
                                   \boldsymbol{\theta}_D)$ 
is defined as follow 
\begin{equation}
\min _{\boldsymbol{\theta}} \mathcal{J}(\boldsymbol{\theta}) = 
\min _{\boldsymbol{\theta}} \frac{1}{\mathcal{N}_{s}} 
\sum_{j=1}^{\mathcal{N}_p} \sum_{i=1}^{\mathcal{N}_{t}}
\left\|\mathbf{V}_{\mathbf{u}}^{T}\mathbf{u}_{h}(t_i, \boldsymbol{\mu}_j) - 
       \tilde{\mathbf{u}}_{\mathcal{N}}(t_i, \boldsymbol{\mu}_j, 
                                           \boldsymbol{\theta}_{E}, 
                                           \boldsymbol{\theta}_D)
\right\|^{2},
\end{equation}
and $\boldsymbol{\theta}=\left(\boldsymbol{\theta}_{E},
\boldsymbol{\theta}_{D}\right)$.   The CAE  is implemented  in PyTorch
\cite{Paszke2019} and the Adam  optimizer \cite{Kingma2014} is used in
the training  stage.  We set  the initial  learning rate to  $\eta_0 =
10^{-4}$,  the mini-batch  size to  $\mathcal{N}_{\rm{mb}} =  50$, and
maximum  number  of  epochs  to $\mathcal{N}_{\rm{epo}}  =  5000$.   A
learning rate decay  with $\eta = \eta_0 /(1+\alpha  * \rm{epoch})$ is
used to  accelerate the training  with hyperparameter $\alpha  = 0.05$
\cite{Wang2019}.  The dataset is  divided into training and validation
set  with  a  proportion  $8:2$,  i.e., $\lambda  =  0.8$.   To  avoid
overfitting, we stop  the training if the loss on  validation set does
not decrease over  500 epochs.  The ELU function defined  as follow is
utilized as the nonlinear activation function
\begin{equation}
\sigma(z) = \begin{cases} z & z \geq 0 \\ \exp (z)-1 & z < 0 \end{cases}.
\end{equation}
No activation function  is applied at the last  convolutional layer of
the decoder neural network.  The weights and biases of the network are
initialized through  the Xavier initialization  \cite{Glorot2010}. The
training process of the CAE is shown in Algorithm \ref{algo:cae}.
\begin{algorithm}[htbp]
\caption{Training process of the CAE}
\label{algo:cae}
\SetAlgoLined
\newcommand{\forcond}{$i=0$ \KwTo $n$}
\SetKwProg{Fn}{Function}{}{end}
\SetKwFunction{CAE}{CAE\_TRAINING}%
\KwIn{Snapshot matrix $\mathbf{S_u}$ and POD basis matrix 
      $\mathbf{V_{u}}$ with $\mathbf{u} \in \{\mathbf{H}, \mathbf{E}\}$,
      training-validation splitting fraction $\lambda$, 
      initial learning rate $\eta_0$, learning rate decay parameter $\alpha$,
      mini-batch size $\mathcal{N}_{\rm{mb}}$ and maximum number of 
      epochs $\mathcal{N}_{\rm{epo}}$.}
\KwOut{Optimal model weights 
       $\boldsymbol{\theta}^{*} = 
        (\boldsymbol{\theta}_{E}^{*}, \boldsymbol{\theta}_{D}^{*})$.}
Compute intrinsic coordinate matrix 
$\mathbf{C}_{\mathbf{u}}^{j} = 
 \mathbf{V}_{\mathbf{u}}^{T} \mathbf{S}_{\mathbf{u}}^{j}$ 
for $j = 1,2,\cdots, \mathcal{N}_{p}, 
     \mathbf{u} \in \{\mathbf{H}, \mathbf{E}\}$\;
Shuffle $\mathbf{C_u} = [\mathbf{C}_{\mathbf{u}}^1, 
         \mathbf{C}_{\mathbf{u}}^2, \cdots, 
         \mathbf{C}_{\mathbf{u}}^{\mathcal{N}_p}]$ randomly by column\;
Split data into training and validation set 
$\mathbf{C_u} = [\mathbf{C}_{\mathbf{u}}^{\rm{tr}}, 
                 \mathbf{C}_{\mathbf{u}}^{\rm{val}}]$ with 
$\mathbf{C}^{\rm{tr}} \in \mathbb{R}^{\mathcal{N}_{h}\times \lambda \mathcal{N}_s}$\;
Obtain the optimal model parameters 
$[\mathbf{W}^{*}, \mathbf{b}^{*}]$ =  
\CAE{$\mathbf{C_u}, \mathcal{N}_{\rm{epo}}, \mathcal{N}_{\rm{b}},\eta_0$}\;
\BlankLine
\Fn{$[\mathbf{W}_{\rm{tr}}, \mathbf{b}_{\rm{tr}}] =$  
\CAE{$\mathcal{D}, \mathcal{N}_{\rm{epo}}, \mathcal{N}_{\rm{b}},\eta_0$}}{
$\mathbf{W}, \mathbf{b} \leftarrow \rm{Initialize}(\mathbf{W}, 
             \mathbf{b})$ \hfill{$\rhd$ initialize weights and biases} \\
$\mathcal{N}_{b} \leftarrow \mathcal{N}_s/\mathcal{N}_{\rm{mb}}$ 
\hfill{$\rhd$ number of batches} \\
\For{$ \rm{epoch} \leftarrow 1$ \KwTo $\mathcal{N}_{\rm{epo}}$}{
     $\eta \leftarrow \rm{Learning\_rate\_decay}(\eta_0, \alpha, \rm{epoch})$ 
     \hfill{$\rhd$ learning rate decay} \\
\For{$s \leftarrow 1$  \KwTo $\mathcal{N}_{\rm{b}}$}{
     $\mathcal{D}_i \leftarrow 
      \mathcal{D}[(s-1)*\mathcal{N}_{\rm{mb}}+1:s*\mathcal{N}_{\rm{mb}}]$   
     \hfill{$\rhd$ the s-th mini-batch} \\
$\mathcal{D}_i \leftarrow 
 \mathrm{Reshape}(\mathcal{D}_i, (\mathcal{N}_{\rm{mb}}, \sqrt{\mathcal{N}}, 
                                  \sqrt{\mathcal{N}}, d))$   
 \hfill{$\rhd$ reshape to $d$ channels} \\
$\mathcal{J} \leftarrow 
 \frac{1}{\mathcal{N}_{s}} 
 \sum\limits_{j=1}^{\mathcal{N}_p} \sum\limits_{i=1}^{\mathcal{N}_{t}}
 \left\|\mathbf{V}_{\mathbf{u}}^{T}\mathbf{u}_{h}(t_i, \boldsymbol{\mu}_j) -
        \tilde{\mathbf{u}}_{\mathcal{N}}(t_i, \boldsymbol{\mu}_j)
\right\|^{2}$ \hfill{$\rhd$ compute loss} \\
$\Delta \mathbf{W} \leftarrow -
 \eta \mathrm{G}_{\mathrm{Adam}}
 \left(\nabla_{\mathbf{W}} \mathcal{J}\right), \Delta \mathbf{b} 
 \leftarrow-\eta \mathrm{G}_{\mathrm{Adam}}
 \left(\nabla_{\mathbf{b}} \mathcal{J}\right)$   
 \hfill{$\rhd$ Adam optimizer step} \\
$\mathbf{W} \leftarrow \mathbf{W} + 
 \Delta \mathbf{W}, \mathbf{b} \leftarrow 
 \mathbf{b}+\Delta \mathbf{b}$ 
 \hfill{$\rhd$ update weights and biases}
}
}$\mathbf{W}_{\rm{tr}} 
  \leftarrow \mathbf{W}, \mathbf{b}_{\rm{tr}} 
  \leftarrow \mathbf{b}$ \hfill{$\rhd$ save optimal weights}
}
\end{algorithm}
\section{Approximation of the reduced-order matrices based on CSI}
\label{sec:csi}

\subsection{Cubic spline interpolation}

Given  some  interpolation nodes  $a=x_{0}<x_{1}<\cdots<x_{n}=b$,  and
their   corresponding   function  values   $f\left(x_{k}\right)=y_{k},
k=0,1,2 \cdots,  n$. A smooth  function $S(x)$ is  said to be  a cubic
spline  function of  $f(x)$ if  $S(x)$ satisfies  the following  three
conditions
\begin{equation}
\left\{
\begin{array}{l}
S(x) \in \mathbb{C}^{2}[a, b], \\
S\left(x_{k}\right)=f\left(x_{k}\right)=y_{k}, k=0,1,2 \cdots, n, \\
S(x) \text { is a cubic polynomial in }
\left[x_{k}, x_{k+1}\right], k=0,1,2, \cdots, n-1.
\end{array}
\right.
\label{eq:csi_define}
\end{equation}
Let  $S_{k}(x)$   be  the  expression   of  $S(x)$  in   the  interval
$\left[x_{k}, x_{k+1}\right], k=0,1,2 \cdots,  n-1$, which includes $4
n$   unknowns.   Based   on  (\ref{eq:csi_define}),   we   have   $2n$
interpolation conditions
\begin{equation}
S_{k}\left(x_{k}\right)=y_{k}, S_{k}\left(x_{k+1}\right)=y_{k+1}, 
k=0,1,2, \cdots, n-1,
\label{eq:csi_condition1}
\end{equation}
and $2(n-1)$ differential continuity conditions
\begin{equation}
S_{k-1}^{\prime}\left(x_{k}^{-}\right) = 
S_{k}^{\prime}\left(x_{k}^{+}\right), 
S_{k-1}^{\prime \prime}\left(x_{k}^{-}\right) = 
S_{k}^{\prime \prime}\left(x_{k}^{+}\right), k=1,2, \cdots, n-1.
\label{eq:csi_condition2}
\end{equation}
We  cannot  obtain  the  cubic spline  function  $S(x)$  according  to
(\ref{eq:csi_condition1})  and  (\ref{eq:csi_condition2}) because  two
more  conditions  are  required.  In   this  study,  we  consider  the
not-a-knot conditions \cite{Behforooz1995}
\begin{equation}
S_{0}^{\prime \prime \prime}\left(x_{1}^{-}\right) = 
S_{1}^{\prime \prime \prime}\left(x_{1}^{+}\right), 
S_{n-2}^{\prime \prime \prime}\left(x_{n-1}^{-}\right) = 
S_{n-1}^{\prime \prime \prime}\left(x_{n-1}^{+}\right).
\label{eq:csi_condition3}
\end{equation}
Combined               with               the               conditions
(\ref{eq:csi_condition1})-(\ref{eq:csi_condition3}),     the     cubic
polynomial  $S_{k}(x)\ (k=$  $0,1, \cdots,  n-1)$ can  be obtained  by
solving a system  of linear equations of order  $n+1$.  In particular,
we can extend the same  multivariate analysis through the widely used
tensor product formulation  \cite{Hasan2018, Hasan2019}.  To determine
the interpolated value at a  desired point, we use cubic interpolation
of the values at the closest knot points in each respective dimension.
The number of independent variables  in the multivariate CSI method is
not limited. For  more detailed definition of  multivariate CSI method
we refer to \cite{Wang2013}.

\subsection{CSI-based approximation of reduced-order matrix}

\begin{figure}[htbp]
\centering
\includegraphics[width=0.9\textwidth]{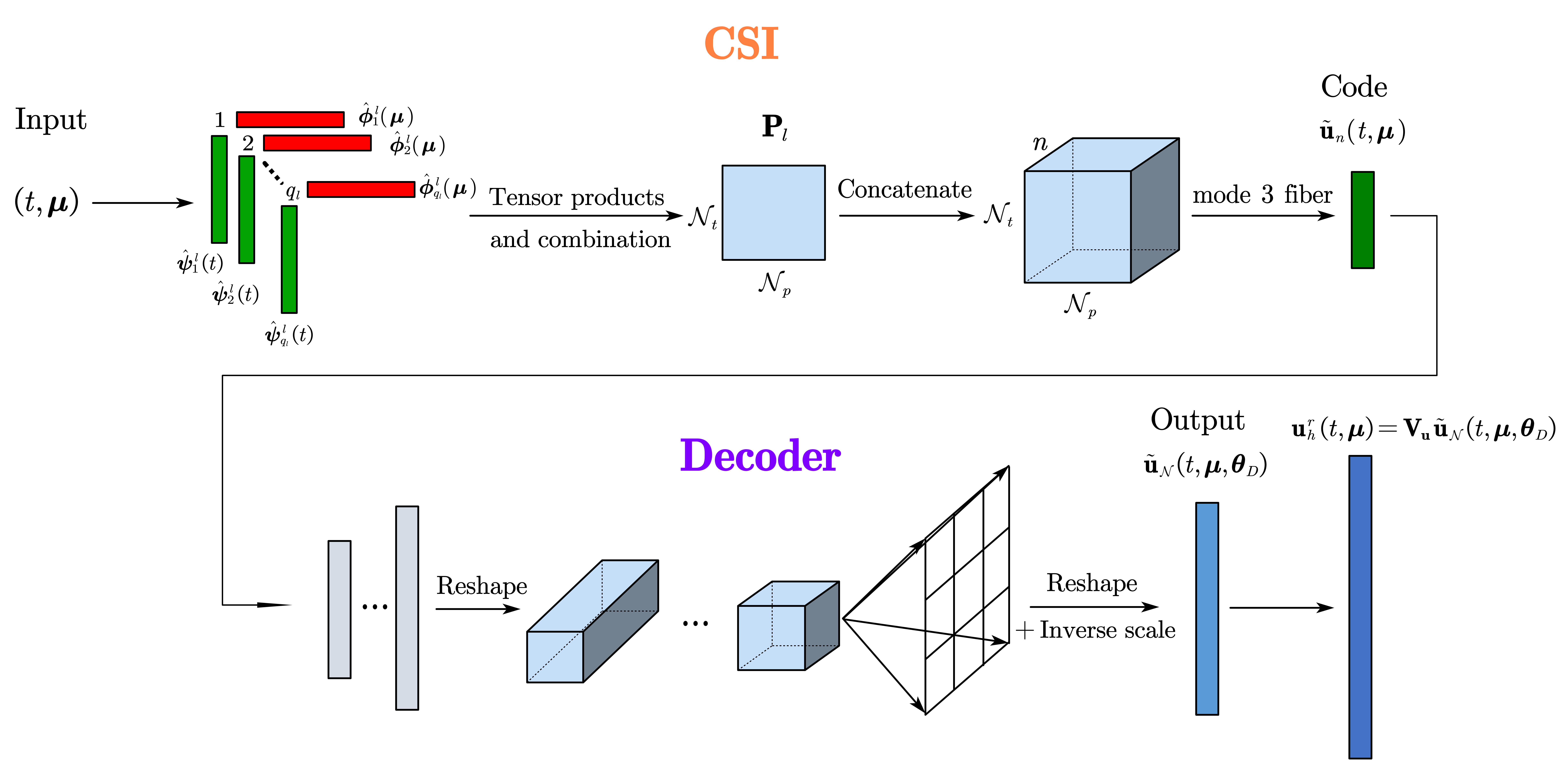}
\caption{Architecture of the  CSI and decoder in  the online stage.The
  CSI-based model produces a low dimensional code
  $\tilde{\mathbf{u}}_n(t, \boldsymbol{\mu})$ for a new input  
  $(t, \boldsymbol{\mu})$, then the decoder reconstructs an approximate
  projection coefficient 
  $\tilde{\mathbf{u}}_{\mathcal{N}}(t, \boldsymbol{\mu}, 
                                     \boldsymbol{\theta}_{D})$ from the code
  $\tilde{\mathbf{u}}_n(t, \boldsymbol{\mu})$.}
\label{fig:online}
\end{figure}

In this  subsection , we  construct a mapping from  time/parameters to
the     low      dimensional     coded      representation     vectors
$\tilde{\mathbf{u}}_{n}(t_i,\boldsymbol{\mu}_j) = 
 [\omega_{1}(t_i, \boldsymbol{\mu}_j), \omega_{2}(t_i, \boldsymbol{\mu}_j), 
  \cdots,
  \omega_{n}(t_i, \boldsymbol{\mu}_j)] \ (i=1,2,\cdots,
 \mathcal{N}_t, j=1,2, \cdots, \mathcal{N}_p)$.  
After the convolutional autoencoder network is trained, we get all the
$\mathcal{N}_s$   low   dimensional   coded   representation   vectors
$\tilde{\mathbf{u}}_{n}(t_i, \boldsymbol{\mu}_j) \in \mathbb{R}^{n}$.
Reshape them to $n$ reduced-order matrices with the form
\begin{equation}
\begin{aligned}
\mathbf{P}_{l} & = 
\left(
\begin{array}{cccc}
\omega_{l}(t_{n_1}, \boldsymbol{\mu}_{1}) & 
\omega_{l}(t_{n_1}, \boldsymbol{\mu}_{2}) & 
\cdots & \omega_{l}(t_{n_1}, \boldsymbol{\mu}_{\mathcal{N}_p}) \\
\omega_{l}(t_{n_2}, \boldsymbol{\mu}_{1}) & 
\omega_{l}(t_{n_2}, \boldsymbol{\mu}_{2}) & 
\cdots & \omega_{l}(t_{n_2}, \boldsymbol{\mu}_{\mathcal{N}_p}) \\
\vdots & \vdots & \ddots & \vdots \\
\omega_{l}(t_{n_{\mathcal{N}_t}}, \boldsymbol{\mu}_{1}) & 
\omega_{l}(t_{n_{\mathcal{N}_t}}, \boldsymbol{\mu}_{2}) & 
\cdots & \omega_{l}(t_{n_{\mathcal{N}_t}}, \boldsymbol{\mu}_{\mathcal{N}_p})
\end{array}
\right) \in \mathbb{R}^{\mathcal{N}_t \times \mathcal{N}_p},
l=1,2, \cdots, n.
\end{aligned}
\end{equation}
It is  difficult to directly  fit the  above matrix in  two dimensions
under  acceptable accuracy  conditions.   So  we resort  to  a SVD  to
decompose  $\mathbf{P}_{l}$  into  several time-  and  parameter-modes
\cite{Guo2019}
\begin{equation}
\mathbf{P}_l \approx \widetilde{\mathbf{P}}_l = 
\sum_{k=1}^{q_l} \sigma_{k}^l 
\boldsymbol{\psi}_{k}^l(\boldsymbol{\phi}_{k}^{l})^{T}, \ 1 \leq l \leq n,
\end{equation}
where $\boldsymbol{\psi}_k^{l}$ and $\boldsymbol{\phi}_k^{l}$ are 
the $k$th discrete time- and parameter- modes for the $l$th reduced-order 
matrix respectively, $\sigma_k^{l}$ is the $k$th singular value, and $q_l$
is the truncation rank corresponding to the error tolerance $\delta$, i.e., 
$q_l = \operatorname{argmin}\left\{\pi(q_l): \pi(q_l) \geq 1-\delta\right\}$
with $\mathcal{\pi}(q_l) = 
      \sum_{k=1}^{q_l}(\sigma_{k}^{l})^{2} / 
      \sum_{k=1}^{r_{l}}(\sigma_k^l)^{2}$ and 
$r_l$ being the rank of $\mathbf{P}_l$. With the discrete modes database, 
CSI models can be trained to approximate the continuous modes as
\begin{equation}
\begin{aligned}
 & t \longmapsto \widehat{\boldsymbol{\psi}}_{ k}^{l}(t), 
   \text { trained from }\{(t_{i},(\boldsymbol{\psi}_{ k}^{l})_{i}), 
   i=1,2, \cdots \mathcal{N}_{t}\}, \\
 & \boldsymbol{\mu} \longmapsto 
   \widehat{\boldsymbol{\phi}}_{ k}^{l}(\boldsymbol{\mu}), 
   \text { trained from }\{(\boldsymbol{\mu}_{j}, 
                           (\boldsymbol{\phi}_{ k}^{l})_{j}), 
           j=1,2, \cdots \mathcal{N}_{p}\}.
\end{aligned}
\end{equation}
Hence, we have
\begin{equation}
\left(\mathbf{P}_{l}\right)_{i j}\approx 
\sum_{k=1}^{q_{l}} \sigma_{k}^{l} 
\widehat{\boldsymbol{\psi}}_{k}^{l}\left(t_{i}\right) 
\widehat{\boldsymbol{\phi}}_{ k}^{l}(\boldsymbol{\mu}_{j}),
\end{equation}
with $1 \leq i \leq \mathcal{N}_{t}, 1 \leq j \leq \mathcal{N}_p$. 
The architecture of the CSI and decoder in the online stage is shown 
in Fig.~\ref{fig:online}. For a new time/parameter value 
$(t^{*}, \boldsymbol{\mu}^{*})$,
we can rapidly get the low-dimensional coded representation
$\tilde{\mathbf{u}}_{n} (t^{*}, \boldsymbol{\mu}^{*}) =
 [\omega_{1}(t^{*}, \boldsymbol{\mu}^{*}),
  \omega_{2}(t^{*}, \boldsymbol{\mu}^{*}),
  \cdots ,
  \omega_{n}(t^{*}, \boldsymbol{\mu}^{*})]$ via CSI.
Then    the    approximation    of    the    projection    coefficient
$\widehat{\alpha}_{\mathbf{u}}(t^{*}, \boldsymbol{\mu}^{*})$ is
obtained   by  the   decoder  network.   The  reduced-order   solution
$\mathbf{u}_{h}^{r}(t, \boldsymbol{\mu})$  served as  an approximation
to the high-fidelity solution $\mathbf{u}_{h}(t,\boldsymbol{\mu})$ can
be recovered by
\begin{equation}
\mathbf{u}_{h}(t^{*}, \boldsymbol{\mu}^{*}) \approx 
\mathbf{u}_{h}^{r}(t^{*}, \boldsymbol{\mu}^{*}) = 
\mathbf{V}_{\mathbf{u}} \widehat{\alpha}_{\mathbf{u}}(t^{*}, 
\boldsymbol{\mu}^{*}), \mathbf{u} \in \{\mathbf{E}, \mathbf{H}\}.
\label{eq:reduced_order_solution}
\end{equation}
\begin{algorithm}[htbp]
\caption{POD-CAE-CSI method for electromagnetic simulations}
\label{algo:pod-cae-csi}
\SetAlgoLined
\newcommand{\forcond}{$i=0$ \KwTo $n$}
\SetKwProg{Fn}{Function}{}{end}
\SetKwFunction{Offline}{POD-CAE-CSI-Offline}
\SetKwFunction{Online}{POD-CAE-CSI-Online}
\Fn{$[\mathbf{V_u}, f_{E}, f_{D}, \widehat{\varpi}]=$\Offline{$\mathcal{P}, 
      \mathcal{T}, \Omega$}}{
Prepare the parameter sampling $\mathcal{P}_{h}^{tr} \subset \mathcal{P} $\;
Calculate the high-fidelity solutions 
$\mathbf{u}_{h}(t_i, \boldsymbol{\mu}_j)$
via the DGTD solver in the time domain 
$\mathcal{T}$, $j= 1,2,\cdots, \mathcal{N}_p$, $\mathbf{u} \in \{\mathbf{E}, 
 \mathbf{H}\}$\;
Form the snapshot matrix 
$\mathbf{S}_\mathbf{u}^{j}\ (\mathbf{u}\in \{\mathbf{E}, \mathbf{H}\})$, and 
prepare	the time sampling $\mathcal{T}_{h}^{tr} \subset \mathcal{T}$\;
Generate the POD basis matrix 
$\mathbf{V_u}\ (\mathbf{u}\in \{\mathbf{E}, \mathbf{H}\})$ via Algorithm 
\ref{algo:2step_pod}\;
Train the autoencoder network 
$f_{E}$ and $f_{D}$ via Algorithm \ref{algo:cae} \;
Build the CSI-based interpolation model $\widehat{\varpi}$.
}
\BlankLine
\Fn{$\mathbf{u}_{h}^{r}(t^{*}, \boldsymbol{\mu}^{*}) = $
    \Online{$\mathbf{V_u}, f_{D}, 
             \widehat{\varpi}, (t^{*}, \boldsymbol{\mu}^{*})$}}{
Compute the approximate low dimensional coded representation 
$\tilde{\mathbf{u}}_{n}(t^{*}, \boldsymbol{\mu}^{*})$ for a new input 
$(t^{*}, \boldsymbol{\mu}^{*})$ through the CSI-based model 
$\widehat{\varpi}$\;
Recover the approximate projection coefficient 
$\widehat{\alpha}_{\mathbf{u}}(t^{*}, \boldsymbol{\mu}^{*})$ through the decoder 
$f_{D}$\;
Evaluate the reduced-order solution 
$\mathbf{u}_{h}^{r}(t^{*}, \boldsymbol{\mu}^{*})$ based on 
\ref{eq:reduced_order_solution}.
}
\end{algorithm}
\section{Numerical experiments}
\label{sec:results}

In this section, numerical  experiments for electromagnetic scattering
problems are  used to evaluate  the proposed CAE-CSI ROM.  We consider
the  2-D  time-domain Maxwell  equations  in  the case  of  transverse
magnetic (TM) waves
\begin{equation}
\left\{
\begin{array}{l}
\nu_{r} \dfrac{\partial H_{x}}{\partial t} + 
\dfrac{\partial E_{z}}{\partial y} = 0, \\ [0.25cm]
\nu_{r} \dfrac{\partial H_{y}}{\partial t} - 
\dfrac{\partial E_{z}}{\partial x} = 0, \\ [0.25cm]
\varepsilon_{r} \dfrac{\partial E_{z}}{\partial t} - 
\dfrac{\partial H_{y}}{\partial x}+\dfrac{\partial H_{x}}{\partial y} = 0.
\end{array}
\right.
\end{equation}
The excitation is an incident plane wave which is defined as
\begin{equation}
\left\{
\begin{array}{l}
H_{x}^{\rm{inc}}(x, y, t) = 0, \\ [0.25cm]
H_{y}^{\rm{inc}}(x, y, t) =-\cos (\omega t-k x), \\ [0.25cm]
E_{z}^{\rm{inc}}(x, y, t) = \cos (\omega t-k x),
\end{array}
\right.
\end{equation}
where $\omega=2  \pi f$  denotes the angular  frequency with  the wave
frequency  $f=30$~GHz, and  $k=\dfrac{\omega}{c}$ is  the wave  number
with $c$ being the wave speed in vacuum.

We  determine the  truncation parameter  $k$ and size parameter  $\mathcal{N}$ in  the
two-step POD method by the following error indicator
\begin{equation}
\bar{e}_{\mathbf{u},\rm{POD}} = 
\dfrac{\sum_{(t,\boldsymbol{\mu}) \in \mathcal{T}_{h}^{tr}\times \mathcal{P}_{h}^{tr}}}
      {\mathcal{N}_s}  
       \dfrac{\left\|\mathbf{u}_{h}(t, \boldsymbol{\mu}) - 
                     \mathbf{V}_{\mathbf{u}} \mathbf{V}_{\mathbf{u}}^{T} 
                     \mathbf{u}_{h}(t, \boldsymbol{\mu})
              \right\|_{\mathbb{R}^{\mathcal{N}_h}}}
      {\left\|\mathbf{u}_{h}(t, \boldsymbol{\mu})
       \right\|_{\mathbb{R}^{\mathcal{N}_h}}}.
\end{equation}
The  relative error  between the  reduced-order solution  generated by
CAE-CSI  and the  high-fidelity  solution  is used  as  the metric  to
evaluate the accuracy of the results
\begin{equation}
e_{\mathbf{u}, \rm{CAE-CSI}}(t, \boldsymbol{\mu}) = 
\dfrac{\left\|\mathbf{u}_{h}(t, \boldsymbol{\mu}) - 
              \mathbf{u}_{h}^{r}(t, \boldsymbol{\mu})
       \right\|_{\mathbb{R}^{\mathcal{N}_h}}}
      {\left\|\mathbf{u}_{h}(t, \boldsymbol{\mu})
       \right\|_{\mathbb{R}^{\mathcal{N}_h}}} = 
\dfrac{\left\|\mathbf{u}_{h}(t, \boldsymbol{\mu}) - 
              \mathbf{V}_{\mathbf{u}} 
              \widehat{\alpha}_{\mathbf{u}}(t, \boldsymbol{\mu})
       \right\|_{\mathbb{R}^{\mathcal{N}_h}}}
      {\left\|\mathbf{u}_{h}(t, \boldsymbol{\mu})
       \right\|_{\mathbb{R}^{\mathcal{N}_h}}}, 
\mathbf{u} \in \{\mathbf{E},\mathbf{H}\},
\end{equation}
which will be compared with the relative projection error
\begin{equation}
e_{\mathbf{u}, \rm{Pro}}(t, \boldsymbol{\mu}) = 
\dfrac{\left\|\mathbf{u}_{h}(t, \boldsymbol{\mu}) - 
              \mathbf{u}_{h}^{p}(t, \boldsymbol{\mu})
       \right\|_{\mathbb{R}^{\mathcal{N}_h}}}
      {\left\|\mathbf{u}_{h}(t, \boldsymbol{\mu})
       \right\|_{\mathbb{R}^{\mathcal{N}_h}}} = 
\dfrac{\left\|\mathbf{u}_{h}(t, \boldsymbol{\mu}) - 
              \mathbf{V}_{\mathbf{u}} \mathbf{V}_{\mathbf{u}}^{T} 
              \mathbf{u}_{h}(t, \boldsymbol{\mu})
       \right\|_{\mathbb{R}^{\mathcal{N}_h}}}
      {\left\|\mathbf{u}_{h}(t, \boldsymbol{\mu})
       \right\|_{\mathbb{R}^{\mathcal{N}_h}}}, 
\mathbf{u} \in \{\mathbf{E}, \mathbf{H}\}.
\end{equation}
The above  errors are evaluated  on a testing  time/parameter sampling
$\mathcal{T}_{h}^{te} \times \mathcal{P}_{h}^{te}$  of   size
$\mathcal{N}_{te}$.  Furthermore,  the average  relative errors  are
used to measure the accuracy of  the ROM in our numerical experiments,
which are defined as
\begin{equation}
\bar{e}_{\mathbf{u},\rm{CAE-CSI}} = 
\dfrac{\sum_{(t,\boldsymbol{\mu}) \in \mathcal{T}_{h}^{te}\times \mathcal{P}_{h}^{te}} 
       e_{\mathbf{u},\rm{CAE-CSI}}(t, \boldsymbol{\mu})}
      {\mathcal{N}_{te}},
\quad 
\bar{e}_{\mathbf{u},\rm{Pro}} = 
\dfrac{\sum_{(t,\boldsymbol{\mu}) \in \mathcal{T}_{h}^{te}\times \mathcal{P}_{h}^{te}} 
       e_{\mathbf{u},\rm{Pro}}(t, \boldsymbol{\mu})}
      {\mathcal{N}_{te}}.
\end{equation}
The CSI-based  approximation of  the reduced-order matrices  are built
via      the      SciPy       functions      \rm{CubicSpline}      and
\rm{RegularGridInterpolator}  in   the  first  and   second  numerical
experiments  respectively.  The  DGTD  and two-step  POD  methods  are
implemented  in MATLAB,  while CAE-CSI  is developed  in Python.  All
simulations  are run  on a  computer equipped  with an Intel Core  i9
10-core 2.8 GHz $\times$ 20 CPU with 64 GB RAM.

\subsection{Scattering of a plane wave by a dielectric disk}

The first numerical experiment is  the electromagnetic scattering of a
plane  wave  by  a  dielectric  disk.   The  computational  domain  is
artificially  bounded  by  a  square  $\Omega=[-2.6  \mathrm{~m},  2.6
  \mathrm{~m}]  \times[-2.6 \mathrm{~m},  2.6  \mathrm{~m}]$ where  we
impose the  Silver-Müller ABC  boundary condition.   The range  of the
relative  permittivity is  $\varepsilon_{r}  \in [1.0  , 5.0]$  (i.e.,
$\mathcal{P} = [1.0 , 5.0]$)  with the relative permeability $\nu_{r} =
1.0$, i.e., we consider a nonmagnetic material.  The medium outside to
the dielectric  disk is assumed  to be vacuum, i.e.   $\varepsilon_r =
\nu_r  = 1.0$.   The  high-fidelity simulations  are  performed on  an
unstructured triangular  mesh with  2575 nodes  and 5044  elements, in
which 1092  elements are located inside  the disk. And we use a DGTD  method with 
  $\mathbb{P}_2$ approximation, 
  resulting in $\mathcal{N}_{h} =
  30264$ DOFs for the FOM.

During  the offline  stage,  the DGTD  solver is  used  to generate  a
collection   of   high-fidelity  solutions   at   $\mathcal{N}_{p}=81$
equidistant  parameter sampling  points  (i.e., $\boldsymbol{\mu}  \in
\mathcal{P}_{h}^{tr} = \{1.0,1.05, \cdots, 4.95, 5.0\}$) with
$\mathcal{N}_{t}=263$ points in the  last oscillation period (i.e., $t
\in \mathcal{T}_{h}^{tr}=$  $\{49.0024, 49.006, \cdots, 49.9623, 49.966\}$).   
A test  parameter set  
$\mathcal{P}_{h}^{te} =  \{1.215, 2.215, 3.215, 4.215\}$ and  
test  time  set $\mathcal{T}_{h}^{te} = \mathcal{T}_{h}^{tr}$ 
are used to evaluate the proposed method.

\begin{figure}[htbp]
\centering
\includegraphics[scale=0.575]{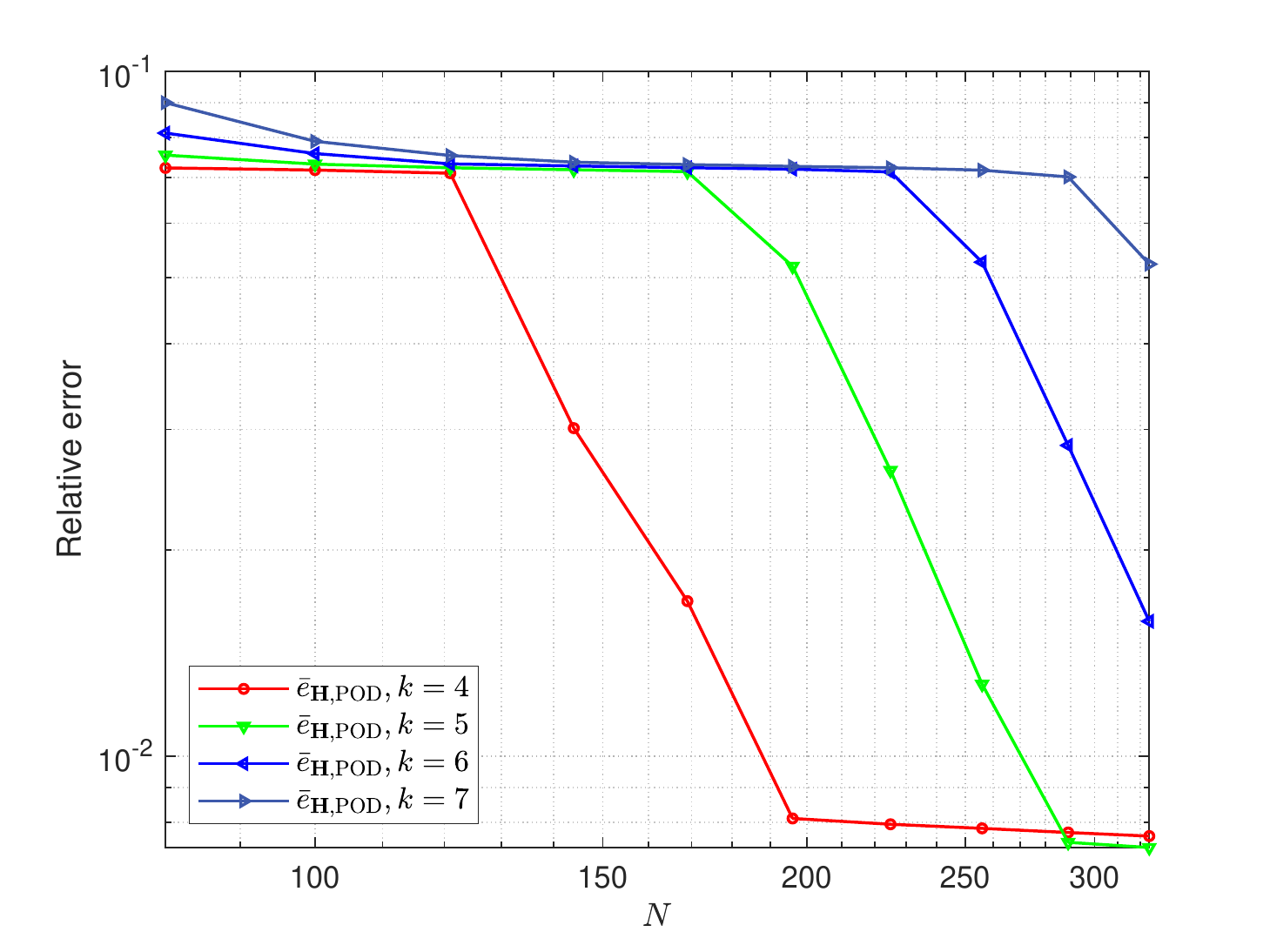}
\hspace{-1.0cm}
\includegraphics[scale=0.575]{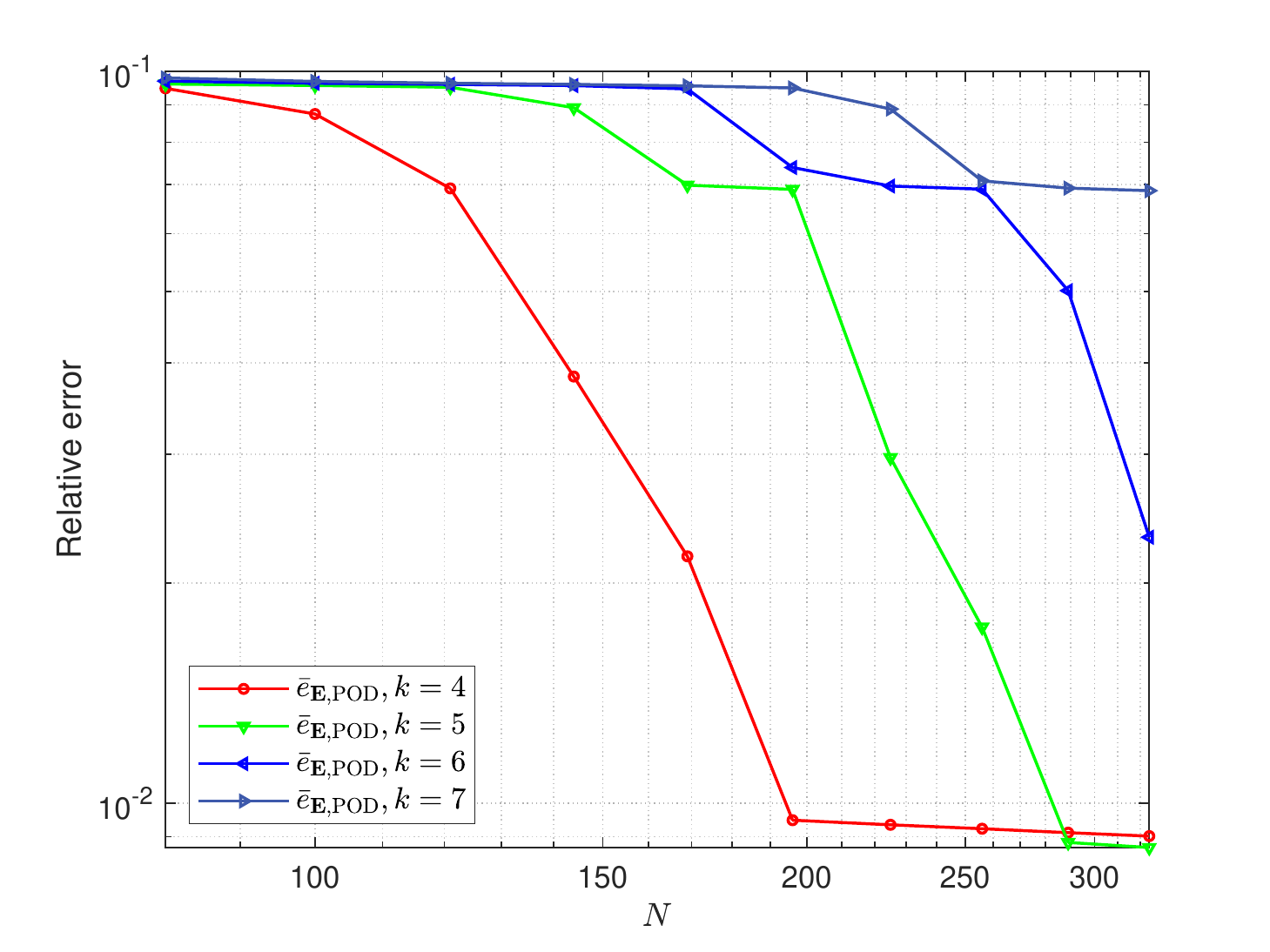}
\caption{Scattering  of  a  plane  wave  by  a  dielectric  disk. 
  Convergence histories of  $\bar{e}_{\mathbf{H},\rm{POD}}$ (left) and
  $\bar{e}_{\mathbf{E},\rm{POD}}$  (right)  with different  truncation
  parameters $k$ and $\mathcal{N}$.}
\label{fig:convergence}
\end{figure}

\begin{figure}[htbp]
\centering
\includegraphics[scale=0.575]{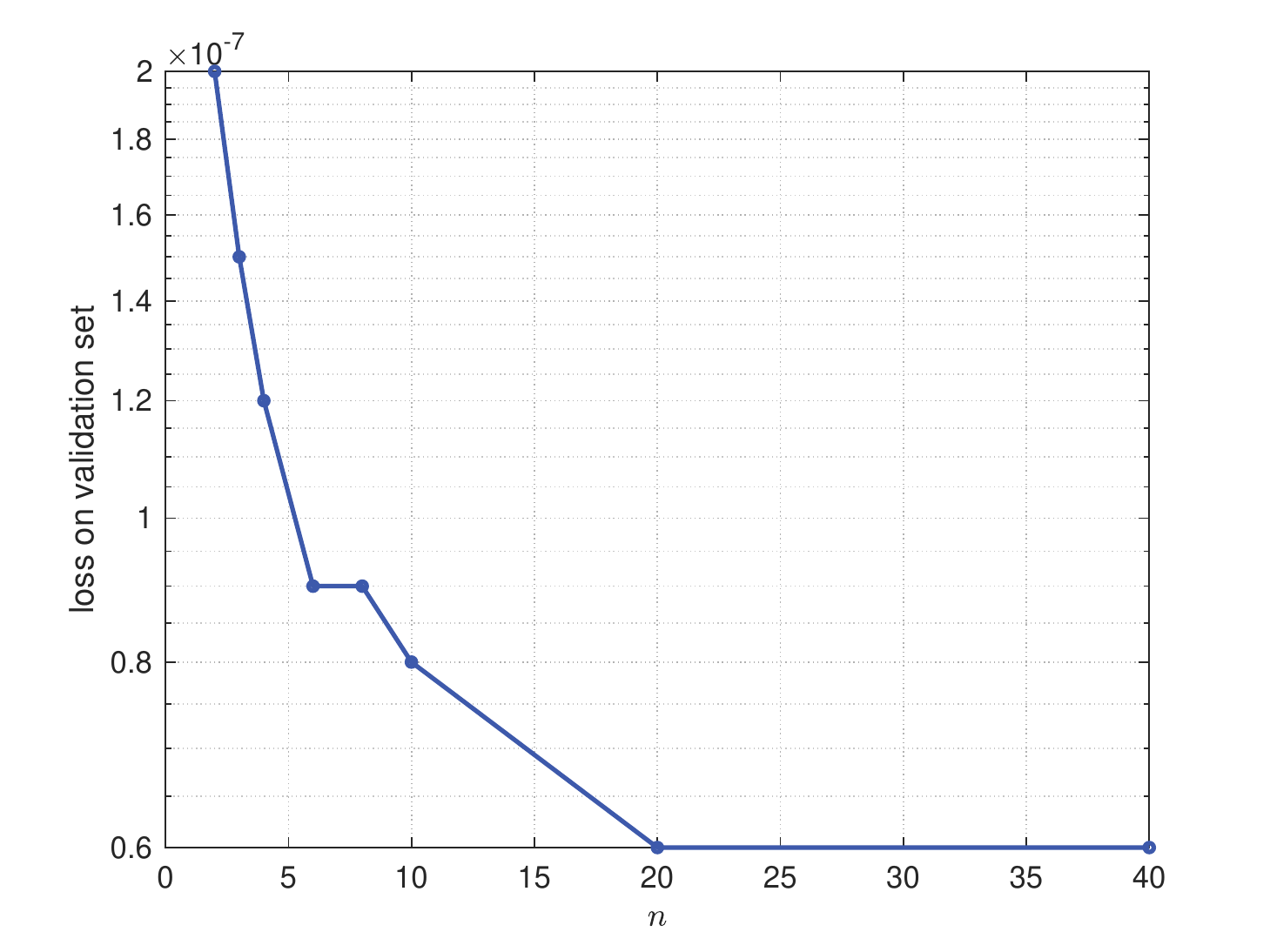}
\hspace{-1.0cm}
\includegraphics[scale=0.575]{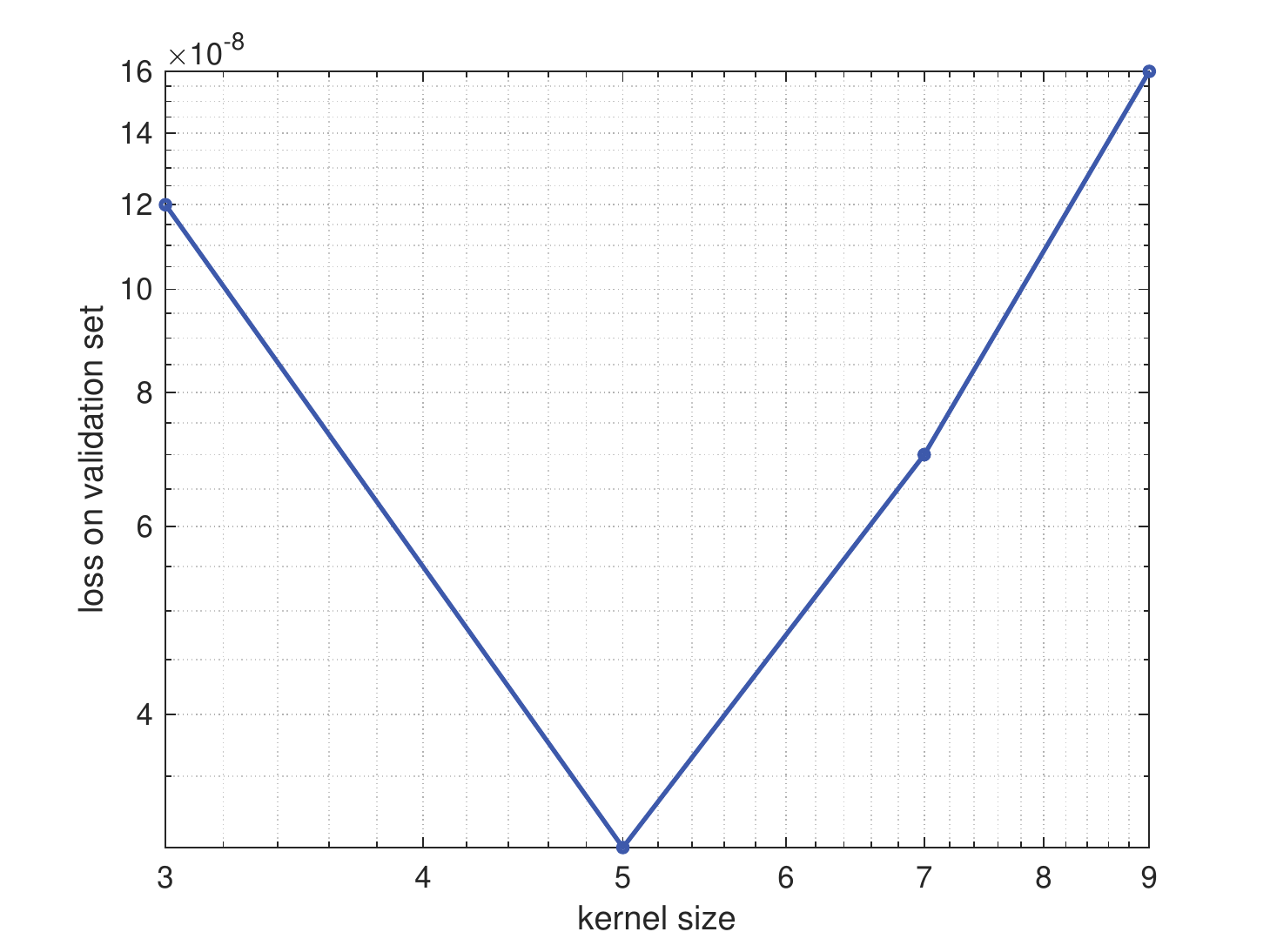}
\caption{Left: loss on validation set versus $n$.
         Right: loss on validation set versus kernel size.}
\label{fig:loss_vs_n_and_kersize}
\end{figure}

\begin{figure}[htbp]
\centering
\includegraphics[width=1\textwidth]
            {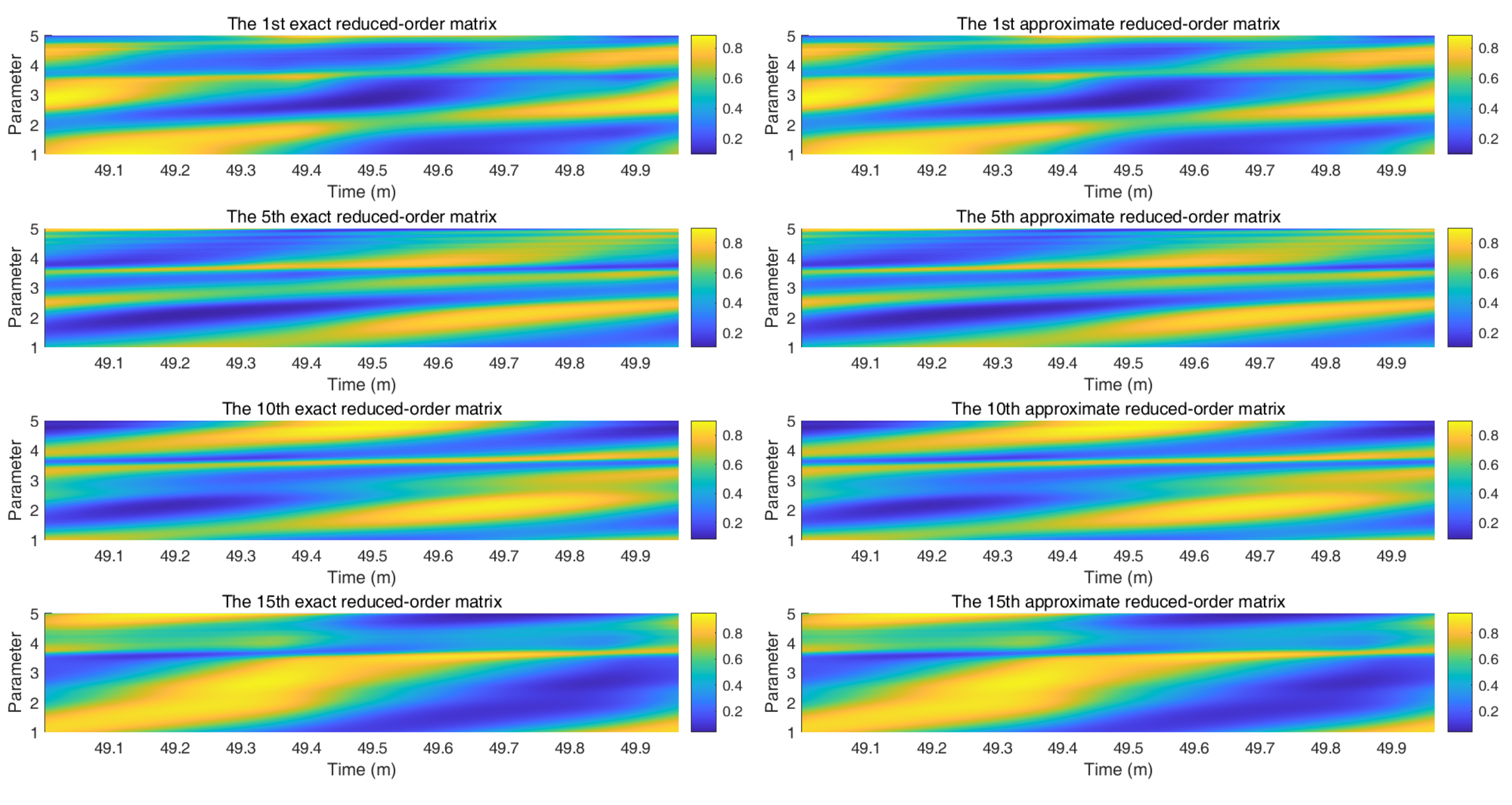}
\caption{Scattering of a plane wave by a dielectric disk. The 1-st, 
         5-th, 10-th, and 15-th exact and approximate reduced-order
         matrices based on CSI.}
\label{fig:reduced_order_matrix}
\end{figure}

\begin{figure}[htbp]
\centering
\begin{minipage}[t]{0.485\linewidth}
\centering
\includegraphics[scale=0.475]{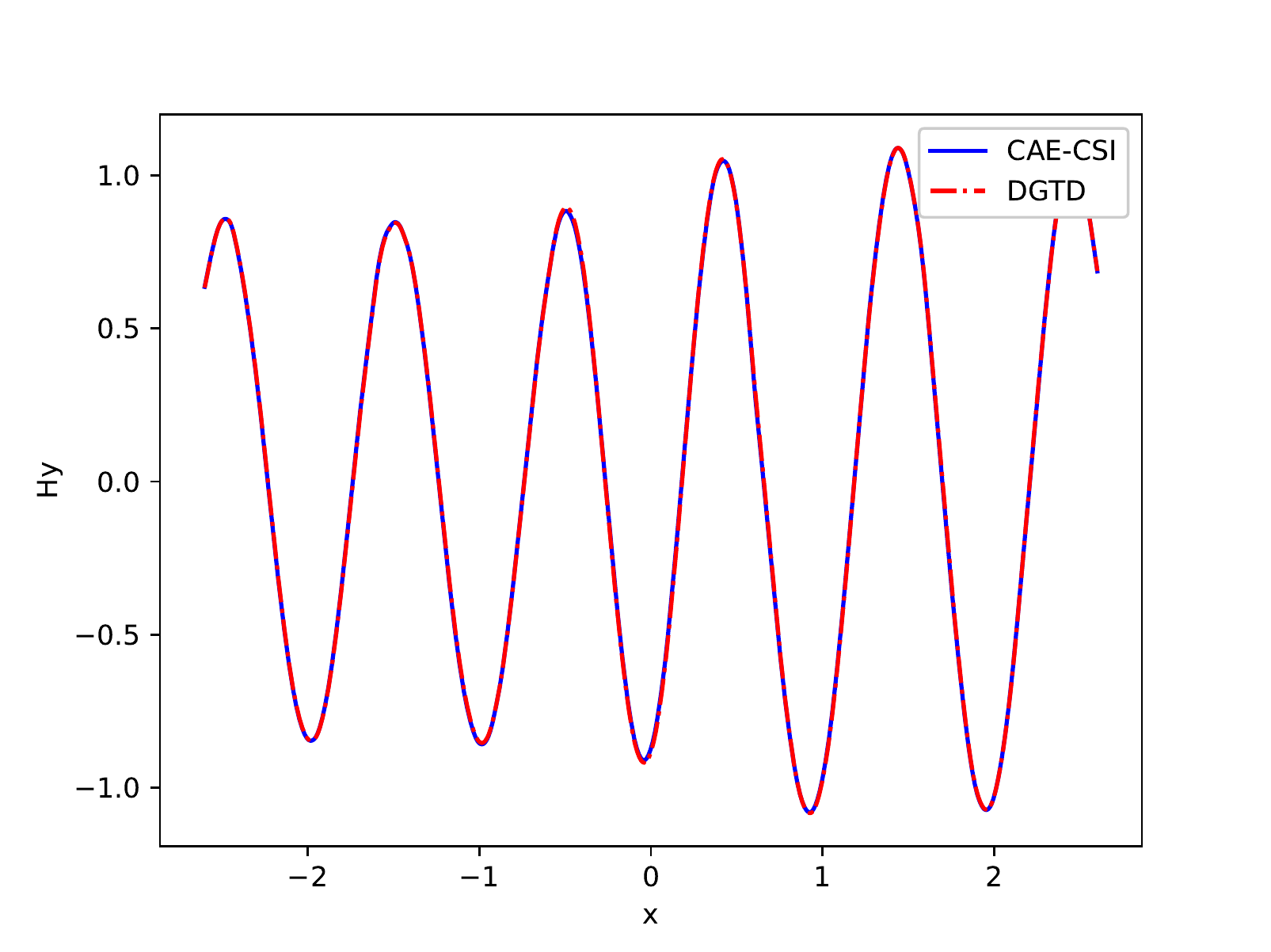}
\end{minipage}
\hspace{-0.5cm}
\begin{minipage}[t]{0.485\linewidth}
\centering
\includegraphics[scale=0.475]{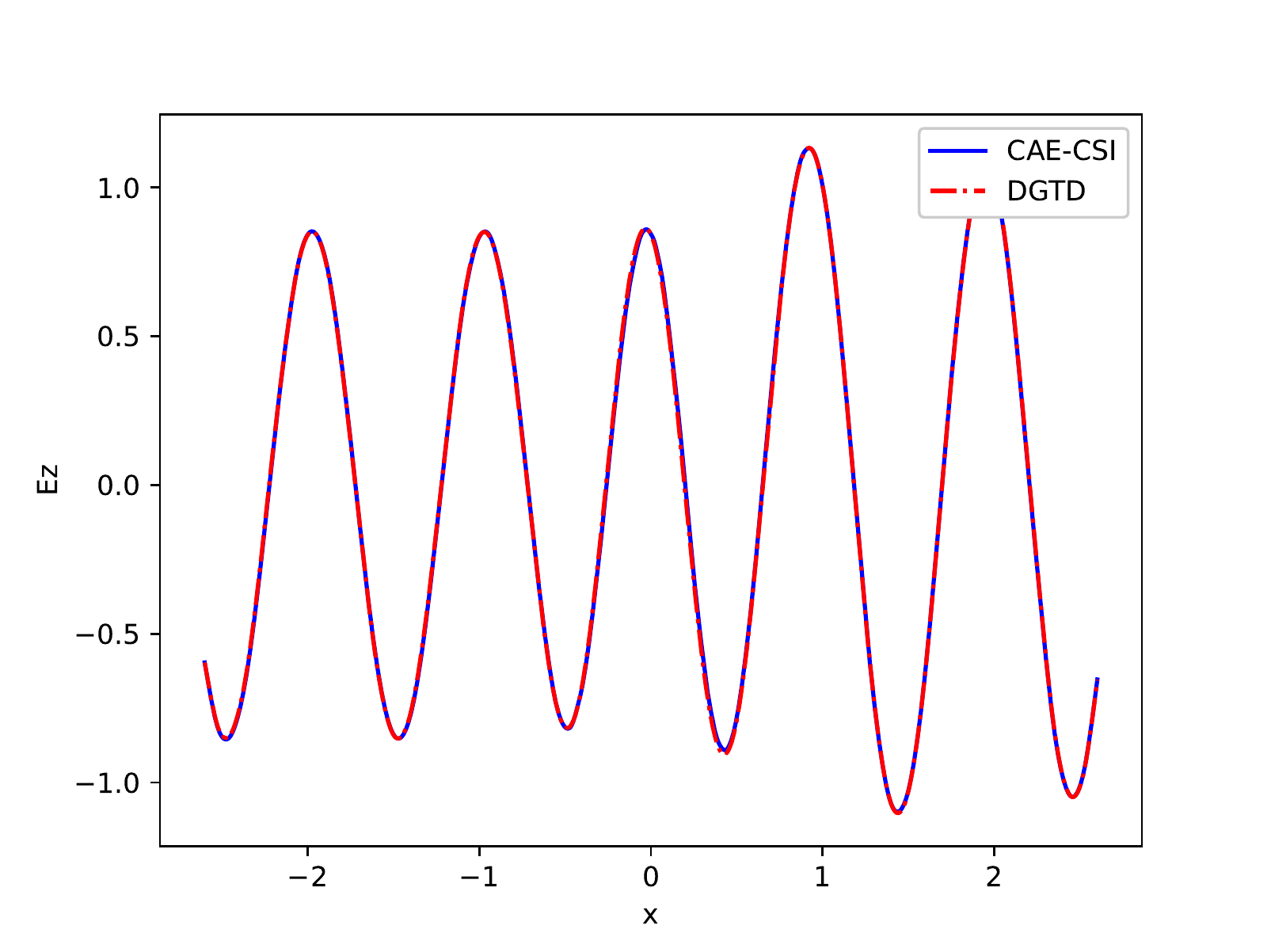}
\end{minipage} \\
\begin{minipage}[t]{0.485\linewidth}
\centering
\includegraphics[scale=0.475]{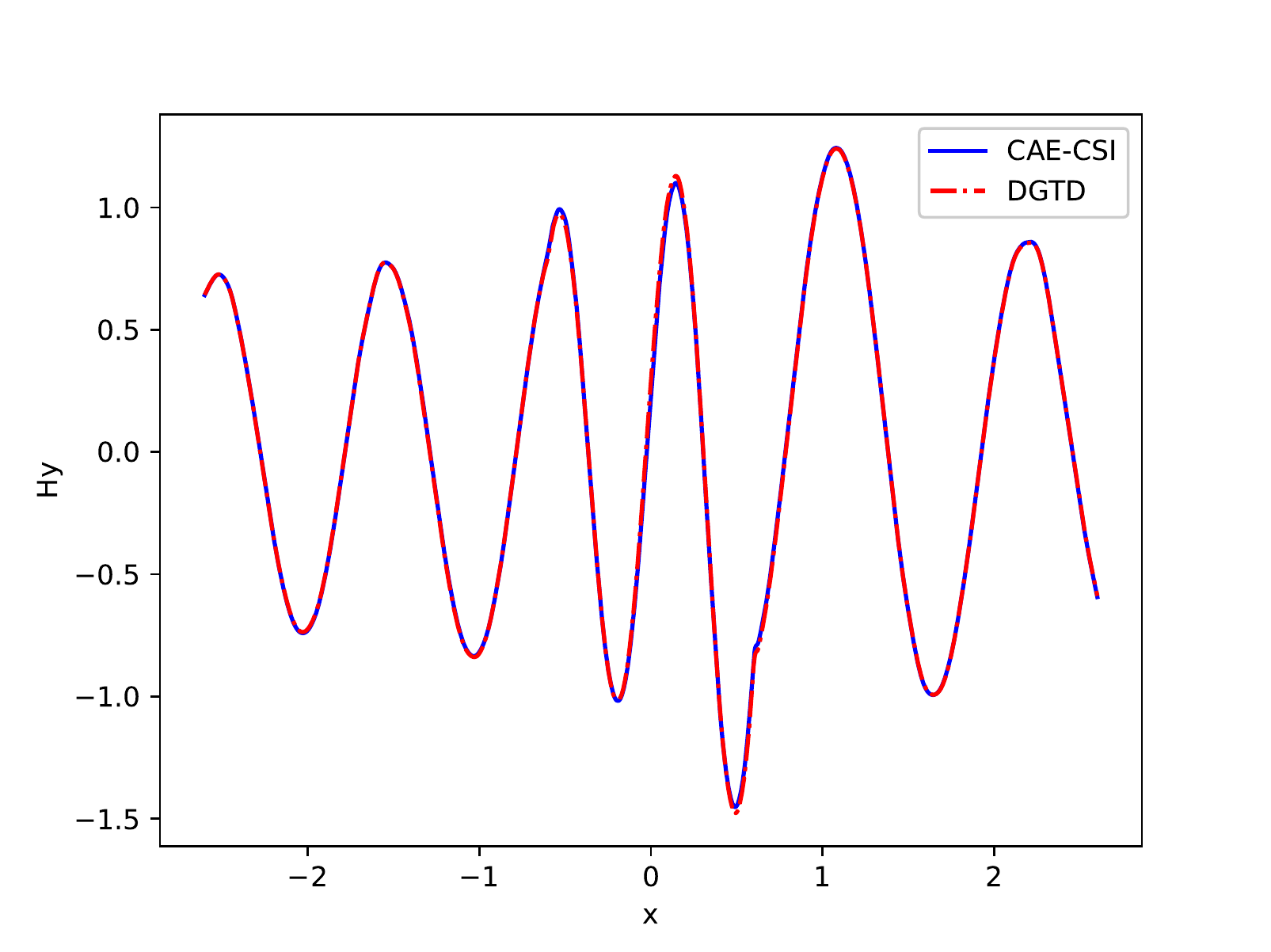}
\end{minipage}
\hspace{-0.5cm}
\begin{minipage}[t]{0.485\linewidth}
\centering
\includegraphics[scale=0.475]{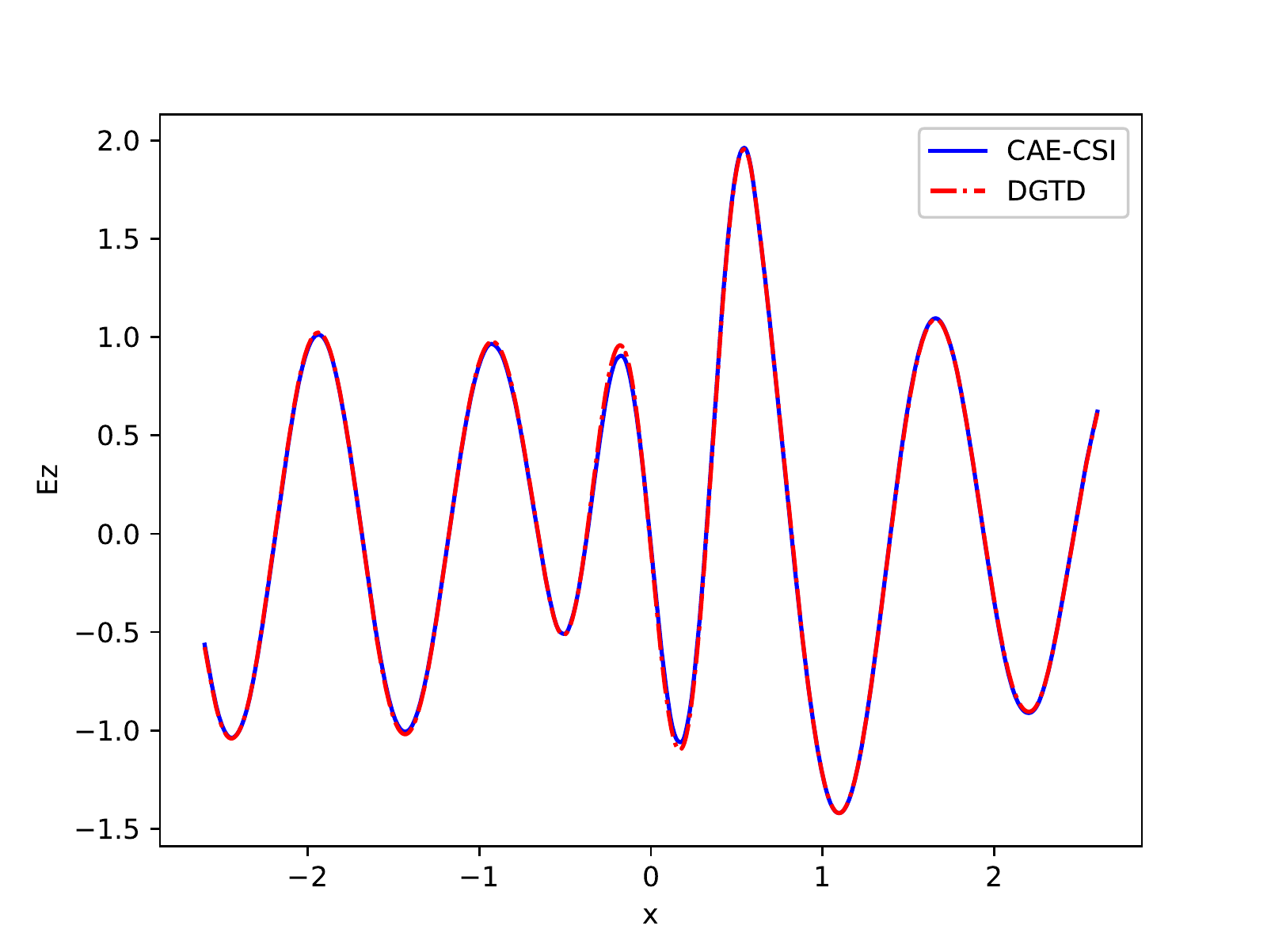}
\end{minipage}
\begin{minipage}[t]{0.485\linewidth}
\centering
\includegraphics[scale=0.475]{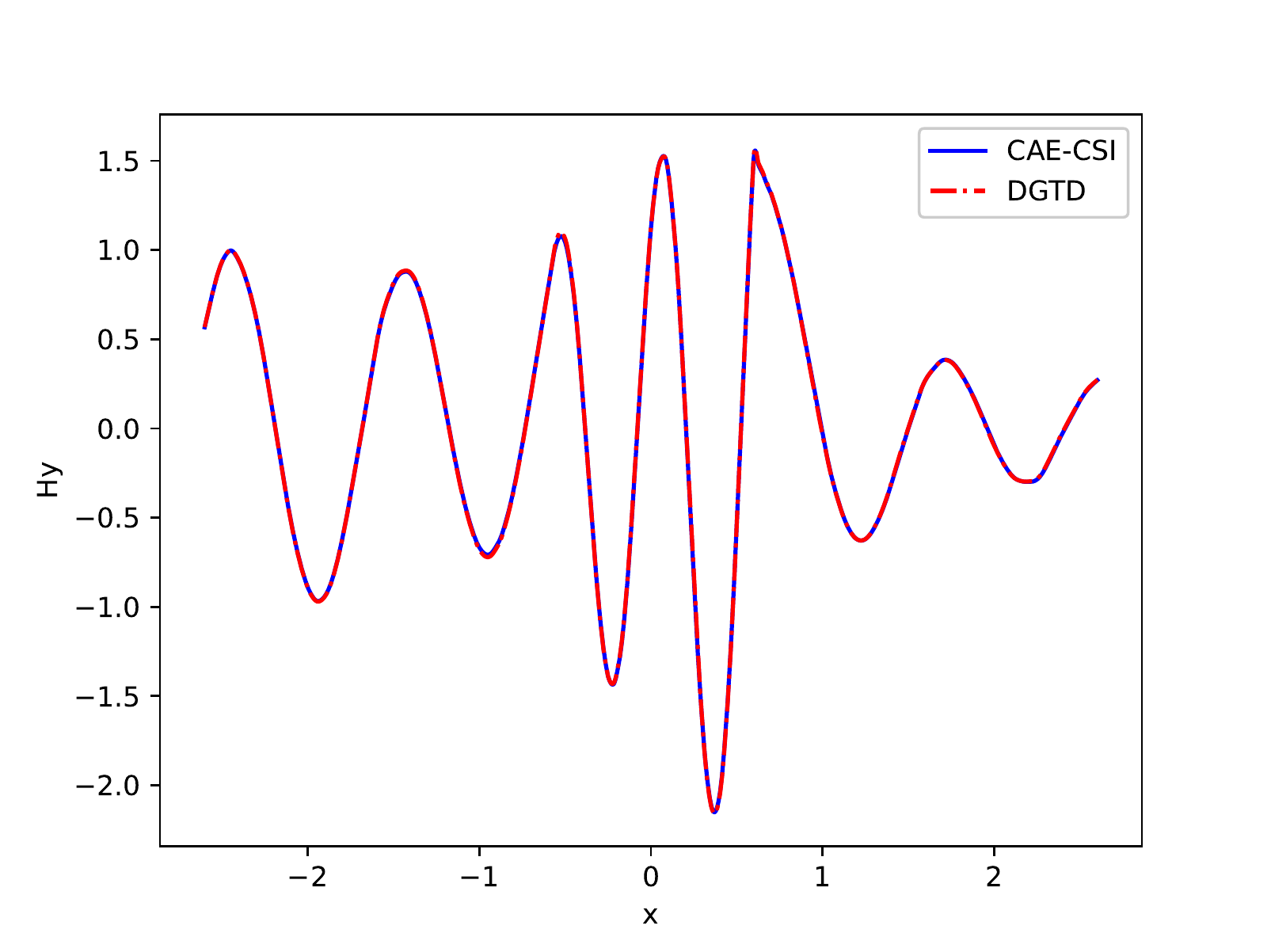}
\end{minipage}
\hspace{-0.5cm}
\begin{minipage}[t]{0.485\linewidth}
\centering
\includegraphics[scale=0.475]{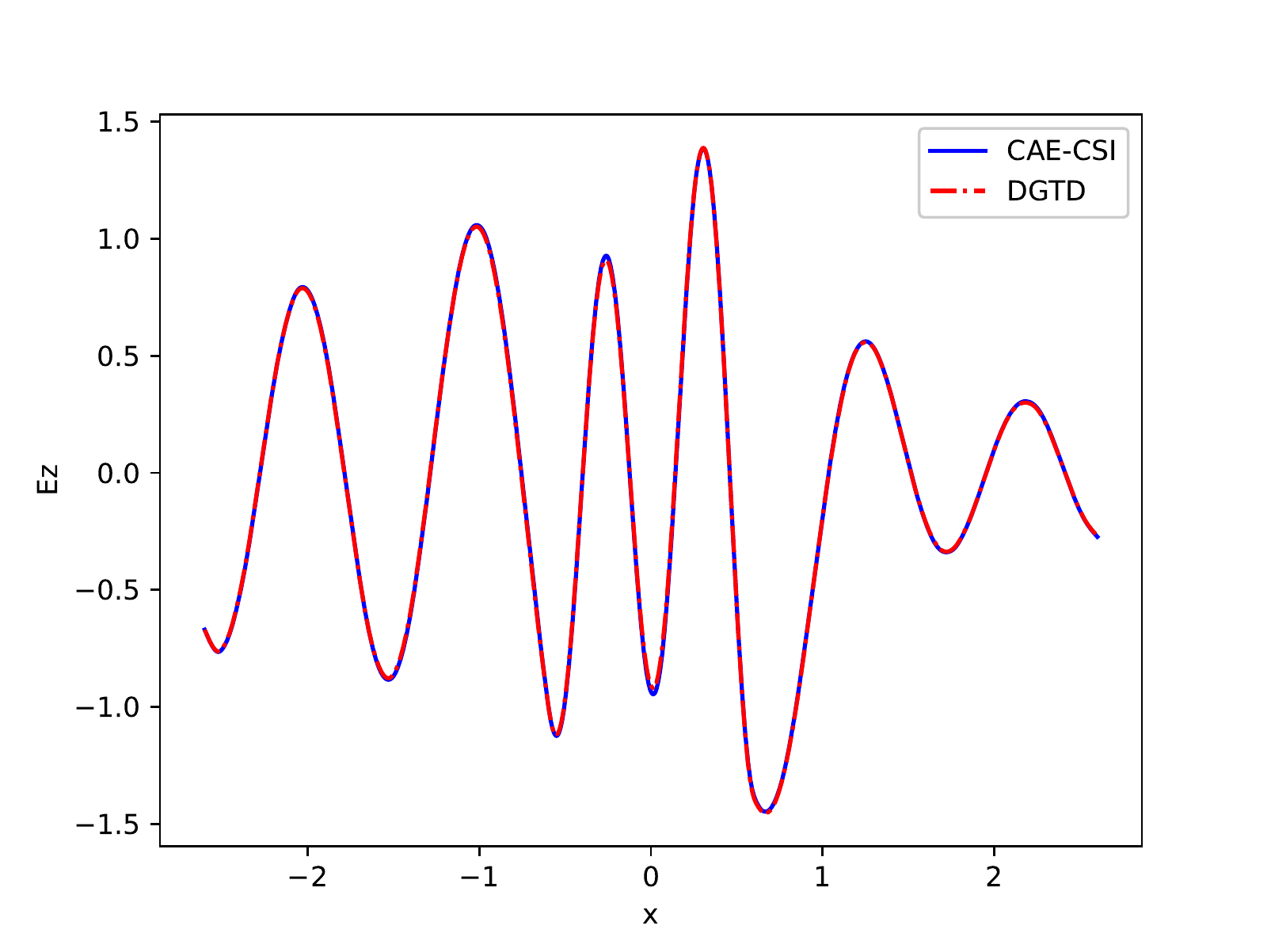}
\end{minipage} \\
\begin{minipage}[t]{0.485\linewidth}
\centering
\includegraphics[scale=0.475]{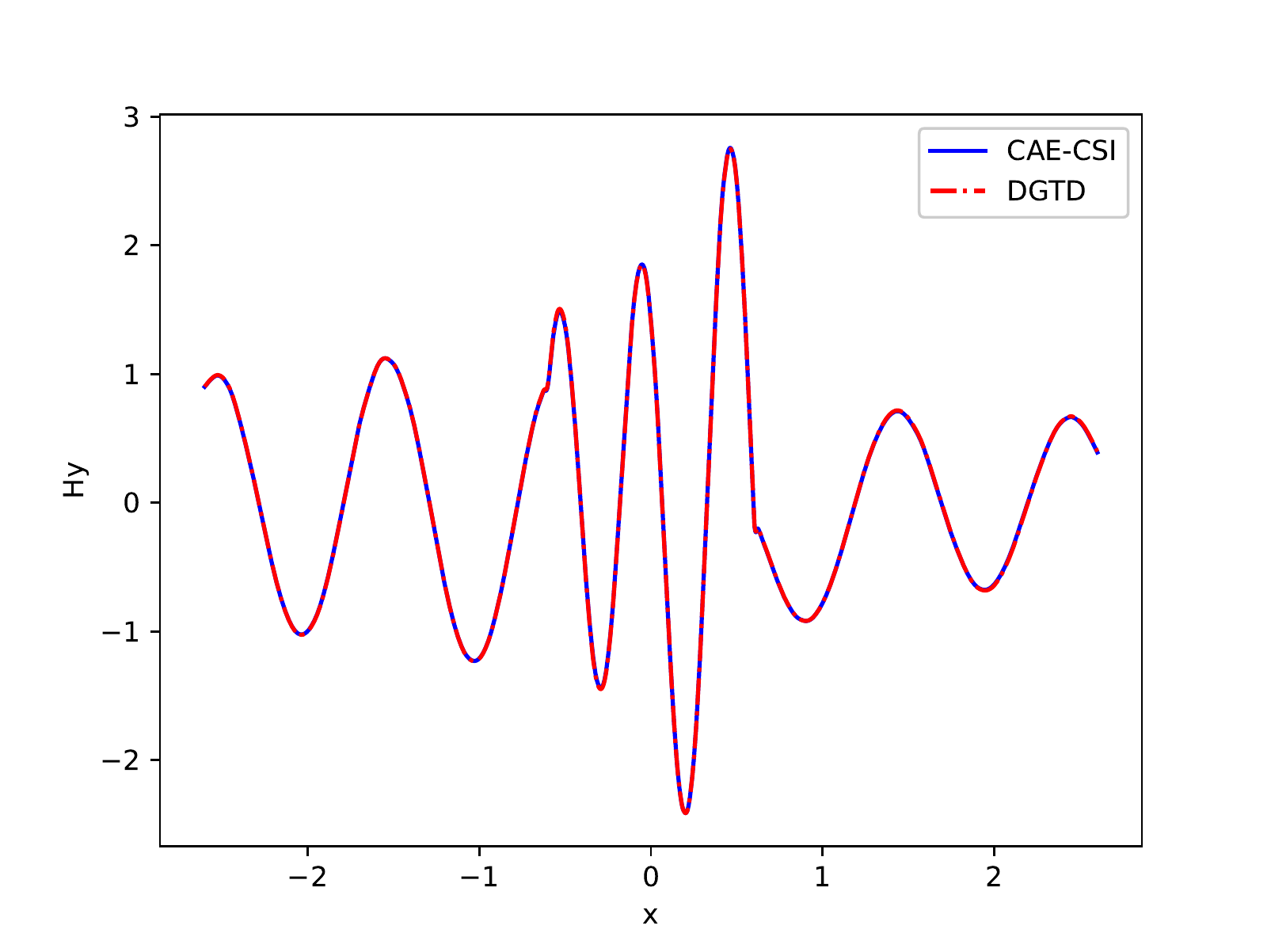}
\end{minipage}
\hspace{-0.5cm}
\begin{minipage}[t]{0.485\linewidth}
\centering
\includegraphics[scale=0.475]{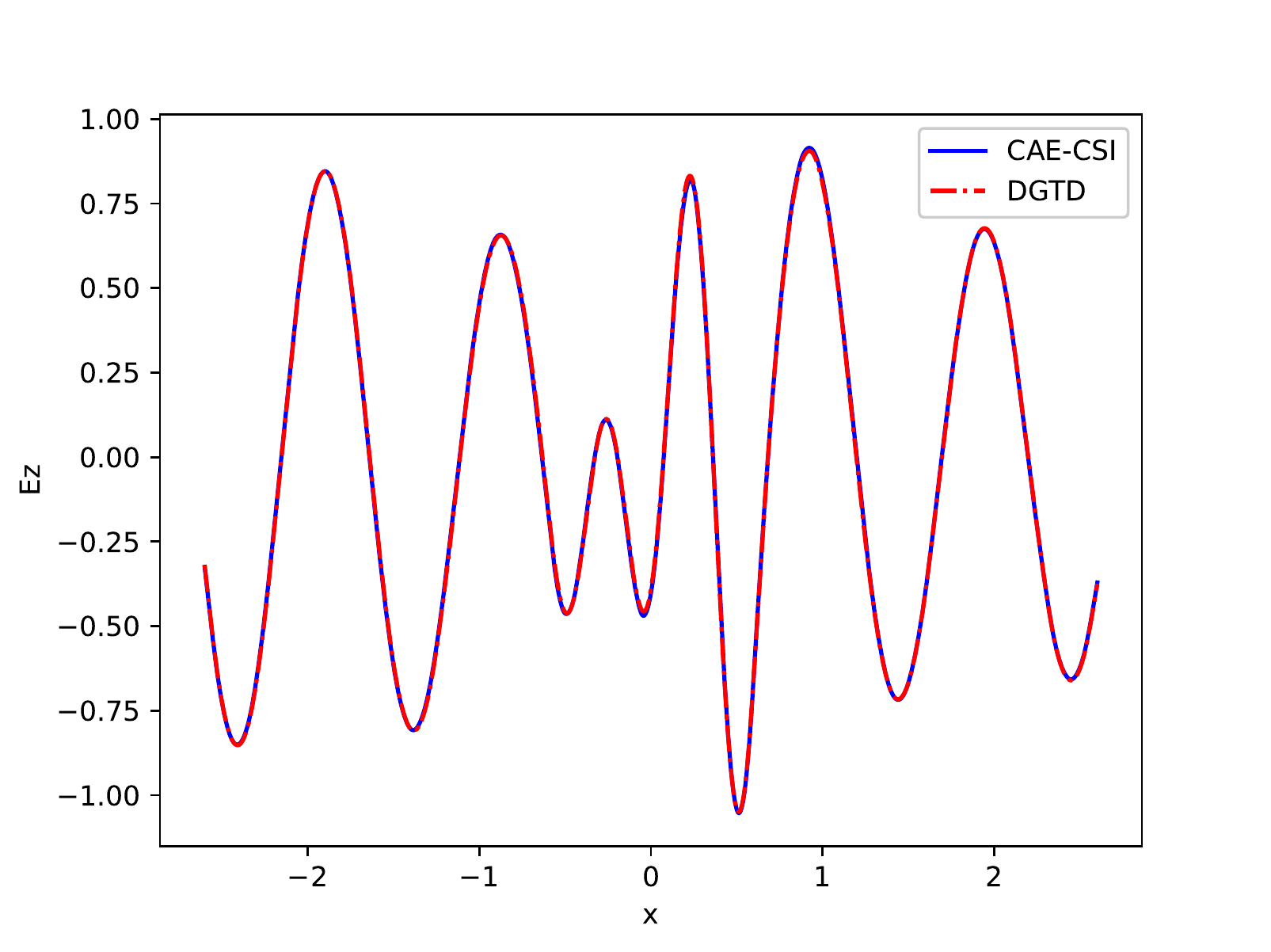}
\end{minipage}
\caption{Scattering of a plane wave by a dielectric disk.
         Comparison of the 1-D $x$-wise distribution along $y=0$ of the 
         real part of $H_y$ (left) and $E_{z}$ (right) for the
	 testing parameter instances: 
 	 $\varepsilon_1=1.215$ (1st row), $\varepsilon_{2}=2.215$ (2th row), 
         $\varepsilon_3=3.215$ (3th row), $\varepsilon_{4}=4.215$ (4th row).}
\label{fig:x_fields}
\end{figure}

\begin{figure}[htbp]
\centering
\begin{minipage}[t]{0.95\textwidth}
\centering
\includegraphics[scale=0.7]{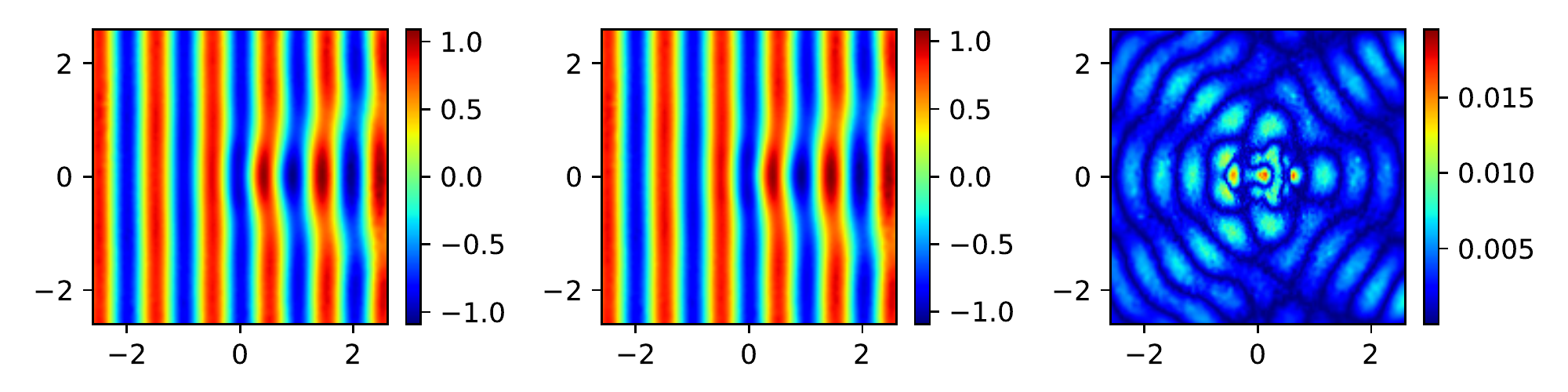}
\end{minipage} \\
\begin{minipage}[t]{0.95\textwidth}
\centering
\includegraphics[scale=0.7]{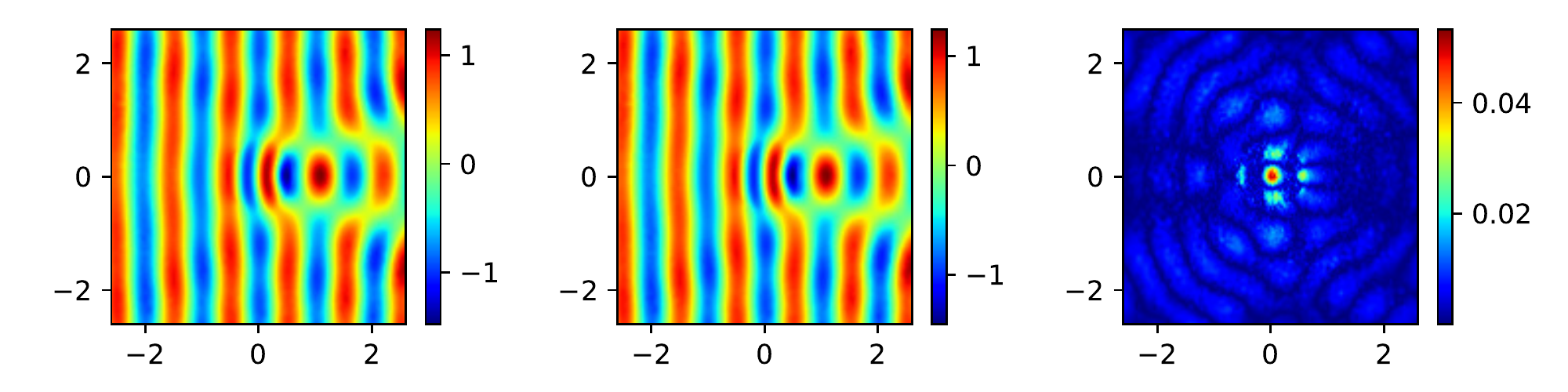}
\end{minipage} \\
\begin{minipage}[t]{0.95\textwidth}
\centering
\includegraphics[scale=0.7]{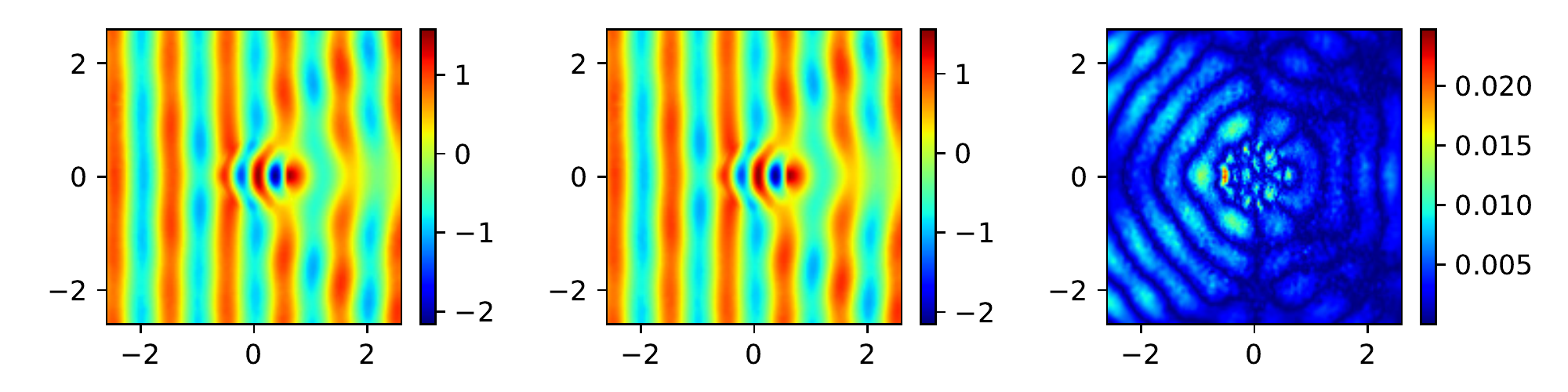}
\end{minipage} \\
\begin{minipage}[t]{0.95\textwidth}
\centering
\includegraphics[scale=0.7]{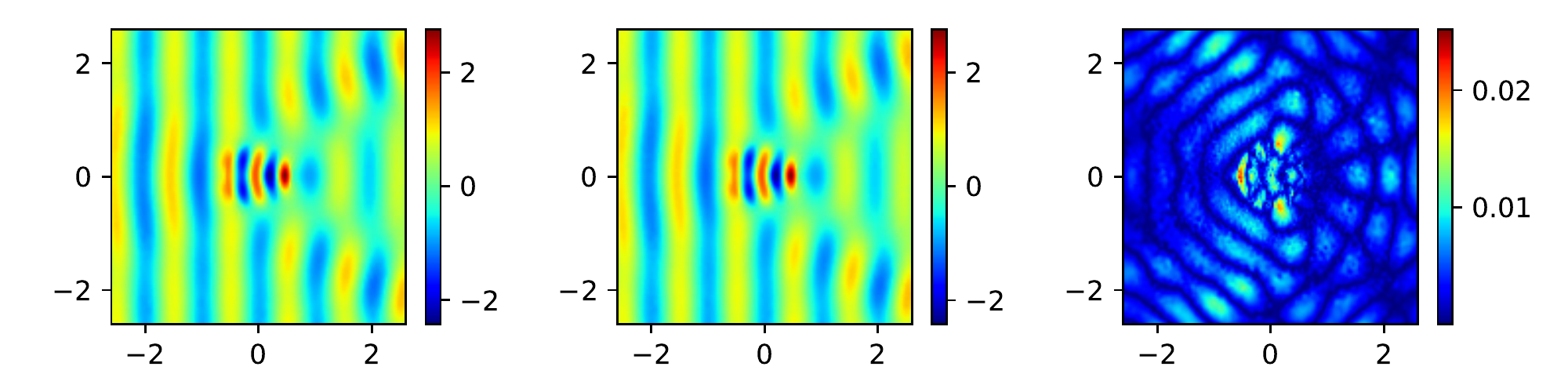}
\end{minipage}
\caption{Scattering of a plane wave by a dielectric disk.
	 Comparison of the 2-D distribution of the real part of $H_{y}$
 	 between DGTD (left), CAE-CSI (middle) and relative error (right) 
         for the testing parameter instances:
 	 $\varepsilon_1=1.215$ (1st row), $\varepsilon_{2}=2.215$ (2th row), 
         $\varepsilon_3=3.215$ (3th row), $\varepsilon_{4}=4.215$ (4th row).}
\label{fig:xy_fields_hy}
\end{figure}

\begin{figure}[htbp]
\centering
\begin{minipage}[t]{0.95\textwidth}
\centering
\includegraphics[scale=0.7]{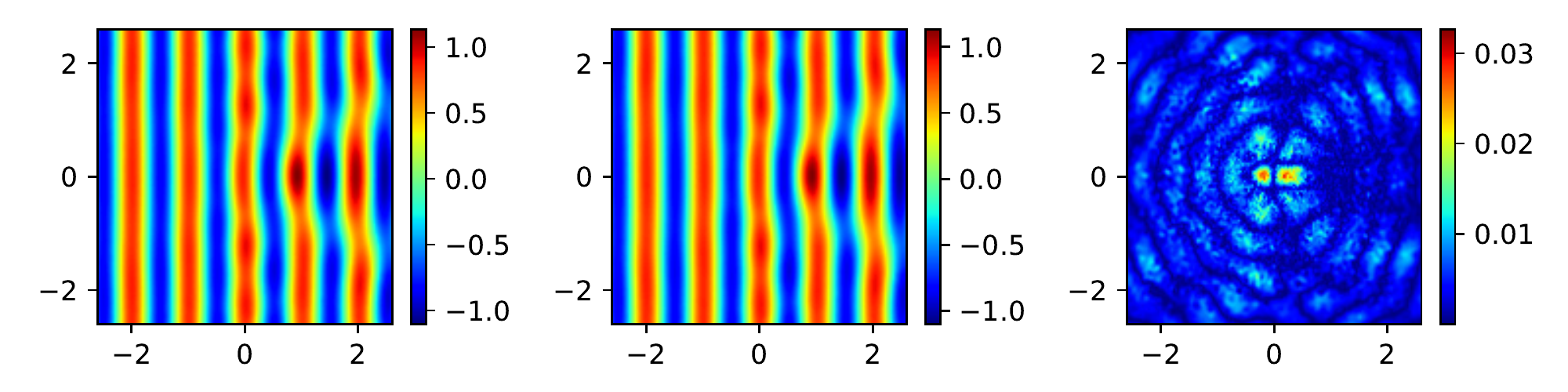}
\end{minipage} \\
\begin{minipage}[t]{0.95\textwidth}
\centering
\includegraphics[scale=0.7]{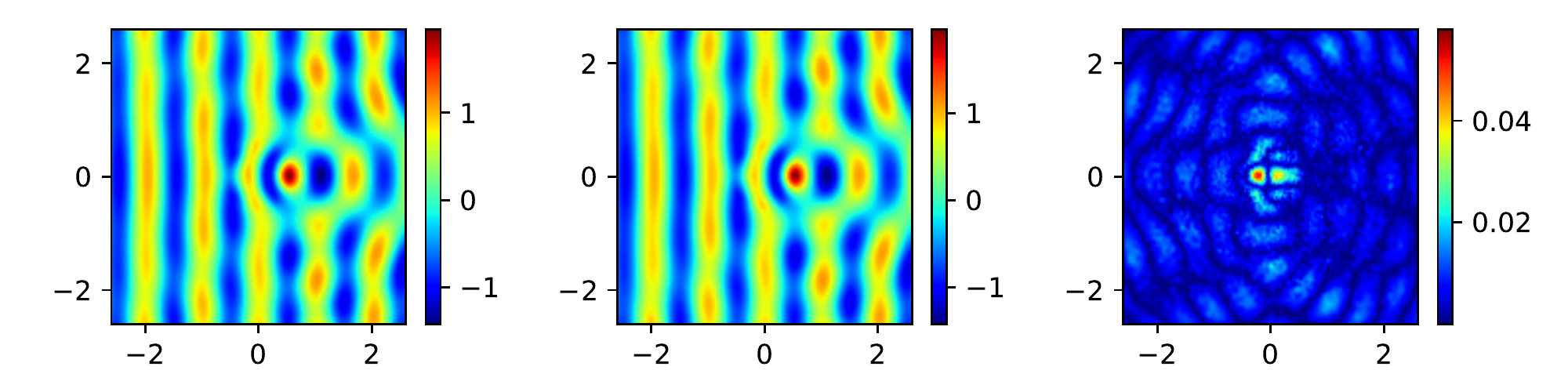}
\end{minipage} \\
\begin{minipage}[t]{0.95\textwidth}
\centering
\includegraphics[scale=0.7]{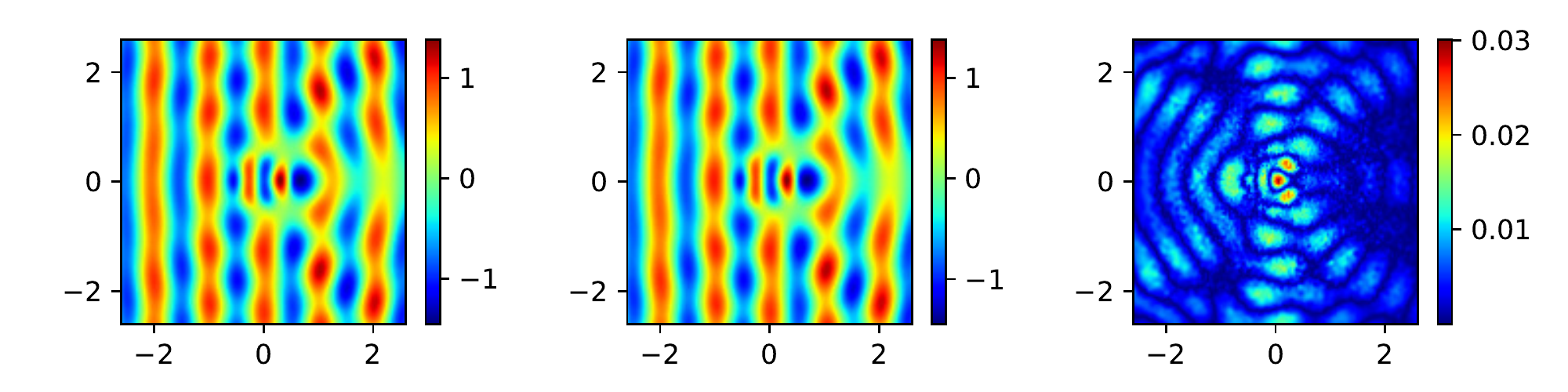}
\end{minipage} \\
\begin{minipage}[t]{0.95\textwidth}
\centering
\includegraphics[scale=0.7]{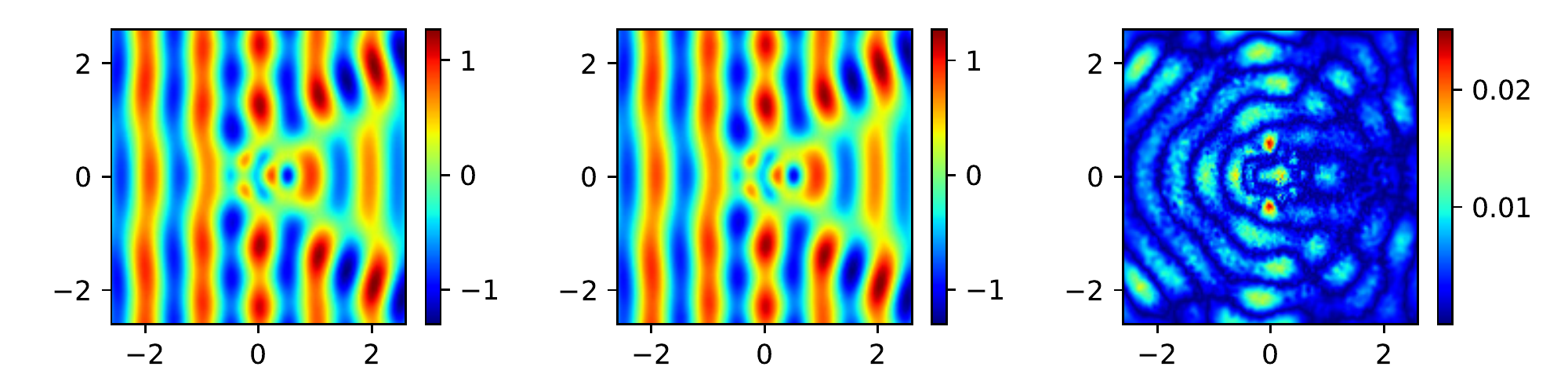}
\end{minipage}
\caption{Scattering of a plane wave by a dielectric disk.
	 Comparison of the 2-D distribution of the real part of $E_z$
	 between DGTD (left), CAE-CSI (middle) and relative error (right) 
         for the testing parameter instances: 
	 $\varepsilon_1=1.215$ (1st row), $\varepsilon_{2}=2.215$ (2th row), 
         $\varepsilon_3=3.215$ (3th row), $\varepsilon_{4}=4.215$ (4th row).}
\label{fig:xy_fields_ez}
\end{figure}

\begin{figure}[htbp]
\centering
\includegraphics[width=1\textwidth]
                {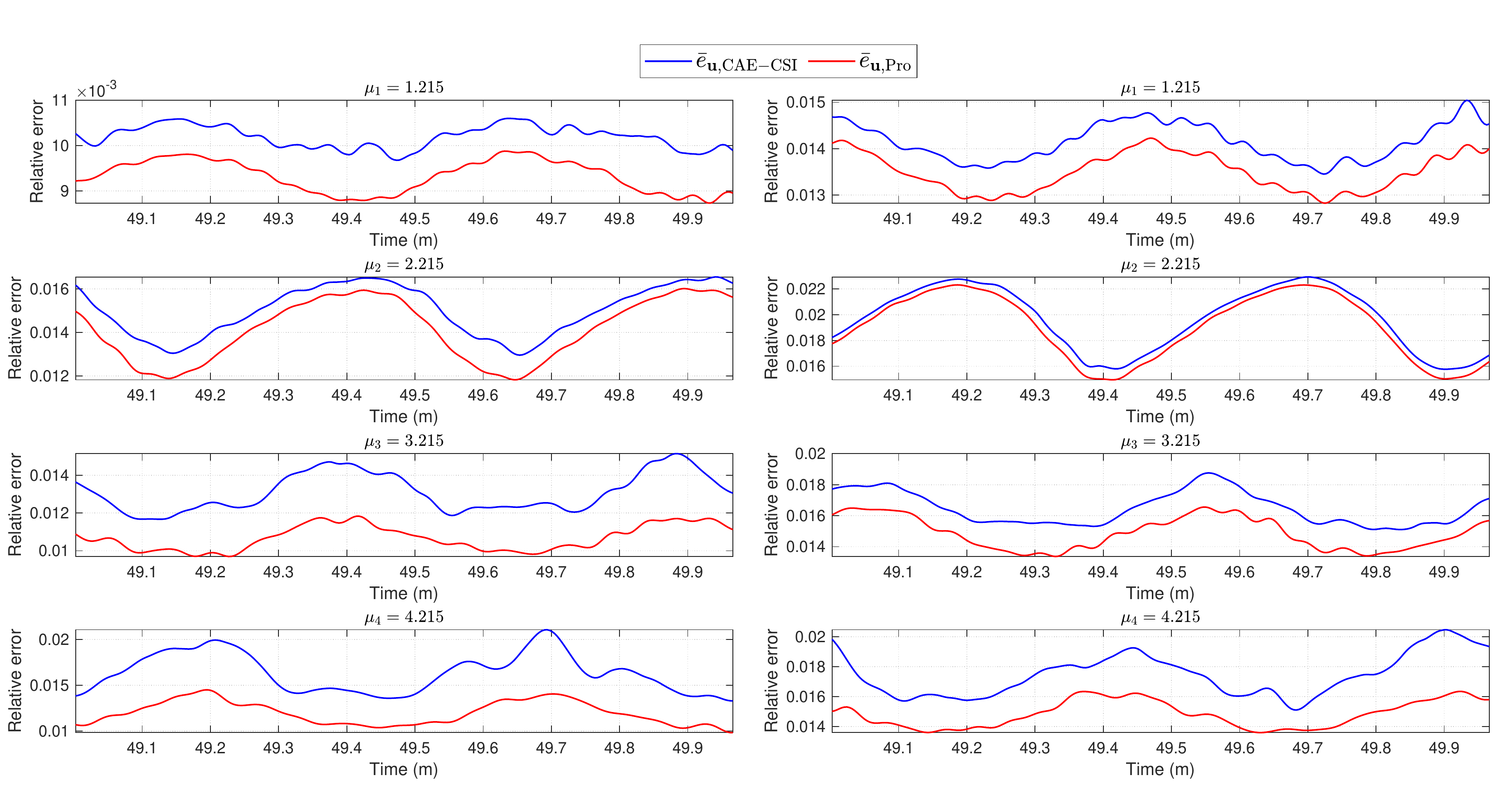}
\caption{Scattering of a plane wave by a dielectric disk. 
	 Comparison of the relative error between DGTD and CAE-CSI for 
         the field  $\mathbf{H}$ (left) and $\mathbf{E}$ (right)
	 with the testing parameter instances 
	 $\varepsilon_1=1.215$ (1st row), $\varepsilon_{2}=2.215$ (2th row), 
         $\varepsilon_3=3.215$ (3th row), $\varepsilon_{4}=4.215$ (4th row).}
\label{fig:time_error}
\end{figure}

To  perform  a  moderate  dimensionality  reduction,  we  utilize  the
two-step  POD method  to  extract  the basis  functions  from the  the
snapshot  matrix.   Fig.~\ref{fig:convergence} shows  the  convergence
histories        of         $\bar{e}_{\mathbf{H},\rm{POD}}$        and
$\bar{e}_{\mathbf{E},\rm{POD}}$   with   the   choice   of   $k$   and
$\mathcal{N}$ in  the two-step  POD method.  We  obtain the  POD basis
matrices $\mathbf{V_u}$ $(\mathbf{u}  \in \{\mathbf{E}, \mathbf{H}\})$
with $k=4$ and $\mathcal{N} = 196$. Then we train the CAE network, and
Fig.~\ref{fig:loss_vs_n_and_kersize}  left shows  the behavior  of the
loss on validation set with respect  to the reduced dimension $n$.  By
increasing   the   dimension   $n$  of   the   low-dimensional   coded
representation  vector  $\tilde{\mathbf{u}}(t,\boldsymbol{\mu})$  from
$3$  to  $20$,   there  is  a  significant  improvement   of  the  CAE
performance. However, when $n>20$, the loss on the validation set does
not  change significantly.   Thus  we set  the  dimensionality of  the
low-representation        coded        vector        to        $n=20$.
Fig.~\ref{fig:loss_vs_n_and_kersize}  right shows  the  impact of  the
size of the convolutional kernels on the loss over the validation set.
Tabs.~\ref{tab:encoder} and \ref{tab:decoder}  summarize attributes of
the convolutional and transposed convolutional layers. As for the training 
of CSI models, we set the SVD truncation tolerance to $\delta = 1\times 10^{-4}$.
\begin{table}[htbp]
\caption{Architecture of convolutional and fully-connected layers 
         in the encoder $f_{E}$.}
\label{tab:encoder} 
\vspace{0.25cm}
\centering
\begin{tabular}{lllllll}
\toprule
 Layer & Input Dim & Output dim  & Kernel size & Number of kernels & 
Stride & Padding \\
\midrule
1      & $[14,14,3]$ & $[14,14,8]$ & $[5,5]$     & 8   & 1  & 2 \\
2      & $[14,14,8]$ & $[8,8,16]$  & $[5,5]$     & 16  & 2  & 3 \\
3      & $[8,8,16]$  & $[4,4,32]$  & $[5,5]$     & 32  & 2  & 2 \\
4      & $[4,4,32]$  & $[2,2,64]$  & $[5,5]$     & 64  & 2  & 2 \\
5      & 256         & 256         &             &     &    &   \\
6      & 256         & 256         &             &     &    &   \\
7      & 256         & $n$         &             &     &    &   \\
\bottomrule
\end{tabular}
\label{encoder}
\end{table}
\begin{table}[htbp]
\caption{Architecture of fully-connected and transposed convolutional 
layers in the decoder $f_{D}$.}
\label{tab:decoder} 
\vspace{0.25cm}
\centering
\begin{tabular}{lllllll}
\toprule
 Layer & Input Dim & Output dim  & Kernel size & Number of kernels & 
Stride & Padding \\
\midrule
1      & $n$          & 256          &             &     &    &   \\
2      & 256          & 256          &             &     &    &   \\
3      & 256          & 256          &             &     &    &   \\
4      & $[2,2,64]$   & $[4,4,64]$   & $[5,5]$     & 64  & 1  & 1 \\
5      & $[4,4,64]$   & $[8,8,32]$   & $[5,5]$     & 32  & 1  & 0 \\
6      & $[8,8,32]$   & $[14,14,16]$ & $[5,5]$     & 16  & 3  & 6 \\
7      & $[14,14,16]$ & $[14,14,3]$  & $[5,5]$     & 8   & 1  & 2 \\
\bottomrule
\end{tabular}
\label{decoder}
\end{table}
In order to assess the performance of CAE-CSI model, the reduced-order
solutions  are  compared  with the  corresponding  DGTD  high-fidelity
solutions on the test parameter set
$\mathcal{P}_{h}^{te} = 
 \{\boldsymbol{\mu}_1, \boldsymbol{\mu}_2, 
               \boldsymbol{\mu}_3, \boldsymbol{\mu}_4\}$ 
where $\boldsymbol{\mu}_1 = \varepsilon_1 = 1.215, 
       \boldsymbol{\mu}_2 = \varepsilon_2 = 2.215, 
       \boldsymbol{\mu}_3 = \varepsilon_3 = 3.215$ and 
$\boldsymbol{\mu}_4 = \varepsilon_4 = 4.215$.
Firstly,   Fig.~\ref{fig:reduced_order_matrix}   shows   the   exact
reduced-order matrices and approximate reduced-order matrices based on
CSI. Over the Fourier domain during the last oscillation period of the
incident  wave, we  display  in Fig.~\ref{fig:x_fields}  the 1-D  $x$-wise
evolution along $y=0$ of the real part of $H_{y}$ and $E_{z}$, as well
as   their  2-D   distribution   in  Fig.~\ref{fig:xy_fields_hy}   and
Fig.~\ref{fig:xy_fields_ez}.   The  time  evolution  of  the  relative
projection   error   $e_{\mathbf{u},\rm{Pro}}$   and   CAE-CSI   error
$e_{\mathbf{u},\rm{CAE-CSI}}               (\mathbf{u}\in\{\mathbf{E},
\mathbf{H}\})$ for  the test parameter  instances $\boldsymbol{\mu}_1,
\boldsymbol{\mu}_2,  \boldsymbol{\mu}_3$ and  $\boldsymbol{\mu}_4$ are
shown  in  Fig.~\ref{fig:time_error}.    The  average  relative  error
$\bar{e}_{\mathbf{u},\rm{Pro}}$ and $\bar{e}_{\mathbf{u},\rm{CAE-CSI}}
(\mathbf{u}\in\{\mathbf{E},\mathbf{H}\})$ for the four test parameters
are shown  in Tab.~\ref{tab:average_relative_error}.   It can  be seen
that the  reduced-order solutions  and the  DGTD solutions  match each
other very well.  Secondly, Tab.~\ref{tab:computational_time} presents
the time performance  comparison of DGTD and CAE-CSI,  where we record
the average online test time of DGTD solver and ROM for the above four
test cases, as well as the training  time of the CAE-CSI model.  It is
worthwhile to take  a long time to build a  surrogate model, since the
online test  time of the ROM  is greatly shorten compared  to the DGTD
solver, achieving a  speed-up of 3387. Note that the  online test time
of  CAE-CSI   is  shortened   compared  to  previous   work  (POD-CSI)
\cite{Li2021}.
\begin{table}[htbp]
\caption{Scattering of a plane wave by a dielectric disk. 
         Average relative error on the test set.}
\label{tab:average_relative_error}
\vspace{0.25cm}
\centering
\begin{tabular}{cccc}
\toprule
 $\bar{e}_{\mathbf{H},\rm{Pro}}$ & $\bar{e}_{\mathbf{H},\rm{CAE-CSI}}$ & 
 $\bar{e}_{\mathbf{E},\rm{Pro}}$ & $\bar{e}_{\mathbf{E},\rm{CAE-CSI}}$ \\
\midrule
 1.15\%  &  1.37\%  & 1.56\%  &  1.69\% \\
\bottomrule
\end{tabular}
\end{table}
\begin{table}[htbp]
\caption{Scattering of a plane wave by a dielectric disk. 
         Comparison between the CAE-CSI (offline and online) and DGTD methods 
         in terms of CPU time. The unit of time is second.}
\label{tab:computational_time}
\vspace{0.25cm}
\centering
\begin{tabular}{cccccc}
\toprule
          & Offline      &         & Online  &         & \\
\midrule
Snapshots & Two-step POD & CAE-CSI  & CAE-CSI & POD-CSI & DGTD \\
  $2.088 \times 10^{4}$ & $6.0530$  & $1.9042 \times 10^{3}$ & $0.0761$ &  
  $0.2561$ & $2.5778 \times 10^{2}$ \\
\bottomrule
\end{tabular}
\end{table}
\subsection{Scattering of a plane wave by a multi-layer heterogeneous medium}

\begin{figure}[htbp]
\centering
\includegraphics[width=0.6\textwidth]{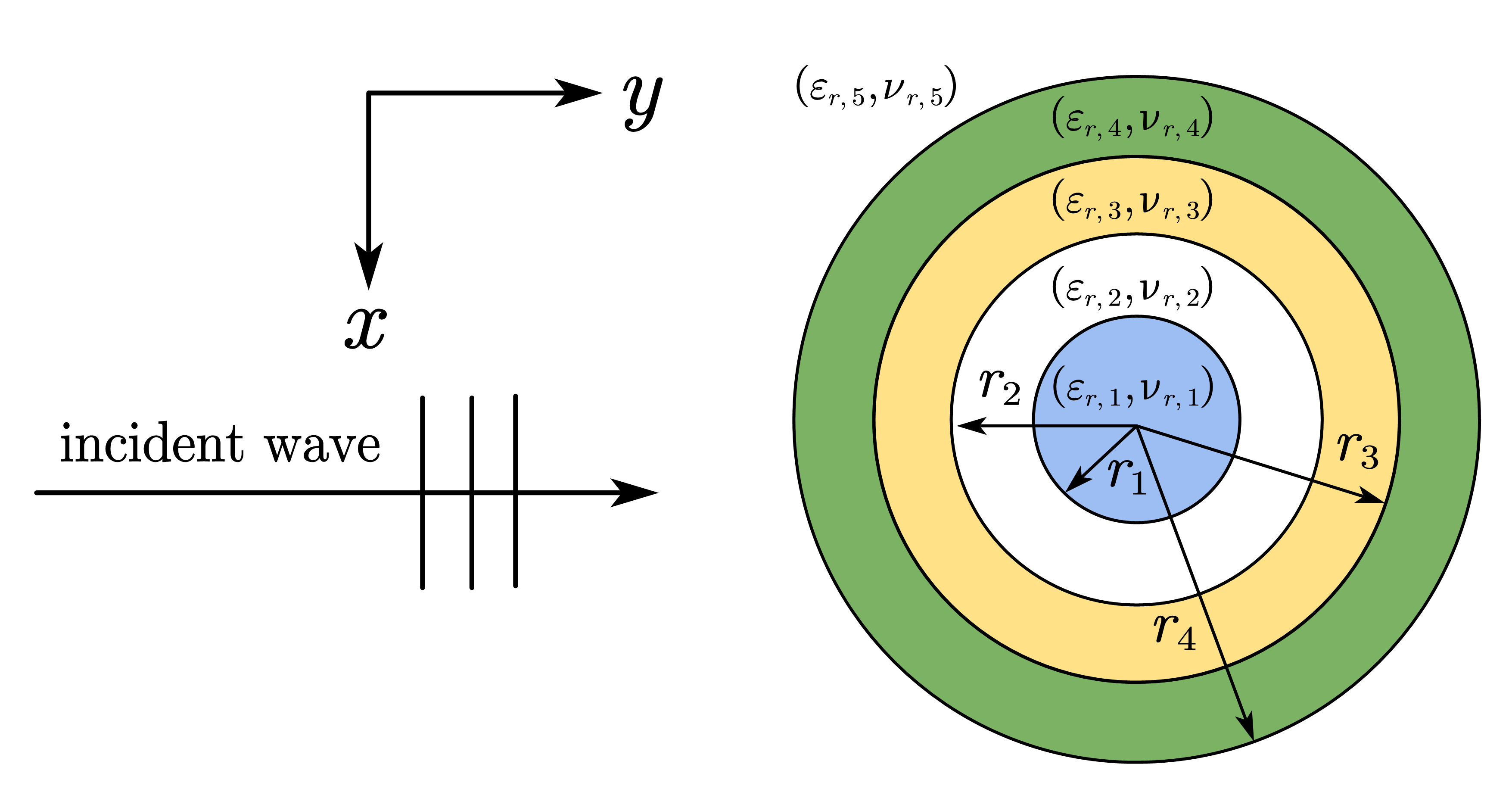}
\caption{Scattering of plane wave by a multi-layer heterogeneous medium. 
         Geometry of the multi-layer heterogeneous medium.}
\label{fig:mutli_layer}
\end{figure}

\begin{figure}[htbp]
\centering
\centering
\includegraphics[scale=0.575]{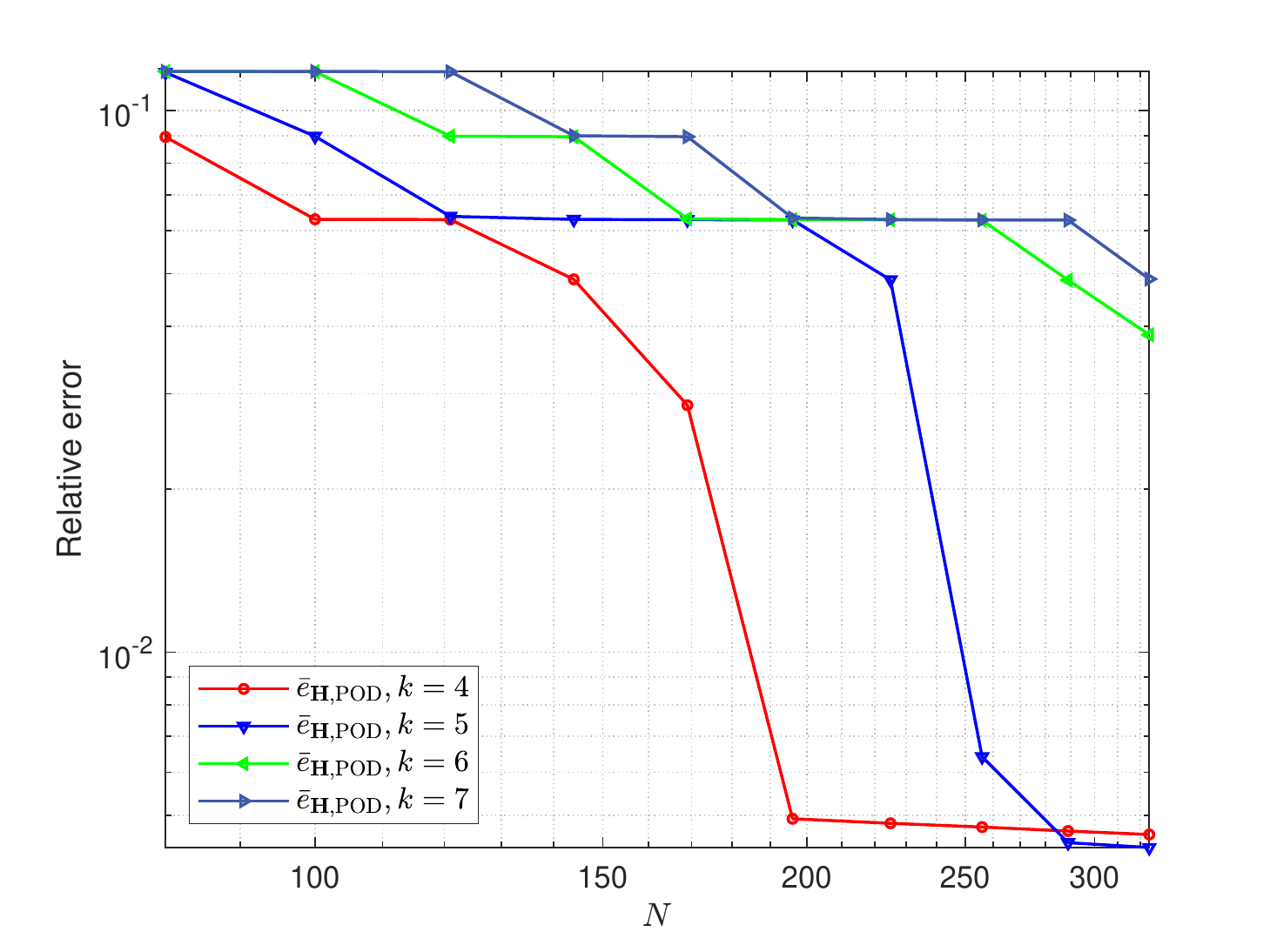}
\hspace{-1.0cm}
\includegraphics[scale=0.575]{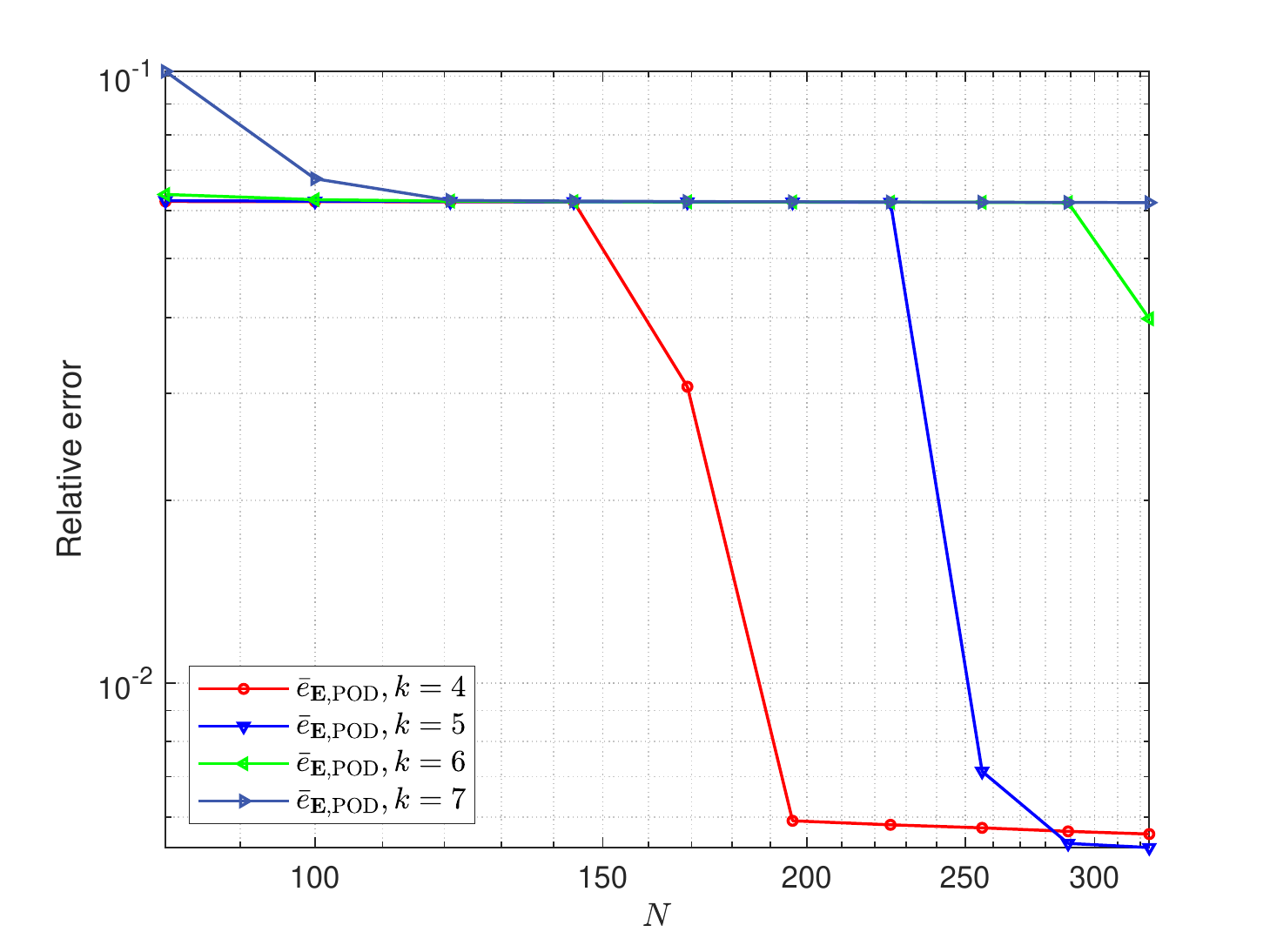}
\caption{Scattering of plane wave by a multi-layer heterogeneous medium. 
         Convergence histories of $\bar{e}_{\mathbf{H},\rm{POD}}$ (left)  
         and $\bar{e}_{\mathbf{E},\rm{POD}}$ (right) with the choice of $k$ 
         and $\mathcal{N}$.}
\label{fig:convergence_multi}
\end{figure}

\begin{figure}[htbp]
\centering
\begin{minipage}[t]{0.49\linewidth}
\centering
\includegraphics[scale=0.45]{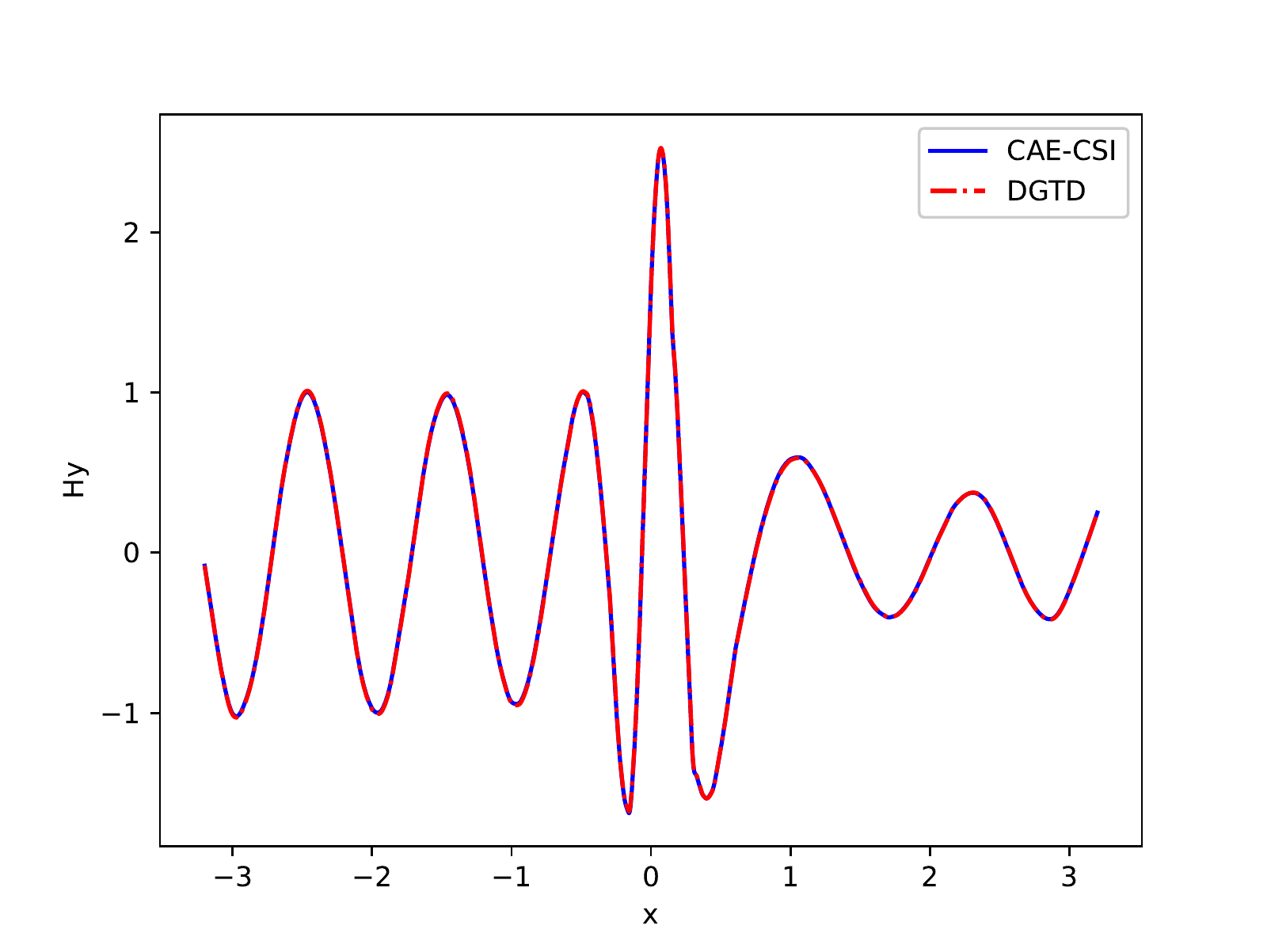}
\end{minipage}
\begin{minipage}[t]{0.49\linewidth}
\centering
\includegraphics[scale=0.45]{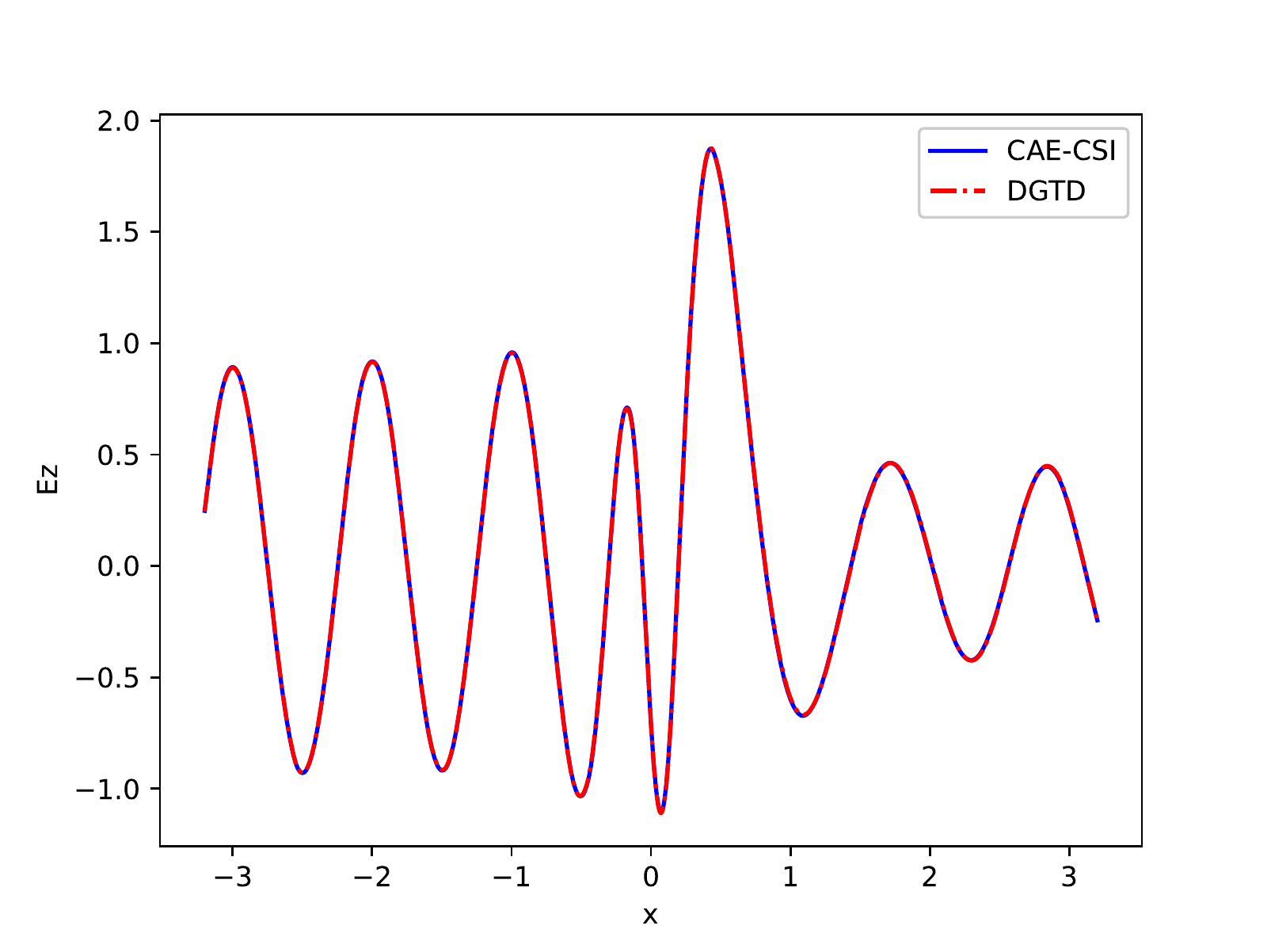}
\end{minipage} \\
\begin{minipage}[t]{0.49\linewidth}
\centering
\includegraphics[scale=0.45]{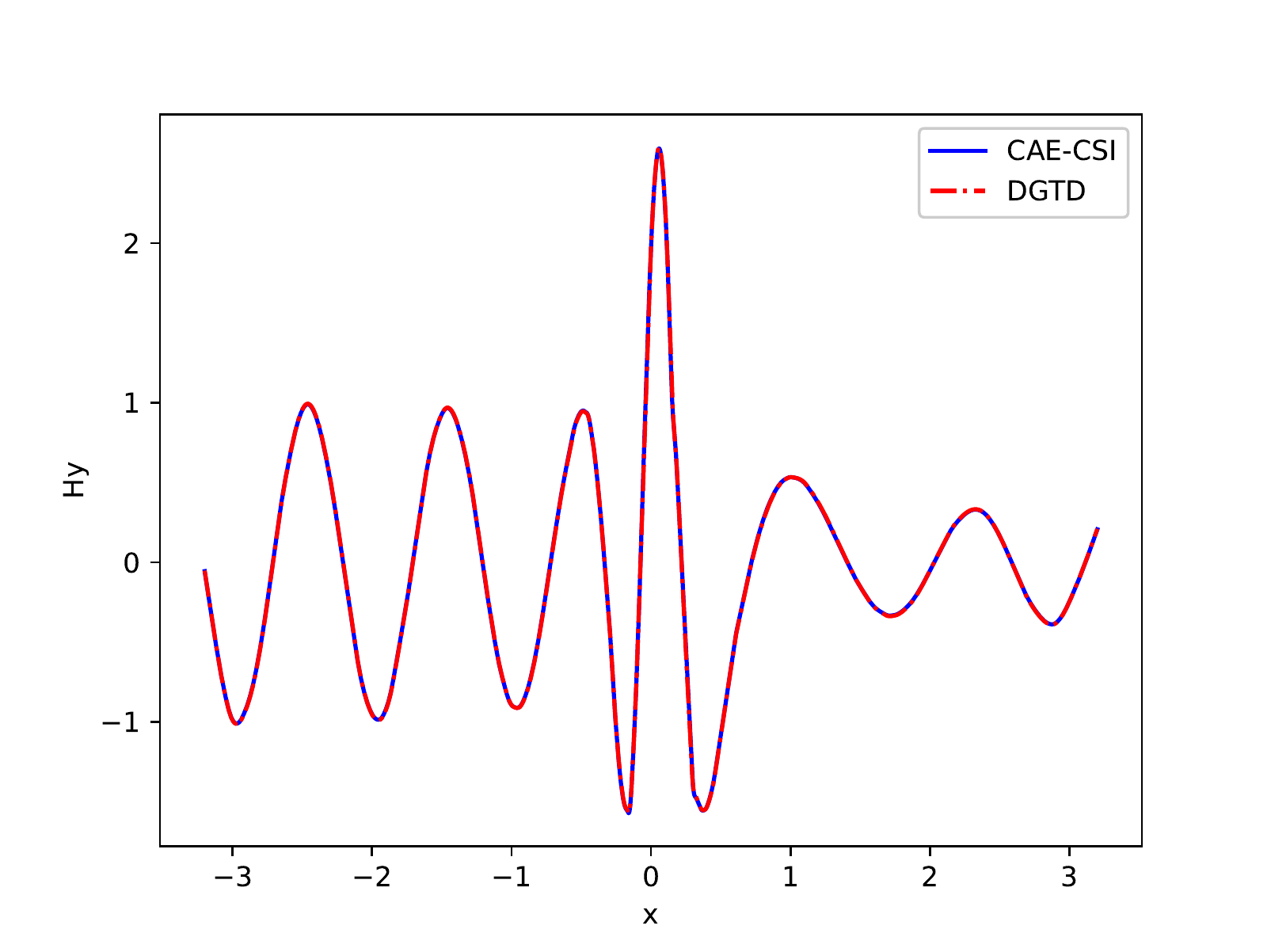}
\end{minipage}
\begin{minipage}[t]{0.49\linewidth}
\centering
\includegraphics[scale=0.45]{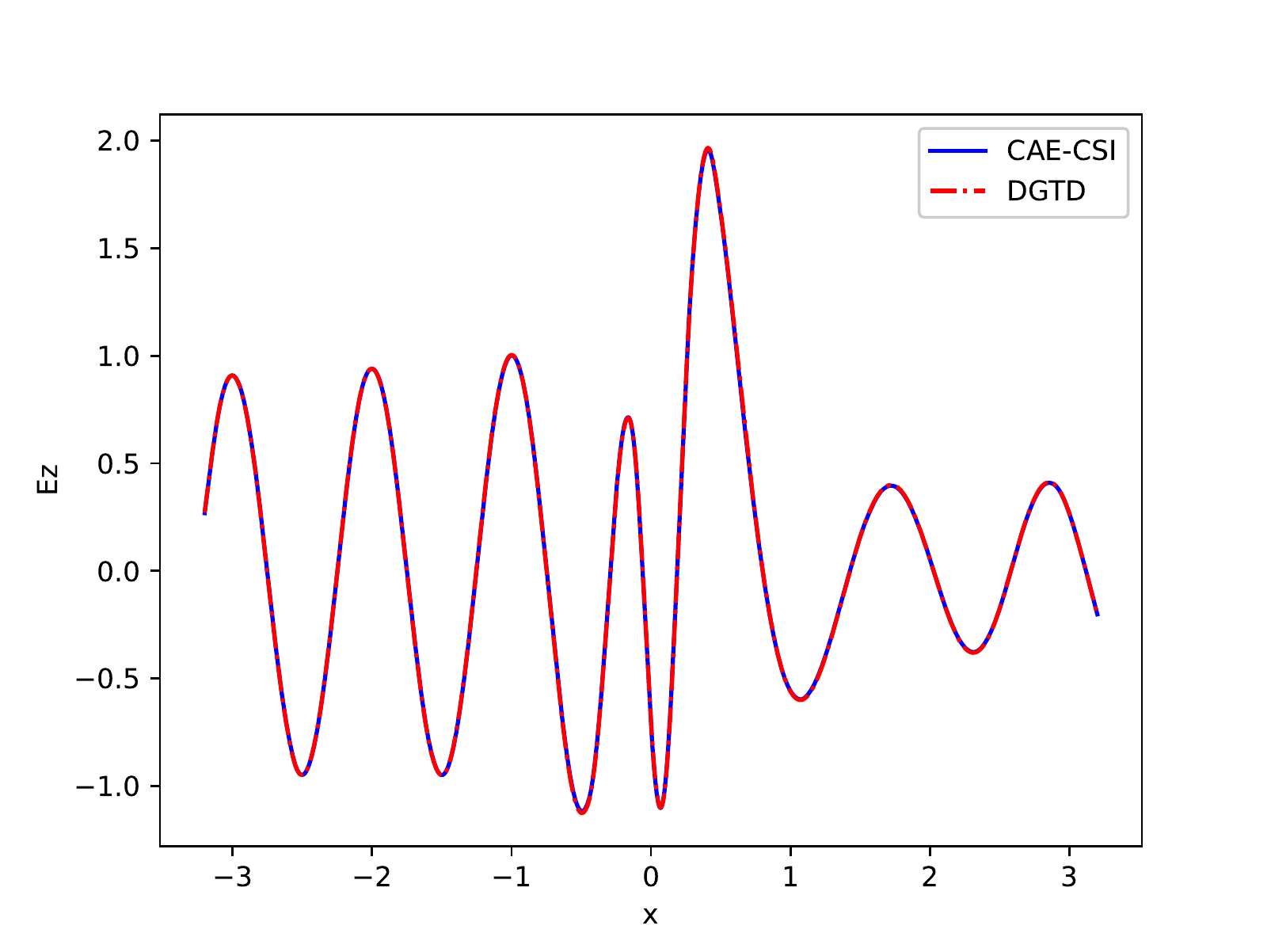}
\end{minipage}
\begin{minipage}[t]{0.49\linewidth}
\centering
\includegraphics[scale=0.45]{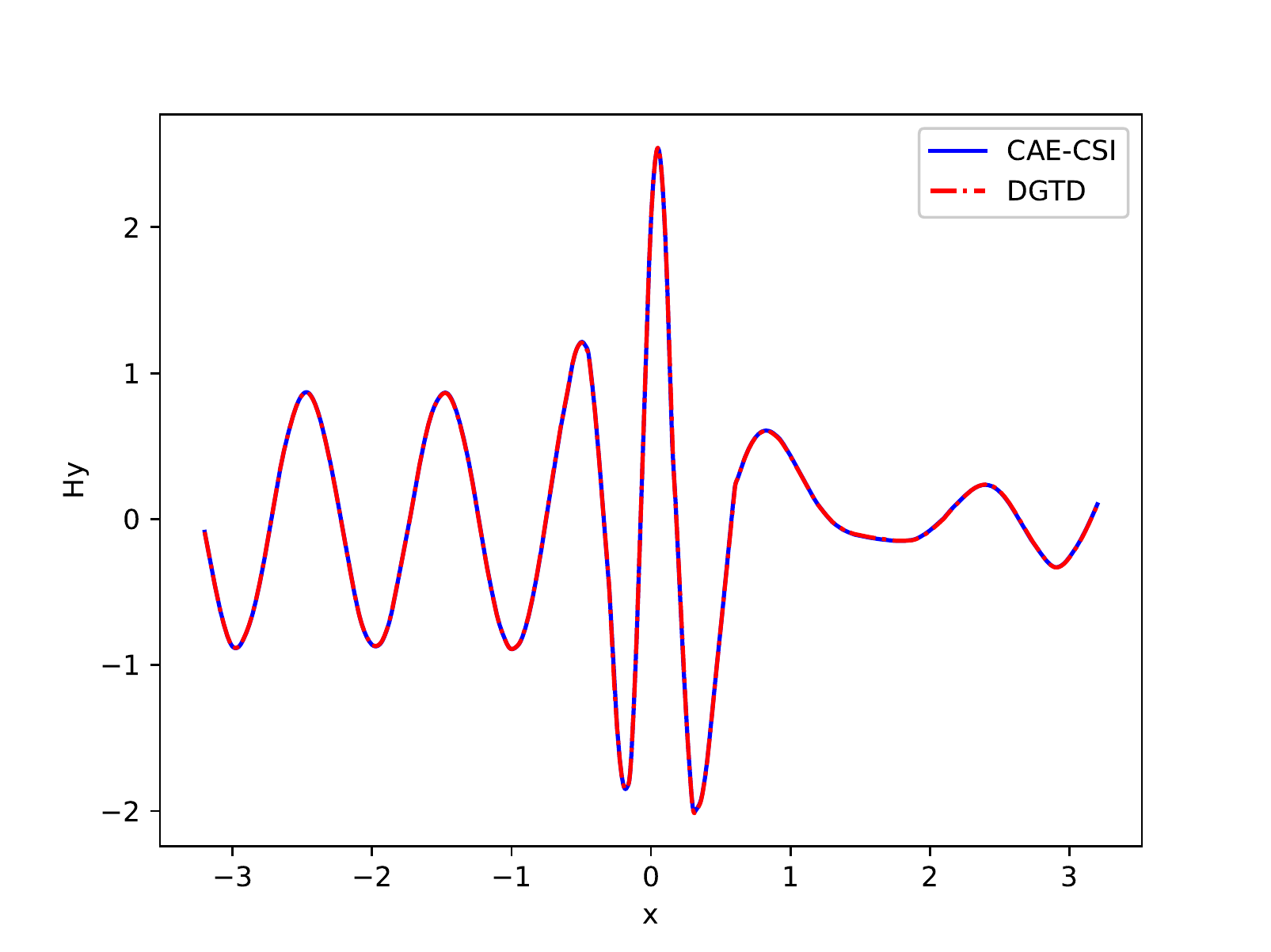}
\end{minipage}
\begin{minipage}[t]{0.49\linewidth}
\centering
\includegraphics[scale=0.45]{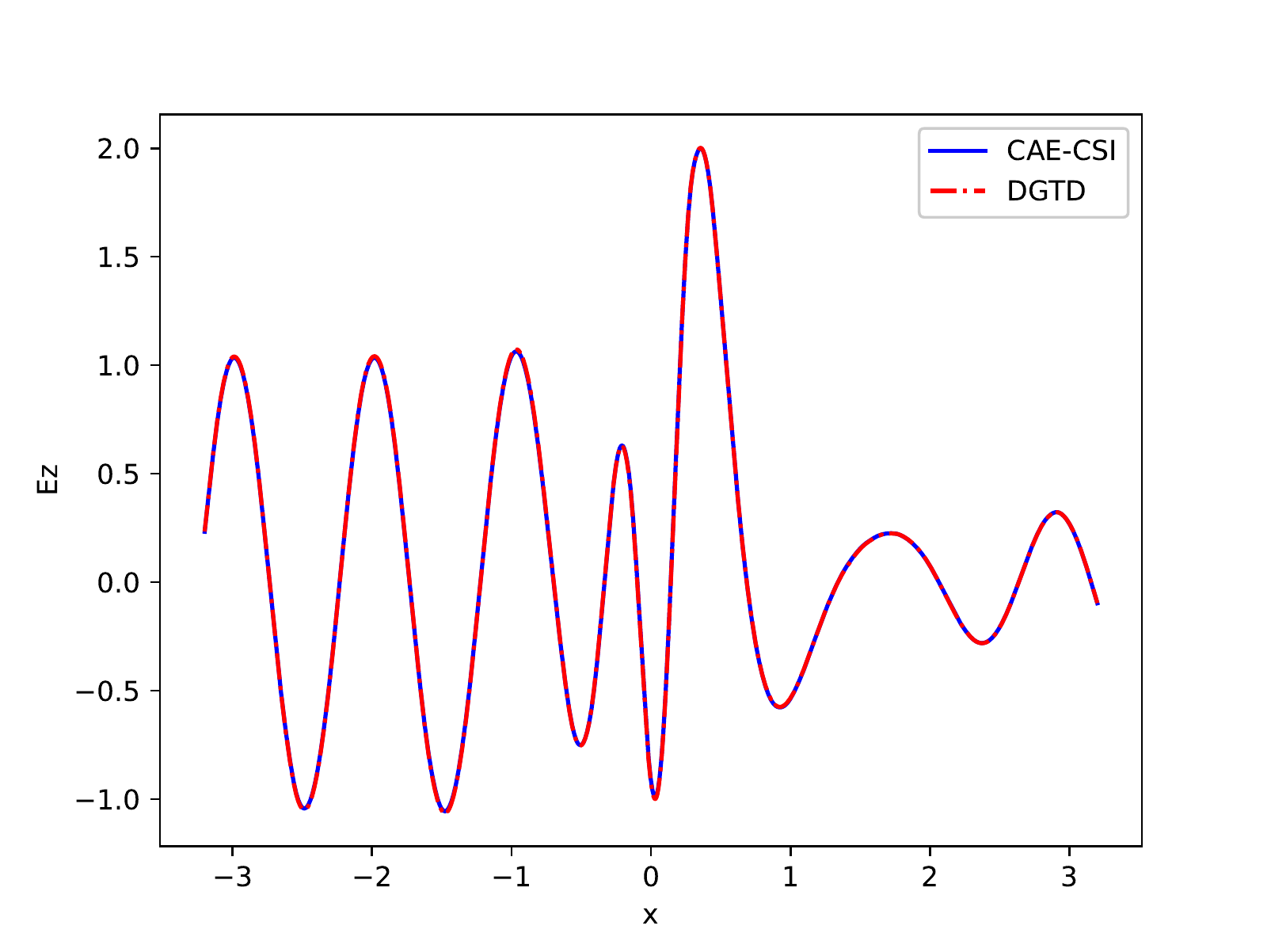}
\end{minipage}
\caption{Scattering of plane wave by a multi-layer heterogeneous medium. 
         Comparison of the 1-D $x$-wise distribution along $y=0$ of the 
         real part of $H_y$ (left) and $E_z$ (right) for the
	 testing parameter instances: 
         $\boldsymbol{\mu}_1=\{(5.1, 3.4, 2.1, 1.4)\}$ (1st row),  
         $\boldsymbol{\mu}_2 = \{(5.4, 3.4, 2.3, 1.3)\}$ (2th row) and 
         $\boldsymbol{\mu}_3 = \{(5.5, 3.7, 2.4, 1.7)\}$ (3th row).}
\label{fig:x_fields_multi}
\end{figure}

\begin{figure}[htbp]
\centering
\begin{minipage}[t]{0.99\textwidth}
\centering
\includegraphics[scale=0.7]{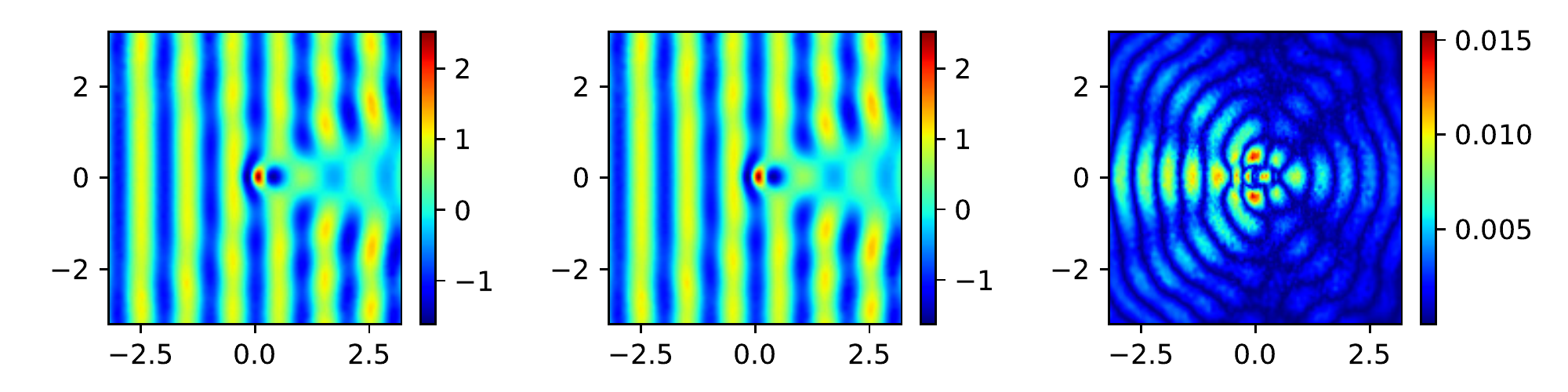}
\end{minipage} \\
\begin{minipage}[t]{0.99\textwidth}
\centering
\includegraphics[scale=0.7]{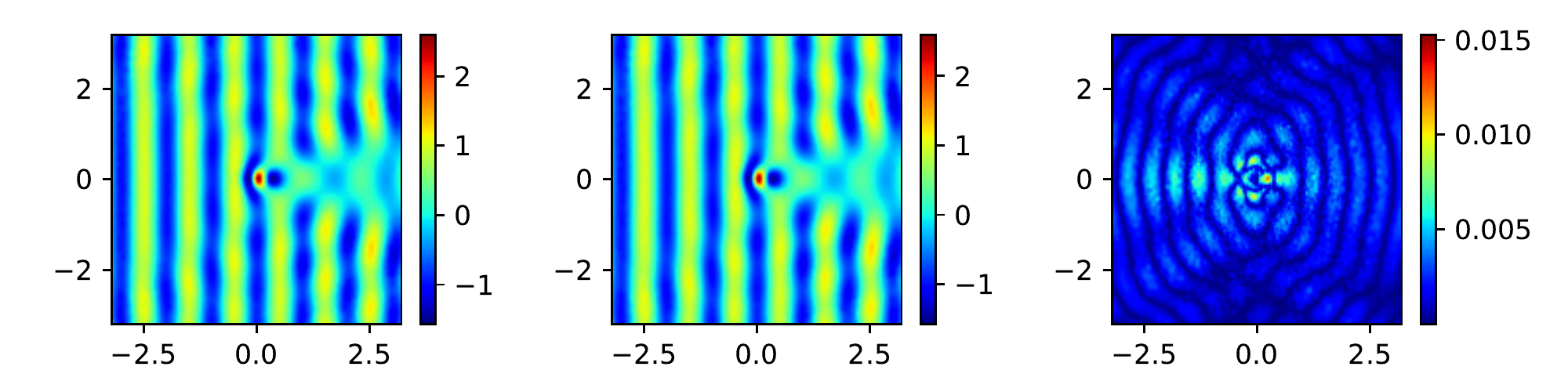}
\end{minipage} \\
\begin{minipage}[t]{0.99\textwidth}
\centering
\includegraphics[scale=0.7]{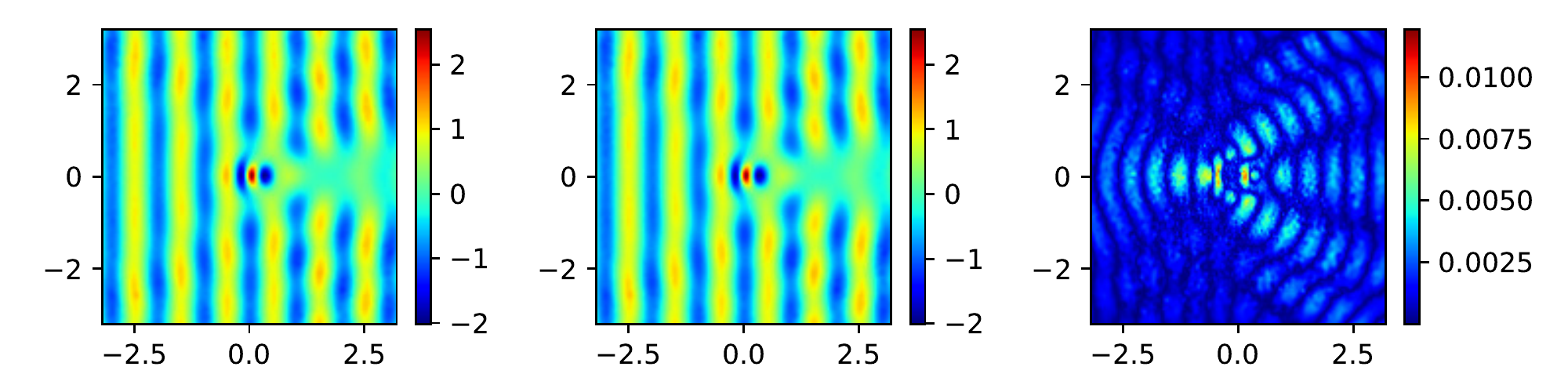}
\end{minipage}
\caption{Scattering of plane wave by a multi-layer heterogeneous medium.
 	 Comparison of the 2-D distribution of the real part of $H_{y}$
 	 between DGTD (left), CAE-CSI (middle) and relative error (right) 
         for the testing parameter instances: 
         $\boldsymbol{\mu}_1=\{(5.1, 3.4, 2.1, 1.4)\}$ (1st row),  
         $\boldsymbol{\mu}_2 = \{(5.4, 3.4, 2.3, 1.3)\}$ (2th row) and 
         $\boldsymbol{\mu}_3 = \{(5.5, 3.7, 2.4, 1.7)\}$ (3th row).}
\label{fig:xy_fields_hy_multi}
\end{figure}

\begin{figure}[htbp]
\centering
\begin{minipage}[t]{1\textwidth}
\centering
\includegraphics[scale=0.7]{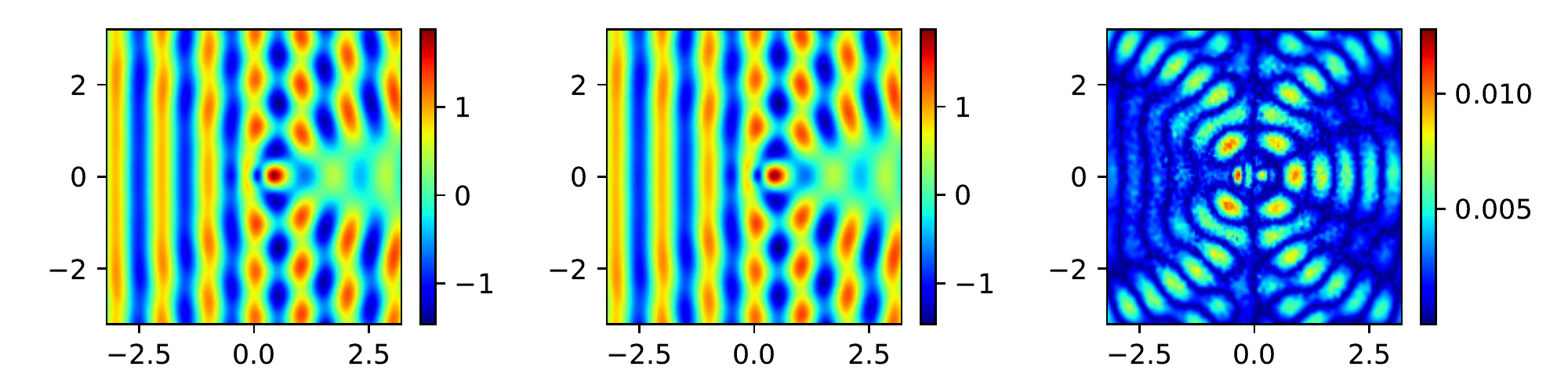}
\end{minipage} \\
\begin{minipage}[t]{1\textwidth}
\centering
\includegraphics[scale=0.7]{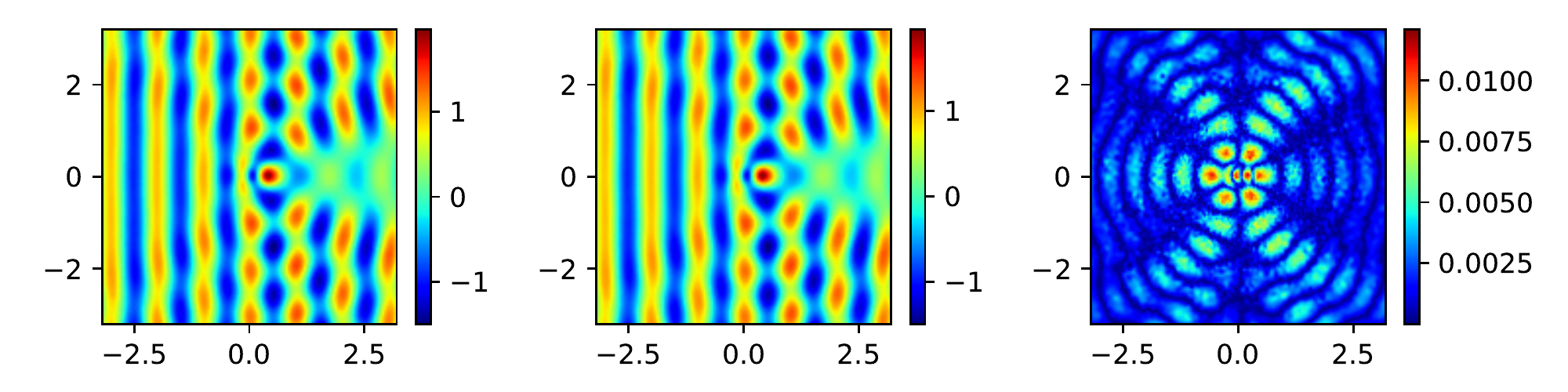}
\end{minipage} \\
\begin{minipage}[t]{1\textwidth}
\centering
\includegraphics[scale=0.7]{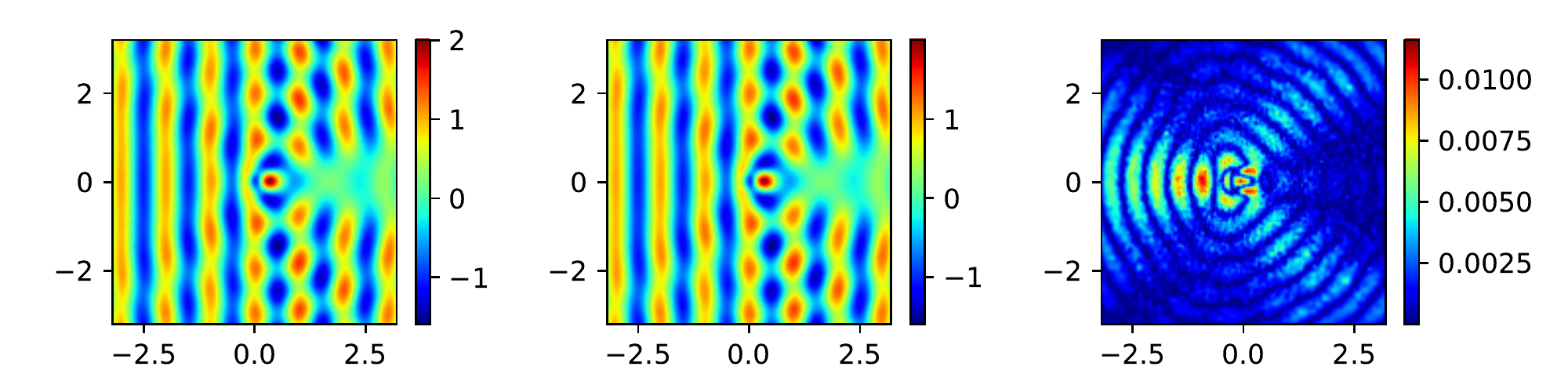}
\end{minipage}
\caption{Scattering of plane wave by a multi-layer heterogeneous medium:
         comparison of the 2-D distribution of the real part of $E_{z}$
         between DGTD (left), CAE-CSI (middle) and relative error (right) 
         for the testing parameter instances: 
         $\boldsymbol{\mu}_1=\{(5.1, 3.4, 2.1, 1.4)\}$ (1st row),  
         $\boldsymbol{\mu}_2 = \{(5.4, 3.4, 2.3, 1.3)\}$ (2th row) and 
         $\boldsymbol{\mu}_3 = \{(5.5, 3.7, 2.4, 1.7)\}$ (3th row).}
\label{fig:xy_fields_ez_multi}
\end{figure}

\begin{figure}[htbp]
\centering
\includegraphics[width=1\textwidth]{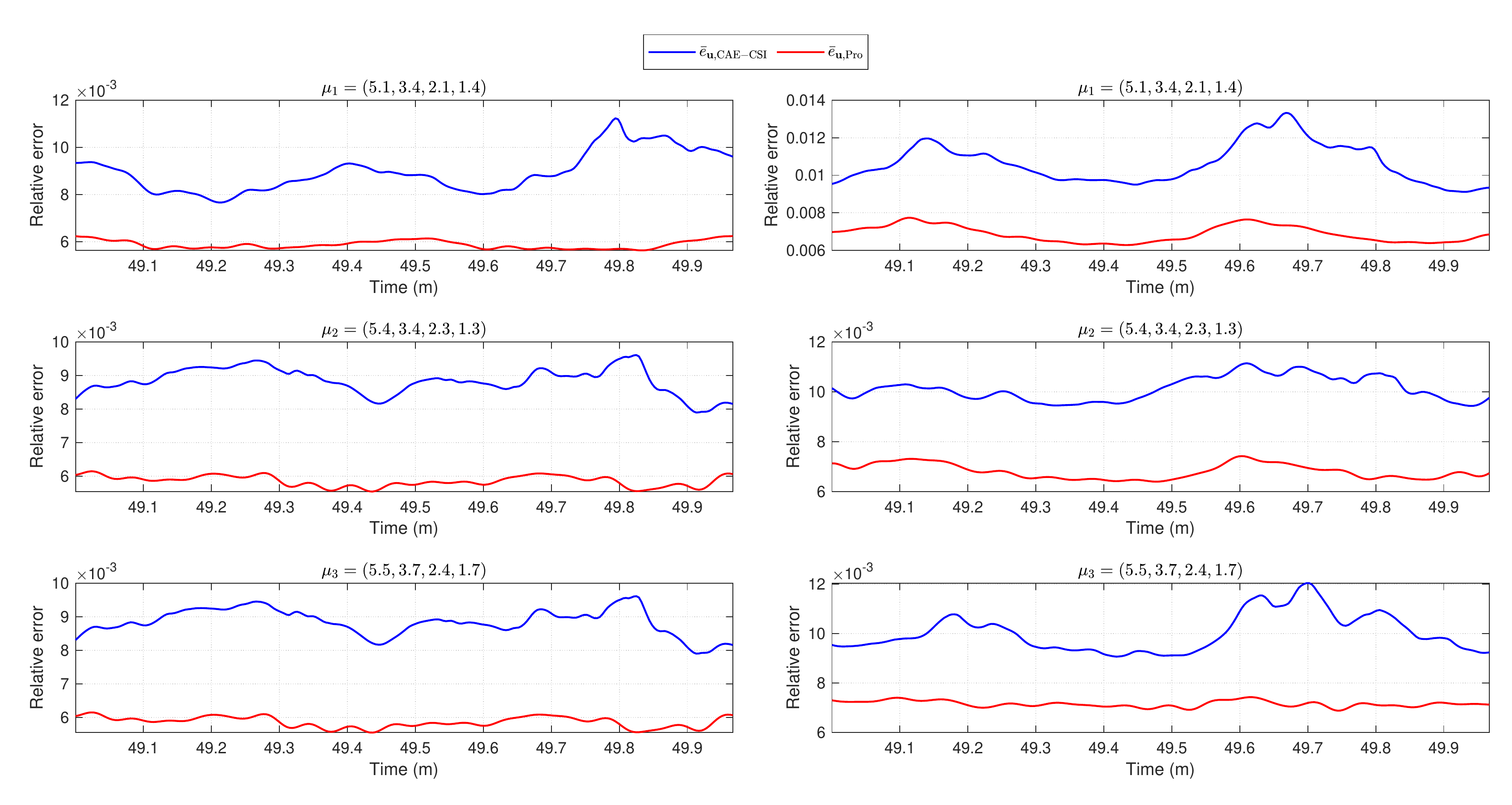}
\caption{Scattering of plane wave by a multi-layer heterogeneous medium. 
         Comparison of the relative projection error $e_{\mathbf{u},\rm{Pro}}$, 
         the CAE-CSI error $e_{\mathbf{u},\rm{CAE-CSI}}$ for 
         $\mathbf{H}$ (left) and $\mathbf{E}$ (right) for the testing parameter
         instances: 
         $\boldsymbol{\mu}_1=\{(5.1, 3.4, 2.1, 1.4)\}$ (1st row),  
         $\boldsymbol{\mu}_2 = \{(5.4, 3.4, 2.3, 1.3)\}$ (2th row) and 
         $\boldsymbol{\mu}_3 = \{(5.5, 3.7, 2.4, 1.7)\}$ (3th row).}
\label{fig:time_error_multi}
\end{figure}

In  this subsection,  we consider  a multi-layer  heterogeneous medium
which  is  illuminated   by  an  incident  plane  wave   as  shown  in
Fig.~\ref{fig:mutli_layer}.  The  computational domain  is artificially
bounded  by  a  square  $\Omega=[-3.2  \mathrm{~m},  3.2  \mathrm{~m}]
\times[-3.2  \mathrm{~m},  3.2  \mathrm{~m}]$   where  we  impose  the
Silver-Müller            ABC            boundary            condition.
Tab.~\ref{tab:range_of_parameters}  summarizes  the  distribution  and
range  of  material  parameters  considered  in  this  study.The  mesh
consists of  3256 nodes and  6206 elements, resulting in  $\mathcal{N}_{h}=37236$ DOFs  of the FOM.
\begin{table}
\caption{Scattering of plane wave by a multi-layer heterogeneous medium. 
         Distribution and range of material parameters.}
\centering
\vspace{0.25cm}
\begin{tabular}{llll}
\toprule
Layer $i$ & $\mathcal{P}^i$ & $\nu_{r,i}$ & $r_i$ \\
\midrule
1         & $\varepsilon_{r,1} \in [5.0, 5.6]$   & 1.0 & 0.15  \\
2         & $\varepsilon_{r,2} \in [3.25, 3.75]$ & 1.0 & 0.3   \\
3         & $\varepsilon_{r,3} \in [2.0, 2.5]$   & 1.0 & 0.45  \\
4         & $\varepsilon_{r,4} \in [1.25, 1.75]$ & 1.0 & 0.6   \\
\bottomrule
\end{tabular}
\label{tab:range_of_parameters}
\end{table}
The DGTD solver  is used to generate the snapshots,  and the numerical
simulation time is set to 50  periods of the incident wave oscillation
with time  step $\Delta t=0.0038$  s for each  4-dimensional parameter
$\boldsymbol{\mu} \in \mathcal{P}$, where $\mathcal{P} = \mathcal{P}^{1}
\times \mathcal{P}^{2}  \times \mathcal{P}^{3} \times \mathcal{P}^{4}
\subset \mathbb{R}^{4}$  with $\varepsilon_{r,i} \in \mathcal{P}^{i}$
$(i=1, 2, 3, 4)$.  We  adopt a  grid  sampling of  tensor product  with
$\mathcal{N}_{\Delta p} = 3$ uniform points for each parameter to form a
training  parameter   samples  $\mathcal{P}_{h}^{tr}$,   resulting  in
$\mathcal{N}_{p} = 81$ points.   Each  choice  of the  parameter  is
sampled  for $\mathcal{N}_{t}  = 253$  snapshots in  time at  the last
period,  i.e.,  $\mathcal{T}_{h}^{t  r} =  \{49.0002, 49.0041,  \cdots,
49.9630, 49.9669\}$.    Fig.~\ref{fig:convergence_multi} shows the
convergence    histories   of    $\bar{e}_{\mathbf{H}, \rm{POD}}$  and
$\bar{e}_{\mathbf{E},\rm{POD}}$  with different  truncation parameter
$k$  and size parameter $\mathcal{N}$  in the  two-step  POD method.   The POD  basis
matrices        $\mathbf{V}_{\mathbf{u}} \  (\mathbf{u} \in
\{\mathbf{E} , \mathbf{H}\})$ are  generated by the two-step  POD method
with $k = 4$ and $\mathcal{N} = 196$.  The architecture of the CAE network
is the same as the one used in the first experiment.  After performing
a SVD to  all reduced-order matrices with a  truncation tolerance $\delta = 1\times 10^{-4}$, the CSI models  are built as the
combination of time- and parameter-modes. 
To evaluate the CAE-CSI method,  we compare the reduced-order solution
with the high-fidelity solution generated by the DGTD method on a test
parameter    set    $\mathcal{P}_{h}^{te} = \{\boldsymbol{\mu}_1,
\boldsymbol{\mu}_2, \boldsymbol{\mu}_3\}$ with  $\boldsymbol{\mu}_1 =
\{(5.1,  3.4, 2.1,  1.4)\}$, $\boldsymbol{\mu}_2 = \{(5.4, 3.4, 2.3,
1.3)\}$    and    $\boldsymbol{\mu}_3 = \{(5.5, 3.7,  2.4,
1.7)\}$.   Fig.~\ref{fig:x_fields_multi}   shows  the   1-D   $x$-wise
distribution along  $y=0$ of the real  part of $H_y$ and  $E_z$ in the
Fourier    domain    on    the     last    period    of    simulation.
Fig.~\ref{fig:xy_fields_hy_multi} and
Fig.~\ref{fig:xy_fields_ez_multi} display the  2-D distribution of the
real part  of $H_y$  and $E_z$.   The time  evolution of  the relative
projection   error   $e_{\mathbf{u},\rm{Pro}}$   and   CAE-CSI   error
$e_{\mathbf{u}, \rm{CAE-CSI}} (\mathbf{u} \in
\{\mathbf{E}, \mathbf{H}\})$ for 3 test parameters $\boldsymbol{\mu}_1,
\boldsymbol{\mu}_2$   and   $\boldsymbol{\mu}_3$  are   displayed   in
Fig.~\ref{fig:time_error_multi}.     The    average    projection    error
$\bar{e}_{\mathbf{u}, \rm{Pro}}$   and   the  average   CAE-CSI   error
$\bar{e}_{\mathbf{u}, \rm{CAE-CSI}} (\mathbf{u}\in\{\mathbf{E}, \mathbf{H}\})$  
for the  3 test  parameters are listed in 
Tab.~\ref{fig:average_relative_error2}.
\begin{table}[htbp]
\centering
\caption{Scattering of plane wave by a multi-layer heterogeneous medium.  
	 Average relative error on the test set.}
\label{fig:average_relative_error2}
\vspace{0.25cm}
\begin{tabular}{cccc}
\toprule
 $\bar{e}_{\mathbf{H}, \rm{Pro}}$ & $\bar{e}_{\mathbf{H}, \rm{CAE-CSI}}$ & 
 $\bar{e}_{\mathbf{E}, \rm{Pro}}$ & $\bar{e}_{\mathbf{E}, \rm{CAE-CSI}}$ \\
\midrule
 0.60\%  &  0.89\%  & 0.69\% &  1.03\% \\
\bottomrule
\end{tabular}
\end{table}
Performance   results   of  the   CAE-CSI   and   DGTD  methods   with
$\mathbb{P}_{2}$       approximation      are       summarized      in
Tab.~\ref{tab:computational_time_multi}.   The CPU  time  of the  DGTD
method is $3.0477 \times 10^{2}$ s  and the online cost of the CAE-CSI
is  0.2923 s,  corresponding to  a speed-up  of 1042.   Note that  the
online test  time of CAE-CSI  is again shortened compared  to previous
work   (POD-CSI)   \cite{Li2021},   which   again   demonstrates   the
significantly enhanced efficiency of the CAE-CSI method.
\begin{table}[htbp]
\centering
\caption{Scattering of plane wave by a multi-layer heterogeneous medium. 
         Comparison between the CAE-CSI (offline and online) and DGTD  
         methods in terms of CPU time.
	 The unit of time is second.}
\label{tab:computational_time_multi}
\vspace{0.25cm}
\begin{tabular}{cccccc}
\toprule
 & Offline & & Online & & \\
	\midrule
Snapshots & Two-step POD & CAE-CSI & CAE-CSI & POD-CSI & DGTD \\
$2.4686 \times 10^{4}$ & 6.5544 & $2.0432 \times 10^{3}$ & $0.2923$ & 
$1.4856$ & $3.0477\times 10^{2}$  \\
\bottomrule
\end{tabular}
\end{table}

\section{Conclusion}
\label{sec:concl}

A data-driven RB method based on POD,  CAE and CSI is
proposed  to  accelerate  the solution  of  parameterized  time-domain
Maxwell  equations.    The  two-step  POD  method   performs  a  prior
dimensionality  reduction  of  the  snapshots, then  the  CAE  network
provides   a  low   dimensional  representation   of  the   projection
coefficients through its  encoder, as well as the  inverse map through
its decoder.  CSI-based models  are trained  to approximate  a mapping
from the  time/parameters to  the low dimensional  representations and
the  decoder is  used to  reconstruct  the system  solutions to  their
original   dimension.    By   combining   CSI  with   a   decoder,   a
`simulation-free' approach is established to obtain a RB solution at a
very low  cost. The results  of the numerical  experiments demonstrate
the powerful dimensionality  reduction and reconstruction capabilities
of the CAE-CSI method, and a highly accurate and cheap surrogate model
for  systems   described  by  PDEs  is   developed.   Future  research
directions include more realistic 3-D simulations and the reduction of
parameterized geometry.
\section*{Acknowledgments} The last author is supported by 
NSFC (Grant No. 12101511).

\bibliography{References}
\end{document}